\documentclass[a4paper,12pt]{article}

\usepackage[T1]{fontenc}
\usepackage[utf8]{inputenc}
\usepackage[english]{babel}
\usepackage[doi=false,isbn=false,url=false]{biblatex}
\usepackage[]{csquotes}
\usepackage{bbm}
\usepackage{float}
\usepackage{amsmath}
\usepackage{amsthm}
\usepackage{amsfonts}
\usepackage{amssymb}
\usepackage{graphicx}
\usepackage{icomma}
\usepackage{latexsym}
\usepackage{textcomp}
\usepackage[all]{xy}
\usepackage{cancel}
\usepackage{dsfont}
\usepackage{caption}
\usepackage{listings}
\usepackage{lmodern}
\usepackage{xcolor}
\usepackage{moreverb}
\usepackage{mathrsfs}
\usepackage{stmaryrd}
\usepackage{thmtools}
\usepackage{hyperref}
\usepackage{cleveref}

\newcommand{\pen}[0]{\textnormal{\bfseries pen}}

\newcommand{\suml}[0]{\sum\limits}

\newcommand{\vect}[0]{\textbf{Span}_{\mathbb{R}}}
\newcommand{\mbf}[1]{\mathbf{#1}}
\newcommand{\mscr}[1]{\mathscr{#1}}
\newcommand{\mcal}[1]{\mathcal{#1}}
\newcommand{\mbb}[1]{\mathbb{#1}}

\newcommand{\cauchy}[0]{\text{Cauchy}}

\newcommand{\slfrac}[2]{\left.#1\middle/#2\right.}

\newtheorem{lemme}{Lemma}

\newtheorem{hypothese}{Assumption}
\newtheorem{prop}{Proposition}
\newtheorem{theoreme}{Theorem}
\newtheorem{coro}{Corollary}
\newtheorem{exemple}{Example}

\providecommand{\keywords}[1]
{
  \small	
  \textbf{\textit{Keywords---}} #1
}

\title{Robust Estimation in Finite Mixture Models}
\author{Alexandre Lecestre\thanks{This project has received funding from the European Union's Horizon 2020 research and innovation programme under grant agreement N° 811017}}
\date{}

\begin{document}
\maketitle
\begin{abstract}
We observe a $n$-sample, the distribution of which is assumed to belong, or at least to be close enough, to a given mixture model. We propose an estimator of this distribution that belongs to our model and possesses some robustness properties with respect to a possible misspecification of it. We establish a non-asymptotic deviation bound for the Hellinger distance between the target distribution and its estimator when the model consists of a mixture of densities that belong to VC-subgraph classes. Under suitable assumptions and when the mixture model is well-specified, we derive risk bounds for the parameters of the mixture. Finally, we design a statistical procedure that allows us to select from the data the number of components as well as suitable models for each of the densities that are involved in the mixture. These models are chosen among a collection of candidate ones and we show that our selection rule combined with our estimation strategy result in an estimator which satisfies an oracle-type inequality.
\end{abstract}
\section{Introduction}
Mixture models are a flexible tool for modeling heterogeneous data, e.g. from a population consisting of multiple hidden  homogeneous subpopulations. Finite mixture models are models containing distribution of the form
\begin{equation}
\label{eq:p_w_f}
P_{w,F} = \suml_{k=1}^K w_k F_k,
\end{equation}
where $K\geq 2$, each $F_k$ belongs to a specific class of probability distributions (e.g. normal distributions in the case of Gaussian mixture models) and $w$ belongs to the simplex $\mcal{W}_K = \left\{ w \in [0,1]^K; w_1 + \dots + w_k = 1 \right\}$. For a complete introduction to mixture models and an overview of the different applications we refer to the books of Mclachlan \& Peel \cite{maclachlan} and Frühwirth-Schnatter \cite{fruhwirth}.\par
Assume we dispose of a sample $\mbf{X} := (X_1,\dots,X_n)$ of i.i.d. data, each coordinate following the probability distribution $P^*$. The majority of the statistical methods based on finite mixture models aim to solve one of the following problems: density estimation (estimation of $P^*$), parameter estimation (estimation of $w^*$ and/or $F^*$ assuming $P^*=P_{w^*,F^*}$) and clustering. The monographs of Everitt \& Hand \cite{everitt} or Titterington \emph{et al} \cite{titterington1985statistical} provide a good overview of the different estimation methods that have been developed for mixture models such as maximum likelihood, minimum chi-square, moments method and Bayesian approaches. Although algorithms are numerous, theoretical guarantees are mostly asymptotic and restricted to very specific situations. To our knowledge, only a few non-asymptotic results have been established in the case of density estimation based on Gaussian Mixture Models (GMMs). The approximation and entropy properties of Gaussian mixture sieves have been investigated by Kruijer \emph{et al} \cite{kruijer}, Ghosal \& van der Vaart \cite{ghosal2001} and Genovese \& Wasserman \cite{genovese} where bounds on the convergence rate are given for the MLE and Bayesian estimators. Similarly, Maugis \& Michel \cite{maugis2011} use a penalized version of the MLE to build a Gaussian mixture estimator with non asymptotic adaptive properties proven in \cite{maugis2012}. However, those results rely on relatively strong assumptions and estimators are not proved to be robust to small departures from those assumptions.\\
This paper aims to provide non-asymptotic results in a very general setting. In our framework, the data are assumed to be independent but not necessarily i.i.d. Our mixture model consists of probabilities of the form (\ref{eq:p_w_f}) where the $F_k$ admit densities, called \emph{emission densities}, that belong to classes of function that are VC-subgraph. We investigate the performances of $\rho$-estimators, as defined by Baraud and Birgé \cite{baraudrevisited},  on finite mixture models. Our main result is an exponential deviation inequality  for the risk of the estimator $\hat{P}$, which is measured with an Hellinger-type loss. We get an upper bound on the risk that is the sum of two terms. The first one is an approximation term which provides a measure of the distance between the true distribution of the data and our mixture model. The second term is a complexity term that depends on the classes containing the emission densities and which is proportional to the sum of their VC-indices. We deduce from this deviation bound that the estimator is not only robust with respect to model misspecification but also to contamination and the presence of outliers among the data set. Dealing with models that may be approximate allows to build estimators that possess properties over wider classes of distribution. Ghosal \& Van der Vaart \cite{ghosal2001} used finite location-scale Gaussian mixtures to approximate continuous Gaussian mixtures with compactly supported mixing distribution. They consider mixtures with scale parameters lying between two constants that depend on the true distribution. By using a similar approximation, we show that our estimator achieves the same rate of convergence but without any restriction on the scale parameters so that the model we consider does not depend on the true mixing distribution. In particular, our result is insensitive to translation or rescaling. \par
Under suitable identifiability assumptions and when the distribution of the data belongs to our model, hence is of the form (\ref{eq:p_w_f}), we also analyze the performance of our estimators of the parameters $w_1,\dots,w_K$ and $F_1,\dots,F_K$. In order to establish convergence rates, we relate the Hellinger distance between the distribution of the data and its estimator to a suitable distance between the corresponding parameters. A general technique is using Fisher's information and results of Ibragimov \& Has’minski\u{\i} \cite{Ibragimov} for regular parametric models. We can also use other results specific to parameter estimation in mixture models such as what Gadat \emph{et al} \cite{gadat2018parameter} proved in the context of two component mixtures with one known component. In both situations, we obtain, up to a logarithmic parameter, the usual $1/\sqrt{n}$-rate of convergence for regular parametric models. We also provide the example of a parametric model for which our technics allow us to establish faster convergence rates while classical methods based on the likelihood or the least-squares fail to apply and hence give nothing.\par
In many applications, starting with a single mixture model may be restrictive and a more reasonable approach is to consider candidate ones for estimating the number of components of the mixture and proposing suitable models for the emission densities.  To tackle this problem, we design a model selection procedure from which we establish, under suitable assumptions, an oracle-type inequality. We consider several illustrations of this strategy. For example, we use a penalized estimator to select the number of components of a Gaussian mixture estimator and obtain similar adaptivity results as Maugis \& Michel \cite{maugis2012}. We also consider a model with a fixed number of components but each emission density can either belong to the Gaussian or to the Cauchy location-scale family. We prove that if we know the number of components, we can estimate consistently the proportions of Gaussian and Cauchy components as well as their location and scale parameters. To our knowledge, this result is the first of its kind.\par
The extension of the theory of $\rho$-estimation to mixture models is based on Proposition \ref{prop:rho_dim_vc} below. The proof of this result relies on an upper bound for the expectation of the supremum of an empirical process over a mixture of VC-subgraph classes. It generalizes the result that was previously established for a single VC-subgraph class. The key argument in the proof is the uniform entropy property of VC-subgraph classes that still holds for the overall density mixture model with lower bounded weights.\par
The paper is organized as follows. We describe our statistical framework in Section \ref{sec:framework}. In Section \ref{sec:main_result}, we present the construction of the estimator on a single mixture model. We state the general result for density estimation on a single model and illustrate the performance of the estimator on the specific example of GMMs. The problem of estimating the parameters of the mixture is addressed in the subsection \ref{sec:param_recover}. Finally, Section \ref{sec:model_selection} is devoted to model selection criterion and the properties of the estimator on the selected model. The appendix contains all the proofs that are gathered in the same sections when they are related. Those sections include the main results, density estimation, the parametric estimation in regular parametric models, the case of two-component mixtures with one known component and the lemmas.
\section{The statistical framework}
\label{sec:framework}
We observe $n$ independent random variables $X_1,X_2,\dots,X_n$ with  respective marginal distributions $P^*_1,P^*_2,\dots,P^*_n$ on the measurable space $(\mscr{X},\mcal{X})$. We model the joint distribution $\mbf{P^*} = P^*_1 \otimes P^*_2 \otimes \dots \otimes P ^*_n$ of $\mbf{X}=(X_1,X_2,\dots,X_n)$ by a probability of the form $\overline{P}^{\otimes n}$ doing as if the observations were i.i.d. with common distribution $\overline{P}$. We assume that $\overline{P}$ is a mixture of the form (\ref{eq:p_w_f}) where $K$ is a positive integer, the $w_k$ some positive weights that satisfy $\sum_{k=1}^K w_k=1$, and $F_k$ probability  distributions. In order to model each of these probabilities we introduce a collection $\left\{ \overline{\mscr{F}}_{k,\lambda}; k\geq 1,\lambda \in \Lambda_k \right\}$ of possible models and assume that for each $k\in\{1,\dots,K\}$, $F_k$ belongs to $\cup_{\lambda\in\Lambda_k} \overline{\mscr{F}}_{k,\lambda}$. We denote by $\mscr{Q}_K$ the family of distributions of the previous form. For each $k\geq 1$, we call $F_k$ an emission probability, $\overline{\mscr{F}}_{k,\lambda}$ an emission model, and $\mscr{E}_k=\left\{ \overline{\mscr{F}}_{k,\lambda}; \lambda\in\Lambda_k\right\}$ an emission family. Based on the observation of $\mbf{X}$, our aim is to design an estimator $\hat{P}$ of $\overline{P}$ of the form
\begin{equation}
\hat{P} = \suml_{k=1}^{\hat{K}} \hat{w}_k \hat{F}_k \in \bigcup_{K\geq 1} \mscr{Q}_K
\end{equation}
where $\hat{K}$, $(\hat{w}_k)_{1\leq k\leq \hat{K}}$ and $(\hat{F}_k)_k$ are estimators of $K$, $(w_k)_k$ and $(F_k)_k$ respectively. The classical situation that has been considered in the literature corresponds to the case where the collection $\left\{ \overline{\mscr{F}}_{k,\lambda}; k\geq 1, \lambda\in\Lambda_k \right\}$ reduces to a single emission model $\mscr{F}$, for example the family of Gaussian distributions, and the problem is to estimate $K$ and the emission probabilities $F_k$ under the assumption that they all belong to $\mscr{F}$. This assumption is quite restrictive and we rather consider a collection $\mscr{E}_k$  of candidate models for $F_k$ that may even depend on $k$. We say that $\mscr{E}_k$ is simple when it reduces to a single emission model $\overline{\mscr{F}}_k$ and composite otherwise.\par
In order to evaluate the performance of the estimator $\hat{P}$, we introduce on the set $\pmb{\mscr{P}}$ of all product probabilities on $(\mscr{X}^n,\mcal{X}^{\otimes n})$ the Hellinger-type distance $\mbf{h}$ defined by
\begin{equation}
\label{eq:h_dist}
\mbf{h}(\mbf{Q},\mbf{Q'}) = \sqrt{ \suml_{i=1}^n h^2(Q_i,Q'_i)}, \quad\text{for } \mbf{Q} = \bigotimes_{i=1}^n Q_i, \mbf{Q'}=\bigotimes_{i=1}^n Q'_i \in\pmb{\mscr{P}},
\end{equation}
where $h$ is the Hellinger distance on the set $\mscr{P}$ of probability distributions on $(\mscr{X},\mcal{X})$. We recall that for $Q$, $Q'$ in $\mscr{P}$
\begin{equation*}
h^2(Q,Q') = \frac{1}{2} \int \left(\sqrt{\frac{dQ}{d\mu}} -\sqrt{\frac{dQ'}{d\mu}}\right)^2 d\mu,
\end{equation*}
where $\mu$ is a measure that dominates both $Q$ and $Q'$, the result being independent of $\mu$. \par
\begin{hypothese}
\label{hyp:countable}
For all $k\geq 1$, the set $\Lambda_k$ is at most countable (which means finite or countable) and that for all $\lambda$ in $\Lambda_k$, $\overline{\mscr{F}}_{k,\lambda}$ contains an at most countable subset $\mscr{F}_{k,\lambda}$ which is dense in $\overline{\mscr{F}}_{k,\lambda}$ with respect to the Hellinger distance $h$.
\end{hypothese}
This condition implies that there exists a $\sigma$-finite measure $\mu$ that dominates all the $\overline{\mscr{F}}_{k,\lambda}$ for $k\geq 1$ and $\lambda\in\Lambda_k$. Throughout this paper, we fix such a measure $\mu$ and associate to each emission model $\overline{\mscr{F}}_{k,\lambda}$ a family of density distributions $\overline{\mcal{F}}_{k,\lambda}$ such that $\overline{\mscr{F}}_{k,\lambda} = \left\{ f \cdot \mu; f\in\overline{\mcal{F}}_{k,\lambda} \right\}$. In all the different examples considered $\mu$ is the Lebesgue measure. We assume the following.
\begin{hypothese}
\label{hyp:vc_density_model}
For all $k\geq 1$ and $\lambda\in\Lambda_k$, the family of density distributions $\overline{\mcal{F}}_{k,\lambda}$ is VC-subgraph with VC-index $V_{k,\lambda}\geq 1$.
\end{hypothese}
For more details on VC-subgraph classes we  refer the reader to Van der Vaart \& Wellner \cite{VanDerVaart} (Section 2.6.5) and Baraud \emph{et al} \cite{baraudinventiones} (Section 8). Throughout this paper we shall use the following notation. For $\mbf{P}=P_1\otimes\dots\otimes P_n \in\pmb{\mscr{P}}$ and $\mscr{A}\subset \mscr{P}$, we write
\begin{equation*}
\mbf{h}^2\left(\mbf{P},\mscr{A}\right) = \inf_{Q\in \mscr{A}} \mbf{h}^2\left(\mbf{P},Q^{\otimes n}\right) = \inf_{Q\in\mscr{A}} \suml_{i=1}^n h^2(P_i,Q).
\end{equation*}
For $x\in\mbb{R}$, $\lfloor x\rfloor$ is the only integer satisfying $\lfloor x\rfloor \leq x < \lfloor x\rfloor +1$ and similarly $\lceil x \rceil$ denotes the integer satisfying $\lceil x\rceil - 1 < x \leq \lceil x\rceil$. Moreover, if $x >0$ we write $\log_+(x)=\log(x)\vee 0$. If $A$ is a finite set, we denote its cardinal by $|A|$ and if $A$ is infinite, we write $|A|=\infty$. The notation $C(\theta)$ will mean that the constant $C=C(\theta)$ depends on the parameter or set of parameters $\theta$.
\section{Estimation on a mixture model based on simple emission families}
\label{sec:main_result}
In this section, we assume that the $\mscr{E}_k=\left\{ \overline{\mscr{F}}_k \right\}$ are simple for all $k\geq 1$ and that $\overline{P}$ belongs to $\mscr{Q}_K$ for some known value of $K\geq 1$. This means that we know that $\overline{P}$ is a mixture of at most $K$ emission probabilities $F_1,\dots,F_K$ and that $F_k$ belongs to $\overline{\mscr{F}}_k$ for all $k\in\{1,\dots,K\}$.  Under Assumption \ref{hyp:vc_density_model}, we denote by $V_k$ the VC-index of $\overline{\mscr{F}}_k$.
\subsection{Construction of the estimator on \texorpdfstring{$\mscr{Q}_K$}{QK}}
For $\delta$ in $(0,1/K]$, we define the subset $\mscr{Q}_{K,\delta}$ of $\mscr{Q}_K$ by
\begin{equation}
\label{eq:mscr_q_delta}
\mscr{Q}_{K,\delta} := \left\{ \suml_{k=1}^K w_k F_k \in\mscr{Q} ; w\in\mcal{W}_K, w_k\geq\delta, w_k\in\mbb{Q}, F_k\in\mscr{F}_k \right\}
\end{equation}
where the $\mscr{F}_k$ are the countable and dense subsets of $\overline{\mscr{F}}_k$ provided by Assumption \ref{hyp:countable}. We associate to $\mscr{Q}_{K,\delta}$ the family $\mcal{Q}_{K,\delta}$ of densities  with respect to $\mu$ and the $\rho$-estimator $\hat{P}_{\delta}$ of $\overline{P}$ based on the family $\mcal{Q}_{K,\delta}$. We recall that $\hat{P}_{\delta}$ is defined as follows. Given 
\begin{equation}
\label{eq:psi_rho}
\psi: \begin{array}{|lcl}
[0,+\infty] & \rightarrow & [-1,1]\\
x & \mapsto & \frac{x-1}{x+1}
\end{array},
\end{equation}
we set for $\mbf{x} = (x_1,...,x_n) \in \mscr{X}^n$ and $q,q'\in \mcal{Q}_{K,\delta}$
\begin{equation}
\label{eq:t}
\mbf{T}(\mbf{x},q,q') := \suml_{k=1}^{n}\psi\left(\sqrt{\frac{q'\left(x_i\right)}{q\left(x_i\right)}}\right),
\end{equation}
with the convention $0/0 = 1$ and $a/0 = +\infty$ for all $a >0$, and
\begin{equation}
\label{eq:upsilon_no_pen}
\mbf{\Upsilon}(\mbf{X},q) := \sup_{q'\in\mcal{Q}_{\delta}} \mbf{T}(\mbf{X},q,q').
\end{equation}
The $\rho$-estimator $\hat{P}_{\delta}$ is any measurable element of the closure (with respect to the Hellinger distance) of the set
\begin{equation}
\label{eq:rho_est}
\pmb{\mscr{E}}(\psi,\mbf{X}) := \left\{ Q = q\cdot\mu ; q\in\mcal{Q}_{\delta}, \mbf{\Upsilon}(\mbf{X},q) < \inf_{q'\in\mcal{Q}_{\delta}} \mbf{\Upsilon}(\mbf{X},q')  + 8.24 \right\}.
\end{equation}
This construction follows \cite{baraudinventiones} and the constant 8.24 is given by the choice of $\psi$.
\subsection{The performance of the estimator}
\label{sec:performance_estimator}
The following result holds.
\begin{theoreme} 
\label{th:single_model}
Let $\delta\in(0,1/K]$ and $\xi>0$. Assume that Assumptions \ref{hyp:countable} and \ref{hyp:vc_density_model} hold and set $\overline{V}=V_1+\dots+V_K$. Any $\rho$-estimator $\hat{P}_{\delta}$ on $\mscr{Q}_{K,\delta}$ satisfies with probability at least $1-e^{-\xi}$,
\begin{align}
\label{eq:th_single_model}
\mbf{h}^2 \left(\mbf{P^*}, \left(\hat{P}_{\delta}\right)^{\otimes n}\right) &\leq c_0 \left[ \mbf{h}^2\left(\mbf{P^*},\mscr{Q}_K\right) + n(K-1)\delta \right]\\
&+ c_1 \overline{V} \left[ 5.82 + \log\left( \frac{(K+1)^2}{\delta} \right) + \log_+ \left( \frac{n}{\overline{V}} \right) \right]\nonumber\\
&+ c_2 (1.49+\xi).\nonumber
\end{align}
where $c_0=300$, $c_1=8.8\times 10^5$ and $c_2=5014$. In particular, for the choice $\delta=1$ for $K=1$ and $\delta=\frac{\overline{V}}{n(K-1)}\bigwedge \frac{1}{K}$ otherwise, the resulting estimator $\hat{P}=\hat{P}_{\delta}$ satisfies 
\begin{equation}
\label{eq:th1_2}
C \mbf{h}^2 \left(\mbf{P^*},\hat{P}^{\otimes n}\right) \leq \mbf{h}^2\left(\mbf{P^*},\mscr{Q}_K\right) +  \overline{V}\left[ 1 + \log\left( \frac{K n}{\overline{V}\wedge n} \right) \right] + \xi,
\end{equation}
with probability at least $1-e^{-\xi}$, where $C$ is a universal constant.
\end{theoreme}
Inequality (\ref{eq:th_single_model}) shows the influence of the choice of the parameter $\delta$ on the performance of the estimator $\hat{P}_{\delta}$. Hereafter, we shall choose $\delta$ as in the second part of Theorem \ref{th:single_model} and therefore only comment on inequality (\ref{eq:th1_2}). Given $\overline{P}$ in $\mscr{Q}_K$, it follows from the triangle inequality and the fact that $(a+b)^2\leq 2a^2 +2b^2$ for all non-negative numbers $a$ and $b$, that
\begin{equation*}
n h^2\left( \overline{P}, \hat{P} \right) =  \mbf{h}^2\left( \overline{P}^{\otimes n}, \hat{P}^{\otimes n} \right) \leq 2 \mbf{h}^2\left( \mbf{P^*}, \hat{P}^{\otimes n}\right) + 2 \mbf{h}^2\left( \mbf{P^*}, \overline{P}^{\otimes n} \right).
\end{equation*}
We immediately derive from (\ref{eq:th1_2}) that on a set of probability at least $1-e^{-\xi}$ 
\begin{equation}
\label{eq:independent}
 C h^2 \left( \overline{P},\hat{P} \right) \leq \frac{1}{n}\suml_{i=1}^n  h^2(P^*_i,\overline{P}) + \frac{\overline{V}\log\left( \slfrac{Kn}{\overline{V}}\right)+\xi}{n} .
\end{equation}
In the ideal situation where the observations are i.i.d. with common distribution $\overline{P}\in\mscr{Q}_K$, we obtain that
\begin{equation*}
 C h^2 \left( \overline{P}, \hat{P} \right) \leq \frac{\overline{V}\log\left( \slfrac{Kn}{\overline{V}}\right)+\xi}{n}.
\end{equation*}
Integrating this result with respect to $\xi$ and the fact that $\overline{P}$ is arbitrary in $\mscr{Q}_K$ leads to the uniform risk bound 
\begin{equation}
\label{eq:simple_esperance}
\sup_{\overline{P}\in\mscr{Q}_K} \mbb{E} \left[ h^2 \left( \overline{P}, \hat{P} \right) \right] \leq C' \frac{\overline{V}\log\left( \slfrac{Kn}{\overline{V}}\right)}{n}.
\end{equation}
where $C'$ is a positive universal constant. This means that up to a logarithmic factor, the estimator $\hat{P}$ uniformly converges over $\mscr{Q}_K$ at the rate $1/\sqrt{n}$ with respect to the Hellinger distance.
Our assumption that the families of density functions $\overline{\mcal{F}}_k$ are VC-subgraph is actually weak since it includes situations where these models consist of unbounded densities or densities which are not in $L_2$ which to our knowledge have never been considered in the literature. A concrete example of such situations is the following one. Let $g$ be some non-increasing function on $(0,+\infty)$ which is unbounded and satisfies $\int_0^{+\infty} g(x) dx = \frac{1}{2}$ and $\overline{\mscr{F}}_k$ is the translation model associated to the family of densities $\left\{ x\mapsto g(|x-\theta|) \mathbbm{1}_{|x-\theta|>0} ; \theta\in\mbb{R} \right\}$ for all $k\in\{1,\dots,K\}$. It follows from Proposition 42-(vi) of Baraud \emph{et al} \cite{baraudinventiones} that the VC-index of $\overline{\mcal{F}}_k$ is not larger than $10$.\par
When the data are independent but not i.i.d., we derive from inequality (\ref{eq:independent}) that the estimator $\hat{P}$ performs almost as well as in the i.i.d. case as long as the marginals $P^*_1,\dots,P^*_n$ are close enough to $\overline{P}$. This means that the estimator is robust with respect to a possible misspecification of the model and the departure from the assumption that the data are i.i.d. In particular, this includes the situations where the dataset contains some outliers or has been contaminated. Consider Hüber’s contamination model where a proportion $\epsilon$ of the data is contaminated, i.e. we have $P^*=(1-\epsilon) \overline{P} + \epsilon Q$, where $\overline{P}$ is the probability distribution we want to estimate and $Q$ is the distribution of the contaminated data. In this situation, for any probability distribution $Q$, using (\ref{eq:independent}) we get
\begin{equation*}
 C h^2 \left( \overline{P},\hat{P} \right) \leq \epsilon + \frac{\overline{V}\log\left( n\right)+\xi}{n}.
\end{equation*}
We can see that there is no perturbation of the convergence rate as long as the contamination rate $\epsilon$ remains small as compared to $\overline{V}\log(n)/n$. Inequality (\ref{eq:hell_mix_parameter}), stated later, also allows to consider misspecification for the emission models for example.
\subsection{The case of totally bounded emission models}
\label{sec:models_with_entropy}
We might also consider emission models for which we do not have any bound on the VC-index. For a subset $\mscr{N}$ of $\mscr{P}$ and $\eta\in[0,1]$, the $\eta$-covering number $N(\eta,\mscr{N},h)$ of $\mscr{N}$, with respect to the Hellinger distance, is the minimum number of balls $\mcal{B}_h(P_i,\eta)$, $i=1,\dots,N$, necessary to cover $\mscr{N}$. In that case, the set $\mscr{N}[\eta]=\{P_i;i=1,\dots,N\}$ constitutes a finite approximation of $\mscr{N}$, i.e. for all $Q$ in $\mscr{N}$ there exists $i\in\{1,\dots,N\}$ such that  $h\left( Q,P_i \right) \leq \eta$. We say that $\mcal{N}$ is totally bounded (for the Hellinger distance) if its $\eta$-covering number is finite for all $\eta\in(0,1]$. A direct consequence of the definition of VC-subgraph classes is that any finite set $\mcal{F}$ of real-valued functions is VC-subgraph with VC-index at most $V(\mcal{F})\leq \log_2 \left( |\mcal{F}| \right)$. Consequently, we can still use $\rho$-estimation for models that are not proven to satisfy Assumption \ref{hyp:vc_density_model} but still are such that emission models are totally bounded.
\begin{theoreme}
\label{th:totally_bounded}
Let $\overline{\mscr{F}}_k$ be a totally bounded class of distributions for all $k\in\{1,\dots,K\}$ with $K\geq 2$. Let $\mscr{Q}$ be the mixture model defined by
\begin{equation*}
\mscr{Q}_K = \left\{ \suml_{k=1}^K w_k F_k ; w\in\mcal{W}_K, F_k \in \overline{\mscr{F}}_k, \forall k \in\{1,\dots,K\} \right\}.
\end{equation*}
Assume there are positive constants $A_k$ and $\alpha_k$ such that $\log_2 N(\eta,\mscr{F}_k,h) \leq \left( \frac{A_k}{\eta} \right)^{\alpha_k}$ for all $k$ in $\{1,\dots,K\}$ and for all $\eta\in(0,1)$. Let $\epsilon$ be in $(0,1)$. For $k$ in $\{1,\dots,K\}$, let $\mscr{F}_k[\epsilon]$ be a minimal $\epsilon$-net of $\overline{\mscr{F}}_k$ such that $|\mscr{F}_k[\epsilon]|=N(\epsilon,\mscr{F}_k,h)$. Let $\mscr{Q}_{K,\delta}[\epsilon]$ be the countable model defined by
\begin{equation*}
\mscr{Q}_{K,\delta}[\epsilon] = \left\{ P_{w,F} ; w\in\mcal{W}_K, w_k\geq \delta, w_k\in\mbb{Q}, F_k\in\mscr{F}_k[\epsilon], \forall k \in\{1,\dots,K\} \right\}.
\end{equation*}
Take $\epsilon = n^{-\frac{1}{\alpha_{\infty}+2}}$ and $\delta=\frac{\overline{V}}{n(K-1)}\wedge \frac{1}{K}$ with $\alpha_{\infty}=\max_{1\leq k\leq K} \alpha_k$. There exists a positive constant $C$ such that for any $\rho$-estimator $\hat{P}=\hat{P}_{\delta}$ on $\mscr{Q}_{K,\delta}[\epsilon]$, for all $\xi>0$, we have
\begin{equation*}
C h^2\left( P^*, \hat{P} \right) \leq h^2\left( P^*, \mscr{Q}_K \right) + n^{-\frac{2}{\alpha_{\infty}+2}} \left( 1 + \suml_{k=1}^K A_k^{\alpha_k} \right) \left[ 1 + \log\left( Kn \right) \right] + \xi,
\end{equation*}
with probability at least $1-e^{-\xi}$.
\end{theoreme}
We illustrate this lemma with the following example. Doss \& Wellner \cite{doss2016} provide a bound on the entropy for classes of $\log$-concave and $s$-concave densities. Let $\mcal{C}=\left\{ \varphi : \mbb{R} \rightarrow [-\infty,\infty); \varphi \text{  is a closed, proper concave function}\right\}$ where \emph{proper} and \emph{closed} are defined in \cite{Rockafellar2015} (Sections 4 and 7). For $0<M<\infty$ and $s>-1$, let $\mcal{P}_{M,s}$ be the class of densities defined by
\begin{equation*}
\mcal{P}_{M,s} = \left\{ p\in\mcal{P}_s ; \sup_{x\in\mbb{R}} p(x) \leq M, 1/M\leq p(x) \text{  for all  } |x|\leq 1 \right\},
\end{equation*} 
where $\mcal{P}_s=\left\{ p : \int pd\lambda =1 \right\} \bigcap h_s \circ\mcal{C}$, $\lambda$ is the Lebesgue measure on $\mbb{R}$ and $h_s:\mbb{R}\rightarrow\mbb{R}$ is defined by
\begin{equation*}
h_s(y) = \begin{cases} e^y, & s=0\\
(-y)_+^{1/s}, &s\in (-1,0),\\
y_+^{1/s}, &s>0.
\end{cases}
\end{equation*}
We fix such values of $M$ and $s$. Let $\mcal{Q}_K$ be the density model of mixtures of $s$-concave densities (or $\log$-concave for $s=0$) defined by
\begin{equation*}
\mcal{Q}_K = \left\{ \suml_{k=1}^K w_k f_k ; w\in\mcal{W}_K, f_k\in\mcal{P}_{M,s} \right\},
\end{equation*}
with $K\geq 2$. Let $\mscr{Q}_K$ be the class of distributions associated to $\mcal{Q}_K$. The class $\mcal{P}_{M,s}$ is not proven to be VC-subgraph but it is totally bounded. As a direct consequence of Theorem 3.1 of Doss \& Wellner \cite{doss2016}, there exists a positive constant $A$ such that for all $\epsilon$ in $(0,1]$, we have 
\begin{equation*}
\log_2 N(\epsilon,\mcal{P}_{M,s},h) \leq A \epsilon^{-1/2}.
\end{equation*}
In particular, it means there exists a $\epsilon$-net $\mcal{P}_{M,s}[\epsilon]$ such that $|\mcal{P}_{M,s}[\epsilon]|\leq 2^{C/\epsilon^{-1/2}}$.  Let $\mcal{Q}_{K,\delta}[\epsilon]$ be the countable density model given by  
\begin{equation*}
\mcal{Q}_{K,\delta}[\epsilon] = \left\{ \suml_{k=1}^K w_k f_k ; w\in\mcal{W}_K, w_k\geq \delta, w_k\in\mbb{Q}, f_k\in\mcal{P}_{M,s}[\epsilon] \right\}.
\end{equation*}
One can check that $\mscr{Q}_{K,\delta}[\epsilon]$ is also a $\epsilon$-net of $\mscr{Q}_{K,\delta}$ with respect to the Hellinger distance.
\begin{coro}
\label{coro:log_concave}
Assume there exists $P^*$ in $\mscr{P}$ such that $\mbf{P^*}=(P^*)^{\otimes n}$. Take $\epsilon= n^{-2/5}$ and $\delta=n^{-4/5}\wedge K^{-1}$. Let $\hat{P}=\hat{P}_{\delta}$ be a $\rho$-estimator on $\mscr{Q}_{K,\delta}[\epsilon]$. For all $\xi>0$, we have
\begin{align*}
C h^2 \left( P^*, \hat{P} \right) &\leq h^2\left( P^*,\mscr{Q}_K\right)+ \frac{K}{n^{4/5}} \left[ 1 + \log\left( K n \right) \right] + \frac{\xi}{n},
\end{align*}
with probability at least $1-e^{-\xi}$.
\end{coro}
This result provides a risk bound over the class of distributions associated to mixtures of $s$-concave densities. Up to a logarithmic factor, the estimator $\hat{P}$ uniformly converges over $\mscr{Q}_K$ at the rate $n^{-2/5}$ with respect to the Hellinger distance, which is the same rate given in Theorem 3.2 of Doss \& Wellner \cite{doss2016} for the MLE over $\mcal{P}_{M,s}$.
\subsection{Application to the estimation of a continuous Gaussian mixture}
\label{sec:gmm_vdv}
We denote by $\phi_{\sigma}$ the density function of the 	normal distribution (with respect to the Lebesgue measure on $\mbb{R}$) with mean $0$ and variance $\sigma^2 >0$, i.e.
\begin{equation}
\label{eq:phi_gauss}
\phi_{\sigma} : x \mapsto \frac{1}{\sqrt{2\pi\sigma^2}}e^{-\frac{x^2}{2 \sigma^2}}.
\end{equation}
We assume $P^*$ is of the following form or is close enough to a distribution of the form
\begin{equation*}
p_H(x) = \int \phi_{\sigma}(x-z) dH(z,\sigma), \forall x\in\mbb{R}.
\end{equation*}
We say that $p_H$ is the Gaussian mixture density with  mixing distribution $H$. We want to approximate any distribution of this form with finite Gaussian mixtures, i.e. distribution with densities of the same form with mixing distribution supported on a finite set. For a mixing measure $H$ on $\mbb{R}\times\mbb{R}^{+*}$, we denote by $\text{supp}(H)$ its support. To obtain an approximation result, we need to consider mixing measures $H$ that are supported on a compact set, i.e. there exist $A\geq 0$ and $R\geq 1$ such that $\text{supp}(H)\subset [-A,A]\times[1,R]$. The Hellinger distance being invariant to translation and rescaling, we consider the following class of densities. For $A>0$ and $R\geq 1$ we define
\begin{equation*}
\mcal{C}(A,R) = \bigg\{ p_H ; \exists l\in\mbb{R},  \exists s >0,  \text{supp}(H) \subset \left[ l-sA, l+sA \right] \times [s,sR]  \bigg\}
\end{equation*}
and we denote by $\mscr{C}(A,R)$ the associated class of distributions. We denote by $\mcal{G}$ the location-scale Gaussian family of probability density functions, i.e.
\begin{equation}
\label{eq:gauss_location_scale}
\mcal{G} = \left\{ x \mapsto \phi_{\sigma}(x-\mu) ; \mu \in\mbb{R}, \sigma > 0 \right\}.
\end{equation}
We denote by $\mscr{G}_K$ the Gaussian mixture model with $K$ components associated to class of densities $\mcal{G}_K$ defined by
\begin{equation*}
\mcal{G}_K := \left\{\suml_{k=1}^K w_k \phi_{\sigma_k}(\cdot-z_k) ; w\in\mcal{W}_K, \sigma_k\in (0,+\infty), z_k\in\mbb{R}, \forall k\in\{1,\dots,K\} \right\}.
\end{equation*}
This situation corresponds to $\overline{\mcal{F}}_k=\mcal{G}$ for all $k\in \{1,\dots,K\}$. We can approximate the class $\mscr{C}(A,R)$ with the model $\mscr{G}_K$ as indicated by the following result.
\begin{prop}
\label{prop:ghosal_approximation}
For $K \geq (2/3)^3 A^4$, we have
\begin{equation*}
\sup_{p_H\in\mcal{C}(A,R)} h^2(P_H, \mscr{G}_K) \leq \frac{1}{2} \exp\left( - K^{1/2} \frac{3 \sqrt{3}}{\sqrt{2} R^2} \right) \left[ K^{1/4} \frac{3 \sqrt{2}}{ \sqrt{e\pi} 7^{1/4}} + R \right].
\end{equation*}
\end{prop}
We can deduce a deviation bound on the estimation over $\mscr{C}(A,R)$ from this last result and Theorem \ref{th:single_model}.
\begin{theoreme}
\label{th:gmm}
Assume  $n\geq \exp(2 (A/R)^2)$ and $\frac{n}{\log^2(n)} \geq 2 R^2 /27$. Let $\hat{P}$ be a $\rho$-estimator on $\mscr{G}_{K,\delta}$ with $\delta$ as in (\ref{eq:th1_2}) and $K= \lceil 2 R^4 \log^2(n) /27\rceil$. Assume the true distribution is i.i.d., i.e. $\mbf{P^*}=(P^*)^{\otimes n}$. There exists a numeric constant $C>0$ such that for all $\xi>0$, with probability at least $1-e^{-\xi}$, we have
\begin{equation}
\label{eq:th_continuous_gaussian_mixture}
C h^2 \left( P^*, \hat{P} \right) \leq h^2(P^*,\mscr{C}(A,R)) + \frac{R^4 \log^3(n)+\xi}{n}.
\end{equation}
\end{theoreme}
Therefore, for a fixed $R$, we obtain a rate of $\log^{3/2}(n)/\sqrt{n}$ over $\mscr{C}(A,R)$ with respect to the Hellinger distance. We can also consider larger classes of distributions, with $R$ increasing as $n$ increases but it would deteriorate this rate. Our result is still an improvement of Theorem 4.2 from \cite{ghosal2001} as it requires weaker assumptions. Their result is sensitive to translation or scaling and they have to specify bounds $0<\underline{\sigma}<\overline{\sigma}$ in the model such that $H^*$ is supported on a compact set $[-a,a]\times [\underline{\sigma},\overline{\sigma}]$. Moreover, our estimator is robust, to contamination for instance. Assume we have an $\epsilon$ contamination rate of our data, i.e. $P^*$ is of the form $P^*=(1-\epsilon)\overline{P} + \epsilon Q$ with $\epsilon\in(0,1)$, $\overline{P}\in\mscr{C}(A,R)$ and $Q$ is any probability distribution. Then, our estimator satisfies $Ch^2(P^*,\hat{P}) \leq \epsilon + \frac{R^4 \log^3(n)+\xi}{n}$ on an event of probability $1-e^{-\xi}$. As long as $\epsilon$ remains small as compared to $R^4\log^3(n)/n$, the rate is not deteriorated by the contamination.
\subsection{Parameter estimation}
\label{sec:param_recover}
We say that $\hat{w}$ and $\hat{F}$ are $\rho$-estimators if the resulting mixture distribution $\hat{P}$ given by
\begin{equation*}
\hat{P}=\suml_{k=1}^K \hat{w}_k \hat{F}_k
\end{equation*}
is a $\rho$-estimator. We have a general result for the performance of $\hat{P}$ but not for $\hat{w}$ and $\hat{F}$. In order to evaluate the performance of these estimators, we first need to ensure that the parameters $w=(w_1,\dots,w_K)$ and $F=(F_1,\dots,F_K)$ are identifiable.
\begin{exemple} Let $\overline{\mscr{F}}$ be the set of uniform distributions $\mcal{U}(a,b)$ the uniform distribution on the interval $(a,b)$ of positive lengths. Then the parameters $w$ and $F$ in the mixture model 
\begin{equation*}
\mscr{Q}_2 = \left\{ w_1 F_1 + (1-w_1) F_2 ; w_1\in(0,1), F_1,F_2\in\overline{F} \right\}
\end{equation*}
are not identifiable since 
\begin{equation*}
\frac{3}{4}\mcal{U}(0,1) + \frac{1}{4}\mcal{U}(1/3,2/3) = \frac{1}{2}\mcal{U}(0,2/3) + \frac{1}{2}\mcal{U}(1/3,1).
\end{equation*}
\end{exemple}
We shall say that $P=P_{w,F}$ is identifiable (with respect to the model) if for all $v$ in $\mcal{W}_K$ and all $G$ in $\mscr{F}_1\times\dots\times\mscr{F}_K$, we have 
\begin{equation*}
P_{w,F}=P_{v,G} \Rightarrow \exists \tau\in\mscr{S}_K, \forall k\in\{1,\dots,K\},w_k=v_{\tau(k)} \text{  and  } F_k=G_{\tau(k)},
\end{equation*}
where $\mscr{S}_K$ denotes the set of all permutations of $\{1,\dots,K\}$.
There is a wide literature about identifiability that includes the works of Teicher \cite{teicher1961}, Sapatinas \cite{sapatinas1995}  and Allman \emph{et al} \cite{allman2009} for example. Identifiability is a minimum requirement for the parameter estimators to be meaningful but we can hardly get more than consistency with it. As mentioned in the introduction, we are looking for a lower bound on the Hellinger distance between mixture distributions. Convexity properties ensure that we always have the upper bound
\begin{equation}
\label{eq:hell_mix_parameter}
h\left( P_{w,F}, P_{v,G} \right) \leq \inf\limits_{\tau\in\mscr{S}_K} \left\{ h( w, v \circ \tau) + \max_{k\in[K]} h\left( F_k, G_{\tau(k)} \right) \right\},
\end{equation}
for all mixing weights and emission distributions (see Lemma \ref{lem:upper_bound_parameters}). On the other hand, obtaining a lower bound on $h^2\left( P_{\pi,F}, P_{\nu,G} \right)$ is quite more complicated unfortunately. There are still some situations where we do have such a lower bound.
\subsubsection*{Regular parametric model}
Let $K$ be an integer larger than 1. We consider parametric emission models associated to density models $\overline{\mcal{F}}_k=\left\{ f_k(\cdot;\alpha), \alpha \in A_k \right\}$, where $A_k$ is a subset of $\mbb{R}^{d_k}$ for all $k\in\{1,\dots,K\}$. It is always possible to find a countable dense subset of $A_k$ with respect to the Euclidean distance on $\mbb{R}^{d_k}$. We assume there is a reasonably good connection between the Hellinger distance on the emission models and the Euclidean distances on the parameter spaces such that a dense subset of $A_k$ would translate into a dense subset of the emission model with respect to the Hellinger distance. This assumption is very weak and does not seem to be restrictive in any way. In the different examples we consider we can always consider $A_k\cap\mbb{Q}^{d_k}$ as a dense subset of $A_k$. Therefore Assumption \ref{hyp:countable} is satisfied with $\mcal{F}_k=\left\{ f_k(\cdot;\alpha), \alpha \in B_k \right\}$. We denote by $\mscr{Q}_K$ the distribution model associated to the mixture density model
\begin{equation*}
\mcal{Q}_K = \left\{ p(\cdot;\theta) = \suml_{k=1}^{K-1} w_k f_k(\cdot;z_k) + (1-w_1-\dots-w_{K-1}) f_K(\cdot;\alpha_K) ; \theta=(w,\alpha) \in \Theta \right\},
\end{equation*}
where $\Theta$ is an open convex subset of $\left\{ w \in (0,1)^{K-1} ; \suml_{k=1}^{K-1} w_k < 1 \right\} \times A_1 \times \dots \times A_K$. We make the following assumptions.
\begin{hypothese}
\label{hyp:regular}
\begin{itemize}
\item[a)] The function $z\mapsto f_k(x;z)$ is continuous on $A_k$ (with respect to the Euclidean distance) for $\mu$-almost all $x\in\mscr{X}$, for all $k\in\{1,\dots,K\}$.
\item[b)] For all $k\in\{1,\dots,K\}$, for $\mu$-almost all $x\in\mscr{X}$ the function $u\mapsto f_k(x;u)$ is differentiable at the point $u=\alpha$ and for all $j\in\{1,\dots,d_k\}$, we have
\begin{equation*}
\int_{\mscr{X}} \left| \frac{\partial f_k(x;\alpha)}{\partial \alpha_j} \right|^2 \frac{\mu(dx)}{f_k(x;\alpha)} <\infty.
\end{equation*}
\item[c)] The function $\theta\mapsto \psi(\cdot;\theta)=\frac{\partial}{\partial \theta} p^{1/2}(\cdot;\theta)$ is continuous in the space $L_2(\mu)$.
\item[d)] The class of densities $\overline{\mcal{F}}_k$ is VC-subgraph with VC-index not larger than $V_k$ for all $k\in\{1,\dots,K\}$. We write $\overline{V}=V_1=\dots+V_k$.
\end{itemize}
\end{hypothese}
We use the approach of Ibragimov and Has’minski\u{\i} \cite{Ibragimov} for regular parametric models to obtain a deviation inequality on the Euclidean distance between parameters using Fisher's information.
\begin{theoreme}{(Theorem 7.6 \cite{Ibragimov})}\\
\label{th:regular_parametric}
Let $\overline{\theta}$ be in $\Theta$. Assume the Fisher's information matrix
\begin{equation*}
I\left(\overline{\theta}\right) = \int_{\mscr{X}} \frac{\partial p\left(x;\overline{\theta}\right)}{\partial \theta} \left( \frac{\partial p\left(x;\overline{\theta}\right)}{\partial \theta} \right)^T \frac{\mu(dx)}{p\left(x;\overline{\theta}\right)}
\end{equation*}
is definite positive and $\inf_{ \substack{||\overline{\theta}-\theta||\geq a\\ \theta \in\Theta}} h^2\left(P_{\overline{\theta}},P_{\theta}\right) > 0$ for all $a>0$. Let $\hat{P}=P_{\hat{w},\hat{F}}$ be a $\rho$-estimator on $\mscr{Q}_{K,\delta}$, with $\delta$ as in (\ref{eq:th1_2}). There exists a positive constant $C\left(\overline{\theta}\right)$ such that for all $\xi>0$, with probability at least $1-e^{-\xi}$, we have
\begin{equation}
\label{eq:th_regular_parametric}
C\left(\overline{\theta}\right) \left( || \overline{w} - \hat{w}||^2 +  \suml_{k=1}^K 1\wedge ||\overline{\alpha}_k-\hat{\alpha}_k||^2 \right) \leq \frac{1}{n} \left[ \mbf{h}^2\left(\mbf{P^*},P_{\overline{\theta}}^{\otimes n}\right) + \overline{V}\log(n)+\xi \right].
\end{equation}
And assuming $P^*=P_{\overline{\theta}}$, we obtain the usual parametric convergence rate up to a logarithmic factor for the parameter estimators.
\end{theoreme}
Inequality (\ref{eq:th_regular_parametric}) proves that even if "true parameters" might not exist the parameter estimators can be meaningful as long as $\mbf{P^*}$ is relatively close to the model. The Gaussian mixture model is the most common mixture model and it is a regular parametric model. Let $K\geq 2$ and take $\mcal{F}_k=\mcal{G}$ for all $k\in\{1,\dots,K\}$. We define a binary relation on $\mbb{R}\times(0,\infty)$ by
\begin{equation}
\label{eq:order_parameter}
(z_1,\sigma_1) > (z_2,\sigma_2) \Leftrightarrow \begin{cases}
\sigma_1 >\sigma_2;\\
\text{or  }\sigma_1=\sigma_2 \text{ and  } z_1>z_2.
\end{cases}
\end{equation}
We consider the parameters $\theta=(w_1,\dots,w_{K-1},z_1,\sigma_1^2,\dots,z_K,\sigma_K^2)$ belonging to the set
\begin{equation*}
\Theta=\left\{ \theta \in(0,1)^{K-1}\times  \left( \mbb{R}\times\mbb{R^*} \right)^K; \suml_{k=1}^{K-1} w_k <1, \left(z_1,\sigma_1\right) > \dots > \left(z_K,\sigma_K\right)  \right\}.
\end{equation*}
\begin{theoreme}
\label{th:parameter_recovery_gmm}
Assume $P^*=P_{\overline{\theta}}=\suml_{k=1}^K \overline{w}_k \mcal{N}(\overline{z}_k,\overline{\sigma}^2_k)$ such that $(\overline{z}_1,\overline{\sigma}_1)>\dots>(z_K,\sigma_K)$ are all distinct and $\inf\limits_{1\leq k\leq K} \overline{w}_k >0$. Let $\hat{P}$ be a $\rho$-estimator on $\mscr{G}_{K,\delta}$, with $\delta$ as in (\ref{eq:th1_2}). There exists a positive constant $C\left(\overline{\theta}\right)$ such that, for all $\xi>0$, we have
\begin{equation}
\label{eq:th_gmm}
C\left(\overline{\theta}\right) \left( \suml_{k=1}^{K-1} ||\overline{w}_k-\hat{w}_k||^2 + \suml_{k=1}^K \left|\left|\left(\overline{z}_k,\overline{\sigma}_k^2\right) - \left(\hat{z}_k,\hat{\sigma}_k^2\right)\right|\right|^2 \wedge 1 \right) \leq \frac{5K\log(n)+\xi}{n},
\end{equation}
with probability at least $1-e^{-\xi}$.
\end{theoreme}
Our estimator reaches the optimal rate of convergence up to a logarithmic factor. One can notice that the assumption of ordered couples of parameters $(z_j,\sigma_j^2)$ can be replace by considering distinct couples only and taking the infimum over permutation of the hidden states in (\ref{eq:th_gmm}).
\subsubsection*{Connection with the $L_2$-distance}
We can use results from the literature that do not apply to the Hellinger distance but to other ones such as the $L_2$-distance between densities. There is a general inequality  between the $L_2$ and Hellinger distances when the density functions are bounded, i.e.
\begin{equation}
\label{eq:l2_hellinger}
||p-q||_2^2 \leq 4 \left(||p||_{\infty} + ||q||_{\infty}\right) h^2(P,Q).
\end{equation}
Assume one can prove an inequality of the following type. For any $w,v$ in $\mcal{W}_K$ and any $f_k,g_k$ in $\overline{\mcal{F}}_k$ for all $k\in\{1,\dots,K\}$ such that the resulting mixtures belong to our model, we have
\begin{equation}
\label{eq:L2_param}
\underline{c} \left( d_{\Pi}^2(w,v) + \max_{k\in[K]} d_F^2(f_k,g_k) \right) \leq \left|\left| \suml_{k=1}^K w_k f_k - \suml_{k=1}^K v_k g_k \right|\right|^2,
\end{equation}
where $d_{\Pi}$ is a distance on $\mcal{W}_K$ and $d_F$ is a distance on $\bigcup_{1\leq k\leq K} \overline{\mcal{F}}_k$. Moreover, assuming the density models $\overline{\mcal{F}}_k$ are uniformly bounded, i.e.
\begin{equation}
\label{eq:bounded_densities}
\sup_{k\in[K]} \sup_{f\in\overline{\mcal{F}}_k} ||f||_{\infty} =: U <\infty,
\end{equation}
we get
\begin{equation*}
d_{\Pi}^2(w,v) + \max_{k\in[K]} d_F^2(f_k,g_{\tau(k)}) \leq \frac{8U}{\underline{c}} h^2 \left( \suml_{k=1}^K w_k F_k, \suml_{k=1}^K v_k G_k \right). 
\end{equation*}
Here again, a density estimation result implies a result for the parameter estimation. We can apply this method to the special case of two-component mixture model with one known component. Let $\phi$ be a density function on $\mbb{R}^d$ with respect to the Lebesgue measure. We consider the 2-component mixture model $\mscr{Q}$ associated to the class of densities
\begin{equation}
\label{eq:contamination_model}
\mcal{Q} = \left\{ x \mapsto p_{\lambda,z}(x) = (1-\lambda) \phi(x) + \lambda \phi(x-z) ; \lambda \in[0,1], z\in\mbb{R}^d \right\},
\end{equation}
with $\overline{\mcal{F}}_1=\left\{ \phi \right\}$ and $\overline{\mcal{F}}_2=\left\{ x\mapsto \phi(x-z) ; z \in \mbb{R}^d \right\}$. We make the following assumptions on $\phi$.
\begin{hypothese}
\label{hyp:phi_contamination}
The function $\phi$ belongs to $\mcal{C}^3\left(\mbb{R}^d\right)\cap\mbb{L}^2\left(\mbb{R}^d\right)$. For any $M>0$, there exists a function $g$ in $\mbb{L}^2\left(\mbb{R}^d\right)$ such that
\begin{equation*}
\forall x\in\mbb{R}^d, \forall z\in[-M,M]^d, |\phi(x)-\phi(x-z)|\leq ||z|| g(x)
\end{equation*}
and
\begin{equation*}
\int g^2(x) \phi^{-1}(x) dx < + \infty.
\end{equation*}
\end{hypothese}
Gadat \emph{et al} proved an inequality such as (\ref{eq:L2_param}) in this situation.
\begin{prop}{(inequality (7.11) in \cite{gadat2018parameter})}\\
\label{prop:gadat}
Under Assumption \ref{hyp:phi_contamination}, for all $M > 0$, there exists a positive constant $c(\phi,M)$ such that for all $z_1,z_2 \in[-M,M]^d$ and $\lambda_1,\lambda_2 \in [0,1]$,
\begin{equation*}
c(\phi,M) ||z_1||^2 \left( ||z_2||^2 \left(\lambda_1-\lambda_2\right)^2 + \left(\lambda_1\right)^2  \left|\left|z_1-z_2\right|\right|^2 \right) \leq || p_{\lambda_1,z_1} - p_{\lambda_2,z_2}||^2.
\end{equation*}
\end{prop}
One can notice that Assumption \ref{hyp:phi_contamination} implies that $\phi$ is bounded (see Assumption ($\mbf{H}_{\mcal{S}}$) in \cite{gadat2018parameter}). Hence, we can deduce a deviation inequality for $\rho$-estimators of parameters.
\begin{theoreme}
\label{th:gadat}
We assume $\overline{\mcal{F}}_2$ has a finite VC-index $V$, $\lambda^*\in(0,1]$ and $z^*\neq 0$. For $\delta$ as in (\ref{eq:th1_2}), there exists a positive constant $C(\phi,\lambda^*,z^*)$ and an integer $n_0=n_0(\phi,\lambda^*,z^*)$ such that for any $\rho$-estimator $\hat{P}=P_{\hat{\lambda},\hat{z}}$ on $\mscr{Q}_{\delta}$, $n\geq n_0$ and for all $\xi\in (0,\xi_n)$, we have
\begin{equation*}
C(\phi,z^*,\lambda^*) \left( \left(\lambda^*-\hat{\lambda}\right)^2 + \left( \left|\left|z^*-\hat{z}\right|\right|^2 \wedge 1 \right) \right) \leq \frac{\xi+(V+1)\log(n)}{n},
\end{equation*}
with probability at least $1-e^{-\xi}$, where $\xi_n= (1+V)[1+\log(2n/(1+V))])$.
\end{theoreme}
This implies the consistency of $\hat{z}$ and consequently the consistency of $\hat{\lambda}$ if $z^*\neq 0$. We can deduce a bound on the convergence rate for $\hat{z}$ and also for $\hat{\lambda}$ but only for $n$ large enough. It is similar to Theorem 3.1 of Gadat \emph{et al.} \cite{gadat2018parameter} with a smaller power for the logarithmic term. This slight improvement is allowed by the VC assumption. Furthermore, we do not need to know a value of $M$ such that $z^*\in[-M,M]$ or to specify it in the model. The examples of translation families taken by Gadat \emph{et al} \cite{gadat2018parameter} (Section 6) all satisfy the VC assumption.
\begin{lemme}
\label{lem:vc}
\begin{itemize}
\item The Cauchy location-scale family $\mcal{C}$ of density functions defined by (\ref{eq:cauchy_location_scale}) is VC-subgraph with VC-index $V(\mcal{C})\leq 5$. 
\item The family of densities $\mcal{G}$ defined by (\ref{eq:gauss_location_scale}) is VC-subgraph with VC-index at most $5$. This bound extend to $3+\frac{d(d+3)}{2}$ for multivariate normal distributions in dimension $d$.
\item The Laplace location family $\mcal{L}$ of density functions defined by
\begin{equation*}
\mcal{L} = \left\{ x\mapsto \frac{1}{2} e^{-|x-z|} ; z\in\mbb{R} \right\}
\end{equation*}
is VC-subgraph with VC-index $V(\mcal{L})\leq $.
\item The location family of densities $\mcal{SG}_{\alpha}$ associated to the skew Gaussian density defined by
\begin{equation*}
\mcal{SG}_{\alpha} = \left\{ x\mapsto 2 \phi_1(x-z) \int_{-\infty}^{x-z} \phi_1(\alpha t) dt ; z\in\mbb{R}\right\}
\end{equation*}
is VC-subgraph with VC-index $V(\mcal{SG}_{\alpha})\leq 10$ for all $\alpha\in\mbb{R}$, where $\phi_1$ is given by (\ref{eq:phi_gauss}).
\end{itemize}
\end{lemme}
By inclusion, if the bound holds for the location-scale family it also holds for the location family with fixed scale parameter.
\subsubsection*{Proving a lower bound for a specific example}
\label{sec:faster_rate}
In some specific situations, it is relatively easy to prove a lower bound on the Hellinger distance. This is what we do in the following example and it allows us to obtain faster rates than the usual parametric one. Let $\alpha$ be in $(0,1)$. We denote by $s_{\alpha}$ the probability density function with respect to the Lebesgue measure on $\mbb{R}$ defined by
\begin{equation*}
s_{\alpha} : x\in\mbb{R} \mapsto \frac{1-\alpha}{2|x|^{\alpha}} \mathbbm{1}_{|x|\in(0,1]}.
\end{equation*}
We consider $\mcal{Q}$ as in (\ref{eq:contamination_model}) with $\phi=s_{\alpha}$ and for $\lambda\in[0,1]$ and $z\in\mbb{R}$, we write
\begin{equation*}
p_{\lambda,z} = (1-\lambda) s_{\alpha} + \lambda s_{\alpha}(\cdot-z).
\end{equation*}
We can prove that the Hellinger distance $h(P_{\lambda,z},P_{\lambda',z'})$ is lower bounded by some distance between the parameters which leads to the following theorem.
\begin{theoreme}
\label{th:faster_rate}
For $\lambda^*>0$ and $z^*\neq 0$, there is a positive constant $C(\alpha,z^*,\lambda^*)$ such that, for any $\rho$-estimator $\hat{P}=P_{\hat{\lambda},\hat{z}}$ on $\mscr{Q}_{\delta}$ with $\delta=10/n$ and $n\geq 20$, for all $\xi>0$, with probability at least $1-e^{-\xi}$ we have
\begin{equation*}
 C(\alpha,z^*,\lambda^*) \left[ 1 \wedge |\hat{z}-z^*|^{1-\alpha} + \left( \lambda^* - \hat{\lambda} \right)^2 \right] \leq \frac{\log(n) + \xi}{n}.
\end{equation*}
\end{theoreme}
We derive from this inequality that our estimators $\hat{\lambda}=\hat{\lambda}_n$ and $\hat{z}=\hat{z}n$ estimate $\lambda^*$ and $z^*$ at a rate which is at least $\sqrt{(\log n)/n}$ and $(n^{-1}\log n)^{1/(1-\alpha)}$ respectively. This latter rate is faster than the usual $1/\sqrt{n}$-rate for all $\alpha\in(0,1)$. Up to the logarithmic factors, these rates are optimal. Moreover, one can notice that both maximum likelihood and least squares approaches do not apply here since we consider density functions that are unbounded, and not even square integrable for $\alpha\in[1/2,1)$.
\section{Model selection}
\label{sec:model_selection}
In Section \ref{sec:main_result} we consider estimation on a model with a fixed order $K$ and simple emission families. We use model selection to overcome this restriction in this section and consider composite emission families and/or models with different orders.  
\subsection{Construction of the estimator}
Let $\Theta$ be a subset of
\begin{equation*}
\bigcup_{K\geq 1} \{K\} \times \prod_{k=1}^K \Lambda_k.
\end{equation*}
Let $\delta : \Theta \rightarrow (0,1]$ be such that for $\theta=(K,\lambda_1,\dots,\lambda_K)\in \Theta$, $\delta(\theta)\in(0,1/K]$. We write
\begin{equation*}
\mscr{Q}_\delta ( \theta ) = \left\{ \suml_{k=1}^K w_k F_k ; w\in\mcal{W}_K, w_k\geq \delta, w_k \in \mbb{Q}, F_k\in\mscr{F}_k,\forall k\in[K] \right\}.
\end{equation*}
We define $\mscr{Q}_{\delta}$ by
\begin{equation*}
\mscr{Q}_{\delta} = \bigcup_{\theta \in \Theta} \mscr{Q}_{\delta}(\theta).
\end{equation*}
We associate to $\mscr{Q}_{\delta}$ the family $\mcal{Q}_{\delta}$ of densities  with respect to $\mu$ and the $\rho$-estimator $\hat{P}_{\delta}$ of $\overline{P}$ based on the family $\mcal{Q}_{\delta}$.
Assuming we have a penalty function $\pen : \mcal{Q}_{\delta} \rightarrow \mbb{R}$, we set
\begin{equation}
\label{eq:upsilon_pen}
\mbf{\Upsilon}(\mbf{X},q) = \sup_{q'\in\mcal{Q}_{\delta}} \left[ \mbf{T}(\mbf{X},q,q') - \pen(q') \right] + \pen(q),
\end{equation}
for all $q\in\mcal{Q}_{\delta}$. The $\rho$-estimator $\hat{P}_{\delta}$ is any measurable element of the closure (with respect to the Hellinger distance) of the set $\pmb{\mscr{E}}(\psi,\mbf{X})$, as defined by (\ref{eq:rho_est}). One can notice that a constant penalty function does not have any impact on the definition of $\mbf{\Upsilon}$ and brings us back to the previous situation.
\subsection{Estimation on a mixture model based on composite emission families}
Let $K$ be larger than or equal to 2. Let $L$ be a subset of $\prod_{k=1}^K \Lambda_k$ and define $\Theta$  by $\Theta=\{K\}\times L$, i.e. $K$ is fixed. For $\lambda=(\lambda_1,\dots,\lambda_K)\in L$, the model $\mscr{Q}(\lambda)$ is a subset of
\begin{equation*}
\left\{ \suml_{k=1}^K w_k F_k ; w\in\mcal{W}_K, F_k\in\overline{\mscr{F}}_{\lambda_k},\forall k\in[K] \right\}
\end{equation*}
and we define its countable subset $\mscr{Q}_{\delta} ( \lambda )$ by
\begin{equation*}
\mscr{Q}_{\delta} = \left\{ \suml_{k=1}^K w_k F_k \in\mscr{Q}(\lambda) ; w\in\mcal{W}_K, w_k\geq \delta(\lambda), w_k \in \mbb{Q}, F_k\in\mscr{F}_{\lambda_k},\forall k\in[K] \right\},
\end{equation*}
where $\delta$ is any function $L \rightarrow (0,1/K]$, and $\mscr{Q}_{\delta}=\bigcup_{\lambda\in L} \mscr{Q}_{\delta}(\lambda)$. Under Assumption \ref{hyp:vc_density_model}, we write $\overline{V}(\lambda) = V(\lambda_1) + \dots + V(\lambda_K)$.
\begin{theoreme}
\label{th:selection_emission_models}
Let $\Delta$ be a mapping $L \rightarrow \mbb{R}^+$ such that $\suml_{\lambda\in  L} e^{-\Delta(\lambda)} \leq 1$. Let $\pen$ be the penalty function defined by
\begin{equation*}
\pen(q) = \kappa \inf\limits_{\lambda \in L|Q\in\mscr{Q}(\lambda)} \left[ 174.1 \overline{V}(\lambda) \left[ 5.82 + \log\left( \frac{(K+1)^2}{\delta(\lambda)} \right) + \log_+ \left( \frac{n}{\overline{V}(\lambda)} \right) \right] + \Delta(\lambda) \right],
\end{equation*}
where $\kappa$ is given by (19) in \cite{baraudinventiones}. Assume there is $P^*$ in $\mscr{P}$ such that $\mbf{P^*}=(P^*)^{\otimes n}$. For the choice $\delta(\lambda)=\frac{\overline{V}(\lambda)}{n(K-1)}\bigwedge \frac{1}{K}$, there is a positive constant $C$ such that the resulting estimator $\hat{P}=\hat{P}_{\delta}$ satisfies the following. For all $\xi>0$, with probability at least $1-e^{-\xi}$ we have
\begin{equation*}
C h^2(P^*,\hat{P}) \leq \inf_{\lambda\in L} \left\{ h^2( P^*, \mscr{Q}(\lambda) ) + \frac{1}{n}  \left( \overline{V}(\lambda) \left[ 1 + \log\left( \frac{ K n }{ \overline{V}(\lambda) \wedge n} \right) \right] + \Delta(\lambda) + \xi \right) \right\}.
\end{equation*}
\end{theoreme}
This is a general result for the situation where you know the number $K$ of subpopulations, or at least want to fix it for the estimation, but are hesitating on the models for the emission distributions. For instance, let us consider Gaussian and Cauchy location-scale families for the composite emission families. For all $k\in\{1,\dots,K\}$, we take $\Lambda_k=\{1,2\}$ with $\overline{\mscr{F}}_1=\mscr{G}$ and $\overline{\mscr{F}}_2=\mscr{C}$, where $\mscr{C}$ is the Cauchy location-scale family of distributions associated to the density class
\begin{equation}
\label{eq:cauchy_location_scale}
\mcal{C}=\left\{ x\mapsto \frac{1}{\pi\sigma} \frac{1}{1+ \left(\frac{x-z}{\sigma} \right)^2}; z\in\mbb{R}, \sigma>0  \right\}.
\end{equation}
We consider the model $\mscr{Q}=\cup_{0\leq j\leq K} \mscr{Q}_j$ with
\begin{equation*}
\mscr{Q}_j = \left\{ \suml_{k=1}^j w_k \mcal{N}(z_k,\sigma^2_k) + \suml_{k=j(\lambda)+1}^K w_k \cauchy(z_k,\sigma_k); \begin{array}{l}
(z_1,\sigma_1)>\dots>(z_j,\sigma_j),\\
(z_{j+1},\sigma_{j+1})>\dots>(z_K,\sigma_K)
\end{array} \right\},
\end{equation*}
where the order $>$ on the parameters $(z_k,\sigma_k)$ is defined by (\ref{eq:order_parameter}) and allows to have identifiability properties again here. We consider a null penalty function.
\begin{theoreme}
\label{th:gauss_cauchy}
Assume $P^*=\suml_{k=1}^{j^*} \overline{w}_k \mcal{N}(\overline{z}_k,\overline{\sigma}^2_k) + \suml_{k=j^*+1}^K \overline{w}_k \cauchy(\overline{z}_k,\overline{\sigma}_k)\in\mscr{Q}_{j^*}$ with $(\overline{z}_1,\overline{\sigma}_1)>\dots>(\overline{z}_{l^*},\overline{\sigma}_{l^*})$ and $(\overline{z}_{l^*+1},\overline{\sigma}_{l^*+1})>\dots>(\overline{z}_K,\overline{\sigma}_K)$. Let $\hat{P}$ be a $\rho$-estimator on $\mscr{Q}_{\delta}$ with $\delta=\frac{5}{n}\bigwedge \frac{1}{K}$ and a null penalty. There exists an integer $n_0(P^*)$ and a positive constant $C(P^*)$ such that for $n\geq n_0(P^*)$ there exists an event of probability $1-(n(K+1))^{-K}$ on which such that $\hat{P}\in\mscr{Q}_{j^*}$ and 
\begin{align*}
C(P^*) \left( ||\overline{w}-\hat{w}||^2 + \suml_{k=1}^{j^*} \left|\left|(\overline{z}_k,\overline{\sigma}^2_k)- (\hat{z}_k,\hat{\sigma}^2_k)\right|\right|^2\wedge 1\right. & + \left. \suml_{k=j^*+1}^{K} \left|\left|(\overline{z}_k,\overline{\sigma}_k)- (\hat{z}_k,\hat{\sigma}_k)\right|\right|^2\wedge 1 \right)\\
&\leq \frac{K\log(n(K+1))}{n}.
\end{align*}
\end{theoreme}
This result shows that it is possible to identify the true emission models for $n$ large enough and if this identification is established we can also estimate the different parameters. This seems to be somehow original as we did not find any result of this kind in the literature.
\subsection{Selection of the order \texorpdfstring{$K$}{K}}
We consider $\Theta$ of the form $\Theta = \bigcup\limits_{K \in \mscr{K}} \{K\}\times \{ \lambda \}^K$, where $\mscr{K}$ is a subset of $\{1,\dots,n\}$. For $K \in \mscr{K}$, we write $\overline{\mscr{F}}=\overline{\mscr{F}}_{\lambda}$ and $\mscr{F}=\mscr{F}_{\lambda}$ its countable and dense subset given by Assumption \ref{hyp:countable}. For $K\in\mscr{K}$, the model $\mscr{Q}(K)$ is a subset of
\begin{equation*}
\left\{ \suml_{k=1}^K w_k F_k ; w\in\mcal{W}_K, F_k\in\overline{\mscr{F}},\forall k\in[K] \right\}.
\end{equation*}
We define $\mscr{Q}_\delta (K) := \left\{ \suml_{k=1}^K w_k F_k \in \mscr{Q}(K) ; w\in\mcal{W}_K, w_k\geq \delta, w_k \in \mbb{Q}, F_k\in\mscr{F},\forall k\in[K] \right\}$ and $\mscr{Q}_{\delta}=\bigcup_{K\in \mscr{K}} \mscr{Q}_{\delta}(K)$, where $\delta : \mscr{K} \rightarrow (0,1]$ satisfies $\delta(K)\leq 1/K$. 
Under Assumption \ref{hyp:vc_density_model}, we denote by $V$ the VC-index of $\overline{\mcal{F}}$. 
\begin{theoreme}
\label{th:selection_K}
Let $\Delta$ be a function $\mscr{K} \rightarrow \mbb{R}^+$ satisfying $\suml_{K \in \mscr{K}} e^{-\Delta(K)} \leq 1$. We consider the penalty function defined by
\begin{equation}
\label{eq:pen_selection_K}
\pen(q) = \kappa \inf\limits_{K \in\mscr{K}|Q\in\mscr{Q}(K)} \left[ 174.1 K V \left[ 5.82 + \log\left( \frac{(K+1)^2}{\delta(K)} \right) + \log_+ \left( \frac{n}{ K V} \right) \right] + \Delta(K) \right].
\end{equation}
Assume there exists $P^*$ in $\mscr{P}$ such that $\mbf{P^*}=(P^*)^{\otimes n}$. For the choice $\delta(1)=1$ and $\delta(K)=\frac{V}{n}\bigwedge \frac{1}{K}$ for $K\geq 2$, there is a positive constant $C$ such that the resulting estimator $\hat{P}=\hat{P}_{\delta}$ sat  Any $\rho$-estimator $\hat{P}_{\delta}$ on $\mscr{Q}_{\delta}$ satisfies the following. For all $\xi>0$, with probability at least $1-e^{-\xi}$ we have
\begin{equation}
\label{eq:th_order_selection_delta}
C h^2(P^*,\hat{P}) \leq \inf_{K \in \mscr{K}} \left\{ h^2( P^*, \mscr{Q}(K) ) + \frac{K V \log(n) + \xi + \Delta(K)}{n} \right\}.
\end{equation}
\end{theoreme}
This result gives an oracle inequality and it provides a way to determine the number of clusters if one wants to use mixture models in order to do clustering. It is also interesting in the context of density estimation. Once again, we take advantage of the approximation properties of GMMs to use our estimator for density estimation on a wider class. We use the approximation result proven by Maugis \& Michel \cite{maugis2012}. Let $\beta> 0$, $r = \lfloor \beta \rfloor$ and $k \in \mbb{N}$ such that $\beta \in (2k, 2k + 2]$. Let also $\mcal{P}$ be the 8-tuple of parameters $(\gamma, l^+, L, \epsilon, C, \alpha, \xi, M)$ where $L$ is a polynomial function on $\mbb{R}$ and the other parameters are positive constants. We define the density class $\mcal{H}(\beta, \mcal{P})$ of all densities $p$ satisfying the following conditions.
\begin{itemize}
\item For all $x$ and $y$ such that $|y - x| \leq \gamma$,
\begin{equation*}
(\log p )^{(r)}(x) - (\log p )^{(r)}(y) \leq r! L(x)|y - x|^{\beta-r}.
\end{equation*}
Furthermore for all $j \in \{0,\dots,r\}$,
\begin{equation*}
|(\log p )^{(j)}(0)| \leq l^+.
\end{equation*}
%$\log f$ is assumed to be locally $\beta$-Hölder: 
%
\item We have
\begin{equation*}
\max\limits_{1\leq j \leq r} \int_{\mbb{R}}{ \left| (\log p)^{(j)}(x)\right|^{\frac{2\beta+\epsilon}{j}} p(x) dx} \vee \int_{\mbb{R}} |L(x)|^{2+\frac{\epsilon}{\beta}} p (x)dx \leq C.
\end{equation*}
\item For all $x \in \mbb{R}, p (x) \leq M \psi(x)$.
\item The function $f$ is strictly positive, non-decreasing on $(-\infty, -\alpha)$ and non-increasing on
$(\alpha,\infty)$. For all $x \in [-\alpha, \alpha]$ we have $p (x) \geq \xi$ .
\end{itemize}
This class of functions can be approximated by Gaussian mixture models, the quality of the approximation depending on the regularity parameter $\beta$.
\begin{lemme}{(Lemma 6.1 in \cite{maugis2012})}\\
For $0 < \underline{\beta} < \overline{\beta}$, there exists a set of parameters $\mcal{P}(\underline{\beta},\overline{\beta})$ and a positive constant $c_{\underline{\beta},\overline{\beta}}$ such that for all $\beta\in \left[ \underline{\beta}, \overline{\beta} \right]$ and for all $p \in \mcal{H}\left(\beta,\mcal{P}(\underline{\beta},\overline{\beta})\right)$, 
\begin{equation*}
h^2 \left( P, \mcal{G}_K \right) \leq c_{\underline{\beta},\overline{\beta}} \frac{ \left( \log K \right)^{3\beta} }{ K^{2\beta} }.
\end{equation*}
\end{lemme}
We consider $\mscr{K}=\{2,\dots,n\}$, $\Delta(K) = K$ and the penalty function $\pen$ as in (\ref{eq:pen_selection_K}).
\begin{theoreme}
\label{th:maugis}
Let $\hat{P}=\hat{P}_{\delta}$ be a $\rho$-estimator on $\mscr{Q}_{\delta}$ with $\delta$ as in (\ref{eq:th_order_selection_delta}). For $0<\underline{\beta}<\overline{\beta}$, there exist a  positive constant $C_{\underline{\beta},\overline{\beta}}$ such that for any $p$ in $\mcal{H}\left(\beta,\mcal{P}(\underline{\beta},\overline{\beta})\right)$ with $\beta\in\left[\underline{\beta},\overline{\beta}\right]$, for all $\xi>0$, we have
\begin{equation*}
h^2(P^*,\hat{P}) \leq C_{\underline{\beta},\overline{\beta}} \left( \frac{(\log n)^{\frac{5\beta}{2\beta+1}}}{n^{\frac{2\beta}{2\beta+1}}} + \frac{\xi}{n}\right),
\end{equation*}
with probability at least $1-e^{-\xi}$.
\end{theoreme}
This theorem provides an upper bound on the convergence rate of our estimator of order $(\log n)^{5\beta/(4\beta+2)} n^{-\beta/(2\beta+1)}$. It is the same rate obtained in Theorem 2.9 of Maugis \& Michel \cite{maugis2012}. Therefore our estimator is minimax adaptive to the regularity $\beta$, up to a power of $\log(n)$.
\section{Acknowledgment}
The author would like to thank Yannick Baraud for his guidance in the redaction of this article.
\printbibliography
\appendix
\section{Main results}
Let $\Theta$ be a subset of
\begin{equation*}
\bigcup_{K\geq 1} \{K\} \times \prod_{i=1}^K \Lambda_i.
\end{equation*}
Let $\delta : \Theta \rightarrow (0,1]$ be such that for $\theta=(K,\lambda_1,\dots,\lambda_K)\in \Theta$, $\delta(\theta)\in(0,1/K]$. We write
\begin{equation*}
\mscr{Q}_\delta ( \theta ) = \left\{ \suml_{k=1}^K w_k F_k ; w_k \in [\delta(\theta),1]\cap \mbb{Q}, F_k\in\mscr{F}_k,\forall k\in[K], \suml_{k=1}^K w_k =1 \right\},
\end{equation*}
We define $\mscr{Q}_{\delta}$ by
\begin{equation*}
\mscr{Q}_{\delta} = \bigcup_{\theta \in \Theta} \mscr{Q}_{\delta}(\theta).
\end{equation*}
\begin{prop}
\label{prop:rho_dim_vc}
Under Assumption \ref{hyp:vc_density_model}, for $\theta=(K,\lambda_1,\dots,\lambda_K)\in\Theta$, we write
\begin{equation*}
V(\theta) = V_{1,\lambda_1} + \dots + V_{K,\lambda_K}.
\end{equation*}
For all $\mbf{P}\in \pmb{\mscr{P}}$ and $\overline{P}\in \mscr{P}$, the $\rho$-dimension admits the following bound
\begin{equation}
\label{eq:prop_rho_dim}
D^{\mscr{Q}_{\delta}(\theta)}\left(\mbf{P},\overline{P}^{\otimes n}\right) \leq D_n\left( \delta, \theta\right) = 818.1 \overline{V} \left[ 5.82 + \log\left( \frac{(K+1)^2}{\delta} \right) + \log_+ \left( \frac{n}{\overline{V}} \right) \right].
\end{equation}
\end{prop}
Let $\Delta$ be a function $\Theta \rightarrow \mbb{R}^+$ satisfying
\begin{equation}
\label{eq:Delta}
\suml_{\theta \in \Theta} e^{-\Delta(\theta)} \leq 1.
\end{equation}
We define the penalty function $\pen$ by
\begin{equation*}
\pen(q) = \kappa \inf\limits_{\theta;Q\in\mscr{Q}(\theta)} \left[ \frac{ D_n(\delta,\theta)}{4.7} + \Delta(\theta) \right],
\end{equation*}
where $D_n(\theta)$ is given by (\ref{eq:prop_rho_dim}). We set
\begin{equation*}
\mbf{\Upsilon}(\mbf{X},q) = \sup_{q'\in\mcal{Q}} \left[ \mbf{T}(\mbf{X},q,q') - \pen(q') \right] + \pen(q).
\end{equation*}
\begin{theoreme}
\label{th:general}
Any $\rho$-estimator $\hat{P}$ on $\mscr{Q}_{\delta}$ satisfies, with probability at least $1-e^{-\xi}$,
\begin{align}
\label{eq:th_general}
\mbf{h}^2\left( \mbf{P^*}, \hat{P}^{\otimes n} \right) \leq \inf_{\theta\in\Theta} &\left[ c_0 \left( \mbf{h}^2( \mbf{P^*}, \mscr{Q}(\theta) ) + n (K-1) \delta(\theta) \right) \right.\\
&+ \left. c_2 \left( 174.1 \overline{V} \left[ 5.82 + \log\left( \frac{(K+1)^2}{\delta} \right) + \log_+ \left( \frac{n}{\overline{V}} \right) \right] + \Delta(\theta) \right) \right]\nonumber\\
&+ c_2 (1.49 + \xi).\nonumber
\end{align}
with $c_0=300$ and $c_2=5014$ (see Baraud \& Chen \cite{baraudglm} for constants).
\end{theoreme}
\subsection*{Proof of Theorem \ref{th:single_model}}
It is a direct application of Theorem \ref{th:general} in the specific situation where $\Theta=\{\theta=(K,\lambda_1,\lambda_2,\dots,\lambda_K)\}$. Then, taking $\Delta(\theta)=0$, inequality (\ref{eq:th_general}) becomes
\begin{align*}
\mbf{h}^2\left( \mbf{P^*}, \hat{P}^{\otimes n} \right) &\leq c_0 \left( \mbf{h}^2( \mbf{P^*}, \mscr{Q} ) +  n(K-1) \delta\right)\\
&+ c_2 174.1 \overline{V} \left[ 5.82 + \log\left( \frac{(K+1)^2}{\delta} \right) + \log_+ \left( \frac{n}{\overline{V}} \right) \right]\\
&+ c_2 (1.49 + \xi).
\end{align*}
Take $\delta = \frac{\overline{V}}{n (K-1)} \wedge \frac{1}{K}$ (considering $K\geq 2$ here). We have $c_1=8.8\times 10^5\geq 5014 \times 174.1$.
\begin{itemize}
\item If $\overline{V} \leq n(K-1)/K$, then
\begin{align*}
\log\left( \frac{(K+1)^2}{\delta} \right) + \log_+ \left( \frac{n}{\overline{V}} \right) &= \log \left( \frac{(K^2-1)  (K+1) n^2}{\overline{V}^2} \right)\\
&\leq 3 \log\left( \frac{K n}{\overline{V}} \right) + \log\left( \frac{(K^2-1)(K+1) \overline{V}}{K^3 n} \right)\\
&\leq 3 \log\left( \frac{K n}{\overline{V}} \right) + \log\left( \frac{(K^2-1)^2}{K^4} \right)\\
&\leq 3 \log\left( \frac{K n}{\overline{V}\wedge n} \right).
\end{align*}
\item Otherwise $\overline{V} > n(K-1)/K$ and
\begin{align*}
\log\left( \frac{(K+1)^2}{\delta} \right) + \log_+ \left( \frac{n}{\overline{V}} \right) &\leq \log\left( \frac{(K+1)^2 K^2}{K-1} \right)\\
&\leq 3 \log(K) + \log\left( \frac{K^2 + 2K + 1}{K(K-1)} \right)\\
&\leq 3 \log(K) + \log\left( 9/2 \right)\\
&\leq \left[ 2 + \frac{3\log(3)}{\log(2)} \right] \log\left( \frac{K n}{\overline{V} \wedge n} \right).
\end{align*}
\end{itemize}
Finally,
\begin{align*}
\mbf{h}^2\left( \mbf{P^*}, \hat{P}^{\otimes n} \right) &\leq c_0 \left( \mbf{h}^2( \mbf{P^*}, \mscr{Q} ) +  n(K-1) \delta\right)\\
&+ c_2 174.1 \overline{V} \left( 2 + 3 \log_2(3) \right) \left[ 5.82 + \log \left( \frac{ K n}{ \overline{V} \wedge n} \right) \right]\\
&+ c_2 (1.49 + \xi).
\end{align*}
We have $c_2 174.1(2+3\log_2(3))\leq 4.52\times 10^6=c'_1$. One can easily check that it still holds for $K=1$ (see \cite{baraudrevisited}).
\subsection*{Proof of Theorem \ref{th:totally_bounded}}
Let $\mscr{Q}_K[\epsilon]$ be the model defined by
\begin{equation*}
\mscr{Q}_K[\epsilon] = \left\{ \suml_{k=1}^K w_k F_k ; w\in\mcal{W}_K, F_k\in\mscr{F}_k[\epsilon], \forall k \in\{1,\dots,K\} \right\}.
\end{equation*}
Since the class $\overline{\mscr{F}}_k$ is totally bounded, the set $\mscr{F}_k[\epsilon]$ is finite for all $k\in\{1,\dots, K\}$. We now satisfy Assumptions \ref{hyp:countable} and \ref{hyp:vc_density_model} and can apply Theorem \ref{th:single_model} with $\overline{V}=\suml_{k=1}^K \log_2(|\mscr{F}_k[\epsilon]|)\leq \suml_{k=1}^K \left( \frac{A_k}{\epsilon}\right)^{\alpha_k}$. Let $\hat{P}=\hat{P}_{\delta}$ be a $\rho$-estimator on $\mscr{Q}_{K,\delta}[\epsilon]$. For all $\xi>0$, we have
\begin{align*}
\mbf{h}^2 \left(\mbf{P^*}, \left(\hat{P}_{\delta}\right)^{\otimes n}\right) &\leq c_0 \left[ \mbf{h}^2\left(\mbf{P^*},\mscr{Q}_K[\epsilon]\right) + n(K-1)\delta \right]\\
&+ c_1 \overline{V} \left[ 5.82 + \log\left( \frac{(K+1)^2}{\delta} \right) + \log_+ \left( \frac{n}{\overline{V}} \right) \right]\\
&+ c_2 (1.49+\xi),
\end{align*}
with probability at least $1-e^{-\xi}$.
\begin{lemme}
\label{lem:upper_bound_parameters}
Let $w$ and $v$ be in $\mcal{W}_K$. Let $F_k$ and $G_k$ be in $\mscr{P}$ for all $k\in\{1,\dots,K\}$. We have
\begin{equation*}
h\left( \suml_{k=1}^K w_k F_k, \suml_{k=1}^K v_k G_k \right) \leq h( w, v) + \max_{k\in[K]} h\left( F_k, G_k \right).
\end{equation*}
\end{lemme}
This lemma implies that $\mscr{Q}_K[\epsilon]$ is a $\epsilon$-net of $\mscr{Q}_K$ with respect to the Hellinger distance, and in particular $\mbf{h}^2\left( \mbf{P^*},\mscr{Q}_K[\epsilon] \right) \leq 2\mbf{h}^2\left( \mbf{P^*},\mscr{Q}_K \right) + 2 n \epsilon^2$. For the choice $\delta=\frac{\overline{V}}{n(K-1)}\wedge \frac{1}{K}$ and $\epsilon = n^{-\frac{1}{\alpha_{\infty}+2}}$, there exists a positive constant $C$ such that for all $\xi>0$, we have
\begin{equation*}
C h^2\left( P^*, \hat{P} \right) \leq h^2\left( P^*, \mscr{Q}_K \right) + n^{-\frac{2}{\alpha_{\infty}+2}} \left( 1 + \suml_{k=1}^K A_k^{\alpha_k} \right) \left[ 1 + \log\left( Kn \right) \right] + \xi,
\end{equation*}
with probability at least $1-e^{-\xi}$.
\subsection*{Proof of Theorem \ref{th:selection_emission_models}}
Applying Theorem \ref{th:general}, we get
\begin{align*}
h^2\left( P^*, \hat{P} \right) \leq \inf_{\lambda\in L} &\bigg[ c_0 \left( h^2( P^*, \mscr{Q}(\lambda) ) + (K-1) \delta(\lambda) \right)\\
+ c_2 \bigg\{ \frac{174.1 \overline{V}(\lambda)}{n} &\left[ 5.82 + \log\left( \frac{(K+1)^2}{\delta(\lambda)} \right) + \log_+ \left( \frac{n}{\overline{V}(\lambda)} \right) \right] + \Delta(\lambda) \bigg\} \bigg]\\
&+ c_2 \frac{1.49 + \xi}{n},
\end{align*}
with probability at least $1-e^{-\xi}$. We deduce from that for the choice $\delta(\lambda) =\frac{\overline{V}(\lambda)}{n(K-1)}\bigwedge \frac{1}{K}$, there is a numeric constant $C>0$ such that, for all $\xi>0$, we have
\begin{equation*}
C h^2\left( P^*, \hat{P} \right) \leq \inf_{\lambda\in L} \bigg[ h^2( P^*, \mscr{Q}(\lambda) ) + \frac{1}{n} \bigg\{ \overline{V}(\lambda) \left[ 1 + \log\left(  \frac{Kn}{\overline{V}(\lambda)\wedge n} \right) \right] + \Delta(\lambda) \bigg\} \bigg] + \frac{\xi}{n},
\end{equation*}
with probability at least $1-e^{-\xi}$.
\subsection*{Proof of Theorem \ref{th:selection_K}}
Applying Theorem \ref{th:general}, we get
\begin{align*}
h^2\left( P^*, \hat{P} \right) \leq \inf_{K\in\mscr{K}} &\bigg[ c_0 \left( h^2( P^*, \mscr{Q}(K) ) + (K-1) \delta(K) \right)\\
+ c_2 \bigg\{ \frac{174.1 K V}{n} &\left[ 5.82 + \log\left( \frac{(K+1)^2}{\delta(K)} \right) + \log_+ \left( \frac{n}{K V} \right) \right] + \Delta(K) \bigg\} \bigg]\\
&+ c_2 \frac{1.49 + \xi}{n},
\end{align*}
with probability at least $1-e^{-\xi}$. We deduce from that for $\delta(1)=1$ and $\delta(K)=\frac{V}{n}\bigwedge \frac{1}{K}$ for $K\geq 2$, there is a numeric constant $C>0$ such that, for all $\xi>0$, we have
\begin{equation*}
C h^2\left( P^*, \hat{P} \right) \leq \inf_{K\in\mscr{K}} \bigg[ h^2( P^*, \mscr{Q}(K) ) + \frac{1}{n} \left\{ KV \left[ 1 + \log\left( \frac{K n}{KV\wedge n} \right) \right] + \Delta(K) \right\} \bigg] + \frac{\xi}{n},
\end{equation*}
with probability at least $1-e^{-\xi}$.
\subsection*{Proof of Theorem \ref{th:general}}
We recall that the function $\psi$ defined by (\ref{eq:psi_rho}) satisfies Assumption 2 of Baraud and Birgé \cite{baraudrevisited} with $a_0 = 4, a_1 = 3/8$ and $a_2^2 = 3\sqrt{2}$ (see Proposition 3 \cite{baraudrevisited}). Using Proposition \ref{prop:rho_dim_vc}, we can apply Theorem 2 \cite{baraudrevisited} with 
\begin{equation*}
D_n(\theta) = 818.1 \overline{V} \left[ 5.82 + \log\left( \frac{(K+1)^2}{\delta} \right) + \log_+ \left( \frac{n}{\overline{V}} \right) \right].
\end{equation*}
There exist constants $\gamma$ and $\kappa$ (given by (19) in \cite{baraudrevisited}) such that, with probability $\geq 1-e^{-\xi}$, we have
\begin{align*}
\mbf{h}^2\left( \mbf{P^*}, \mbf{\hat{P}} \right) &\leq \inf_{\theta\in\Theta} \left[ \gamma \mbf{h}^2( \mbf{P^*}, \mscr{Q}_\delta(\theta) ) + \frac{4 \kappa}{a_1} \left( \frac{D_n(\theta)}{4.7} + \Delta(\theta) \right) \right]\\
&+ \frac{4\kappa}{a_1} (1.49+\xi).
\end{align*}
\begin{lemme}
\label{lem:dist_m_k_delta}
For $K\geq 2$, $\delta\in [0,1/K]$,
\begin{equation}
\label{eq:approx_delta}
\forall P\in\mscr{P}, h^2(P,\mscr{Q}_{\delta}) \leq (K-1)\delta + h^2(P,\mscr{Q}).
\end{equation}
\end{lemme}
Using this inequality, we get
\begin{align*}
\mbf{h}^2\left( \mbf{P^*}, \mbf{\hat{P}} \right) &\leq \inf_{\theta\in\Theta} \left[ 2 \gamma \left( \mbf{h}^2( \mbf{P^*}, \mscr{Q}(\theta) ) + n (K(\theta)-1) \delta(\theta) \right)\right.\\
&+ \left. \frac{4 \kappa}{a_1} \left( 174.1 \overline{V} \left[ 5.82 + \log\left( \frac{(K+1)^2}{\delta} \right) + \log_+ \left( \frac{n}{\overline{V}} \right) \right] + \Delta(\theta) \right) \right]\\
&+ \frac{4\kappa}{a_1} (1.49+\xi).
\end{align*}
From Baraud \& Chen \cite{baraudglm}, we get that $\gamma<150$ and $4\kappa/a_1<5014$.
\subsection*{Proof of Proposition \ref{prop:rho_dim_vc}}
By (34) in Baraud \& Birgé \cite{baraudrevisited}, we always have the bound $D^{\mscr{Q}_{\delta(\theta)}}\left(\mbf{P^*},\overline{P}^{\otimes n}\right) \leq n/6$. We follow the proofs of Proposition 7 \cite{baraudrevisited} and Theorem 12 \cite{baraudinventiones}. T he key point here is to prove that we still have uniform entropy when we consider mixtures of VC-subgraph classes. For a metric space $(\mscr{A},d)$ and $\epsilon >0$, we denote by $N(\epsilon,\mscr{A},d)$ the minimal number of balls of radius $\epsilon$ needed to cover $\mscr{A}$. The next lemma is an intermediate result in the proof of Theorem 2 \cite{baraudglm}. 
\begin{lemme}
\label{lem:aux_rho_dim}
Let $\mcal{F}$ be a set of measurable functions $\mscr{X}\rightarrow [-1,1]$ such that for any product probability distribution $\mbf{P}=P_1\otimes \dots\otimes P_n$, we have
\begin{equation*}
\log(N(\epsilon,\mcal{F},||\cdot||_{2,P}))\leq a + b\log(1/\epsilon).
\end{equation*}
We define $Z(\mcal{F})$ by
\begin{equation*}
Z(\mcal{F}) = \sup_{f\in\mcal{F}} \left| \suml_{i=1}^n \left( f(X_i)-\mbb{E}\left[ f(X_i) \right] \right) \right|
\end{equation*}
and assume $\sup_{f\in\mcal{F}} \frac{1}{n} \suml_{i=1}^n \mbb{E}\left[ f^2(X_i)\right] \leq \sigma^2 \leq 1$. Let $q\in(0,1)$. We have
\begin{align*}
\mbb{E} \left[ Z(\mcal{F}) \right] \leq \mbb{E} \leq 32 A^2 + A 2 \sqrt{2n\sigma^2},
\end{align*}
with $A=\frac{1+q}{1-q}\left( 1 + \frac{b}{\log 2 + 2 a + b\log(1/q)} \right) \sqrt{\log 2 + 2 a + b\log(1/q) + 2b\log(1/\sigma)}$.
\end{lemme}
We now need a bound on the covering number in our situation. We use the following lemma.
\begin{lemme}
\label{lem:unif_entropy_vc_mix}
Let $\overline{P}\in \mscr{P}$ fixed. We define
\begin{equation*}
\mcal{F}^{\mcal{Q}_{\delta}}\left(\overline{P}\right) := \left\{ \psi\left(\sqrt{\frac{q}{\overline{p}}}\right) ; Q\in\mscr{Q}_{\delta}\right\}.
\end{equation*}
For any probability distribution $R$, we have
\begin{equation}
\label{eq:entropy_mixture}
\forall \epsilon \leq 2, \log N\left(\epsilon,\mcal{F}^{\mscr{Q}_{\delta}}, ||\cdot||_{2,R} \right) \leq \overline{V} \log\left( \frac{ e^{1+1/e}  8 (K+1)^2 }{ \delta } \right) + 2\overline{V} \log( 1/\epsilon).
\end{equation}
\end{lemme}
Let $y$ be a positive real number. We set
\begin{equation*}
\mcal{F} = \left\{ \psi\left(\sqrt{\frac{q}{\overline{p}}}\right) ; Q\in\mscr{Q}_{\delta}, \mbf{h}^2(\mbf{P^*},\mbf{\overline{P}}) + \mbf{h}^2(\mbf{P^*},\mbf{\hat{P}}) < y^2 \right\} \subset \mcal{F}^{\mcal{Q}_{\delta}}(\overline{P}).
\end{equation*}
Since $\psi$ satisfies Assumption 2 \cite{baraudrevisited}, we can apply Lemma \ref{lem:aux_rho_dim} with 
$\sigma^2 = (a_2^2 y^2/n) \wedge 1$, $a=\overline{V} \log\left( \frac{ e^{1+1/e}  8 (K+1)^2 }{ \delta } \right)$ and $b=2\overline{V}$. Following notation from Baraud \& Birgé \cite{baraudrevisited}, we get
\begin{align*}
\mbf{w}^{\mscr{Q}_{\delta}}\left(\mbf{P},\mbf{\overline{P}},y\right) \leq \mbb{E} \left[ Z(\mcal{F}) \right] &\leq 32 A^2 + A 2 \sqrt{2n\sigma^2}.
\end{align*}
We have
\begin{align*}
\frac{b}{\log 2 + 2 a + b\log(1/q)} &= \frac{2\overline{V}}{\log 2 + 2 \overline{V} \log\left( \frac{ e^{1+1/e}  8 (K+1)^2}{ \delta } \right) + 2\overline{V} \log(1/q)}\\
&\leq \frac{1}{\log\left( \frac{ e^{1+1/e}  8 (K+1)^2}{ \delta q} \right)}\\
&\leq \frac{1}{\log\left( \frac{e^{1+1/e} 8 K(K+1)^2}{q} \right)} \leq \frac{1}{\log\left( \frac{e^{1+1/e} 2^4\times 3^2}{q} \right)},
\end{align*}
hence
\begin{align*}
A \leq \frac{1+q}{1-q}\left( 1 + \frac{1}{\log\left( \frac{ e^{1+1/e}  2^4 \times 3^2}{q} \right)} \right) \sqrt{ 2 \overline{V} \left[ \log\left( \frac{ e^{1+1/e}  2^{13/4}}{ q } \right) + \log\left( \frac{ (K+1)^2}{\delta \sigma^2}\right) \right] }.
\end{align*}
For $q=1/9$, we have
\begin{align*}
A &\leq \frac{5}{4} \left( 1 + \frac{1}{1+\frac{1}{e} + 4 \log(6)} \right) \sqrt{ 2 \overline{V} \left[ \frac{1}{e} + 1 + \log\left(  2^{13/4}\times 9 \right) + \log\left( \frac{(K+1)^2}{\delta \sigma^2}\right) \right] }\\
&\leq \frac{5}{4} \times 1.12 \sqrt{ 2 \overline{V} \left[ 5.82 + \log\left( \frac{(K+1)^2}{\delta \sigma^2}\right) \right] }\\
&= 2.8 \sqrt{ 2 \overline{V} \left[ 5.82 + \log\left( \frac{ (K+1)^2 }{ \delta \sigma^2}\right) \right] }.
\end{align*}
Finally,
\begin{equation}
\label{eq:aux_rho_dim}
\mbf{w}^{\mscr{Q}_{\delta}}(\mbf{P},\mbf{\overline{P}},y) \leq C_0 \sqrt{n \overline{V} \sigma^2 \mscr{L}(\sigma,K,\delta)} + C_1 \overline{V} \mscr{L}(\sigma,K,\delta)
\end{equation}
with $\mscr{L}(\sigma,K,\delta) = 5.82 + \log\left(\frac{(K+1)^2}{\delta \sigma^2} \right)$, $C_0= 2.8\times 4=11.2$ and $C_1=2^6 \times 2.8^2$. Then we follow the proof of Proposition 6 \cite{baraudglm}. For $D\geq 2^{-11} \overline{V}$ and $y\geq \beta^{-1} \sqrt{D}$,
\begin{align*}
\mscr{L}(\sigma,K,\delta) &= 5.82 + \log\left( \frac{(K+1)^2}{\delta} \right) + \log_+ \left( \frac{n}{a_2^2 y^2} \right)\\
&\leq 5.82 + \log\left( \frac{(K+1)^2}{\delta} \right) + \log_+ \left( \frac{a_1^2 n}{ 16 a_2^4 D} \right)\\
&= 5.82 + \log\left( \frac{(K+1)^2}{\delta} \right) + \log_+ \left( \frac{n}{ 2^{11} D} \right)\\
&\leq 5.82 + \log\left( \frac{(K+1)^2}{\delta} \right) + \log_+ \left( \frac{n}{\overline{V}} \right) = L.
\end{align*}
We combine it with (\ref{eq:aux_rho_dim}) and we get
\begin{align*}
\mbf{w}^{\mscr{Q}_{\delta}}(\mbf{P},\mbf{\overline{P}},y) &\leq 11.2  a_2 y \sqrt{\overline{V} L} + 2^6 \times 2.8^2 \overline{V} L\\
&= \frac{a_1 y^2}{8} \left[ \frac{ 8\times 11.2 a_2 \sqrt{\overline{V} L}}{a_1 y} + \frac{2^9 \times 2.8^2 \overline{V} L}{a_1 y^2} \right]\\
&\leq \frac{a_1 y^2}{8} \left[ \frac{ 8\times 11.2 a_2 \sqrt{\overline{V} L}}{a_1 \beta^{-1} \sqrt{D}} + \frac{2^9 \times 2.8^2 \overline{V} L}{a_1 \beta^{-2} D} \right]\\
&= \frac{a_1 y^2}{8} \left[ 2 \times 11.2 \frac{ \sqrt{\overline{V} L} }{\sqrt{D}} + 2^9 \times 2.8^2 \frac{ a_1 \overline{V} L}{16 a_2^2 D} \right]\\
&= \frac{a_1 y^2}{8} \left[ 22.4 \frac{ \sqrt{\overline{V} L} }{\sqrt{D}} + 2^{9/2} \times 2.8^2 \frac{\overline{V} L}{D} \right]\\
&\leq \frac{a_1 y^2}{8} \left[ 22.4 \frac{ \sqrt{\overline{V} L} }{\sqrt{D}} + 177.4 \frac{\overline{V} L}{D} \right].
\end{align*}
One can check that for $D = 818.1 \overline{V} L  \geq \overline{V} L \left[ \sqrt{ 177.4 + 22.4^2/4 } + 22.4/2 \right]^2$ we have $\mbf{w}^{\mscr{Q}_{\delta}}(\mbf{P},\mbf{\overline{P}},y) \leq \frac{a_1 y^2}{8}$. Since $L\geq 5.82$, we also have $D\geq 2^{-11} \overline{V}$.
\section{Density estimation}
\subsection*{Proof of Theorem \ref{th:gmm}}
The Gaussian location-scale family of density functions is VC-subgraph (see Lemma \ref{lem:vc}). Proposition \ref{prop:ghosal_approximation} provides an approximation bound for $\mscr{C}(A,R)$.
We can now apply Theorem \ref{th:single_model} with those two propositions. With (\ref{eq:th1_2}), there exists a universal constant $C$ such that for $\mbf{P^*}=(P^*)^{\otimes n}$, $\xi>0$, with probability at least $1-e^{-\xi}$, we have
\begin{align*}
C h^2 \left(P^*,\hat{P}\right) &\leq h^2\left(P^*,\mscr{C}(A,R)\right) + \exp\left( - K^{1/2} \frac{3 \sqrt{3/2}}{R^2} \right) \left[ K^{1/4} \frac{3 \sqrt{2}}{\sqrt{e\pi} 7^{1/4}} + R \right]\\
&+ \frac{K \log\left( n \right) + \xi}{n}\\
C' h^2(P^*\hat{P}) &\leq h^2\left(P^*,\mscr{C}(A,R)\right) + \frac{R}{n} \left[ 2 ^{3/4} (3/7)^{1/4} \sqrt{\log(n)} \frac{1}{\sqrt{e\pi} } + 1 \right]\\
&+ \frac{R^4 \log\left( n \right) 2/27 + \xi}{n}.
\end{align*}
Finally, there exists a numeric constant $C>0$ such that, for $n \geq K= \left\lceil \frac{2 R^4 \log^2(n)}{27} \right\rceil \geq (2/3)^3 A^4$, for all $\xi>0$, with probability at least $1-e^{-\xi}$, we have
\begin{equation*}
C h^2 \left( P^*, \hat{P} \right) \leq h^2(P^*,\mscr{C}(A,R)) + \frac{R^4 \log(n)+\xi}{n}. 
\end{equation*}
The different conditions are satisfied for $n\geq \exp(2 (A/R)^2)$ and $\frac{n}{\log^2(n)} \geq 2 R^2 /27$.
\subsection*{Proof of Theorem \ref{th:maugis}}
The Gaussian location-scale family of density functions is VC-subgraph (see Lemma \ref{lem:vc}). For $0 < \underline{\beta} < \overline{\beta}$ and $\beta \in[\underline{\beta},\overline{\beta}]$, let $\mcal{H} \left( \beta, \mcal{P}(\underline{\beta},\overline{\beta}) \right)$ be the class of density functions defined in Maugis \& Michel \cite{maugis2012}. One can check that
\begin{equation*}
\suml_{k\in\mscr{K}} e^{-\Delta(K)} \leq 1,
\end{equation*}
for $\Delta(K)=K$. Applying Theorem \ref{th:selection_K}, for $\xi>0$, with probability at least $1-e^{-\xi}$, we have
\begin{align*}
C h^2(P^*,\hat{P}) &\leq \inf_{K \in \mscr{K}} \left\{ h^2( P^*, \mscr{G}_K ) + \frac{K (5 \log(n) + 1) + \xi}{n} \right\}\\
&\leq 2 h^2( P^*, \mcal{H}\left(\beta,\mcal{P}(\underline{\beta},\overline{\beta})\right) ) + \inf_{K \in \mscr{K}} \left\{ 2 c_{\underline{\beta},\overline{\beta}} \frac{(\log K)^{3\beta}}{K^{2\beta}} + \frac{K (5 \log(n) + 1)}{n} \right\} + \frac{\xi}{n}.
\end{align*}
Therefore, following the proof of Theorem 2.9 of Maugis \& Michel \cite{maugis2012}, we have
\begin{align*}
\inf_{K \in \mscr{K}} \left\{ 2 c_{\underline{\beta},\overline{\beta}} \frac{(\log K)^{3\beta}}{K^{2\beta}} + \frac{K (5 \log(n) + 1)}{n} \right\} &\lesssim  c_{\underline{\beta},\overline{\beta}}\inf_{K\in\mscr{K}} \left\{ \frac{(\log K)^{3\beta}}{K^{2\beta}} + \frac{K \log(n)}{n} \right\}\\
&\lesssim c_{\underline{\beta},\overline{\beta}} \frac{(\log n)^{\frac{5\beta}{2\beta+1}}}{n^{\frac{2\beta}{2\beta+1}}}.
\end{align*}
Finally, there exists $C_{\underline{\beta},\overline{\beta}}$ such that for all $\xi>0$, with probability at least $1-e^{-\xi}$, we have
\begin{equation*}
h^2(P^*,\hat{P}) \leq C_{\underline{\beta},\overline{\beta}} \left( \frac{(\log n)^{\frac{5\beta}{2\beta+1}}}{n^{\frac{2\beta}{2\beta+1}}} + \frac{\xi}{n}\right).
\end{equation*}
\subsection*{Proof of Proposition \ref{prop:ghosal_approximation}}
\begin{lemme}
\label{lem:ghosal_dtv}
Let $k$ be a positive integer. For any probability distribution $H$ on $[-a,a]\times [\underline{\sigma},\overline{\sigma}]$, there is a discrete probability distribution $H'$ supported by $k(2k-1)+1$ points in $[-a,a]\times [\underline{\sigma},\overline{\sigma}]$ such that
\begin{equation*}
d_{TV}\left( P_H, P_{H'} \right) \leq \inf\limits_{m>1} \left\{ \frac{\sqrt{2/\pi}}{\underline{\sigma}} am  \left( \frac{ e a^2 (1+m)^2}{2 k \underline{\sigma}^2} \right)^k + \frac{\overline{\sigma}}{2 \underline{\sigma}} \exp\left( - \frac{(m-1)^2 a^2}{2\overline{\sigma}^2}\right) \right\}.
\end{equation*}
\end{lemme}
Let $A$ and $R$ be two real numbers respectively greater than 0 and 1. As a direct consequence of this lemma, for any $l\in\mbb{R}$, any probability distribution $H$ on $[l \pm \underline{\sigma}A] \times [\underline{\sigma},R\underline{\sigma}]$ and for $K\geq k(2k-1)+1$, we have
\begin{equation*}
h^2(P_H, \mscr{G}_K) \leq \inf\limits_{m>0} \left\{ \sqrt{2/\pi} A (1+m)  \left( \frac{ e  A^2 (2+m)^2}{2 k} \right)^k + \frac{R}{2} \exp\left( - \frac{m^2 A^2}{2 R^2}\right) \right\}. 
\end{equation*}
Now
\begin{align*}
h^2(P_H, \mscr{G}_K) &\leq \inf\limits_{m \geq 2} \left\{ \sqrt{2/\pi} A \frac{3}{2}m  \left( \frac{ e  A^2 4m^2}{2 k} \right)^k + \frac{R}{2} \exp\left( - \frac{m^2 A^2}{2 R^2}\right) \right\}\\
&= \inf\limits_{m \geq 2} \left\{ \frac{3}{\sqrt{2\pi}} A m  \left( \frac{ 2e  A^2 m^2}{ k} \right)^k + \frac{R}{2} \exp\left( - \frac{m^2 A^2}{2 R^2}\right) \right\}. 
\end{align*}
For $m = \frac{\sqrt{2W(1/4eR^2)}R}{A}k^{1/2}$ and $k\geq \frac{2 A^2}{W(1/4eR^2)R^2}$, we get
\begin{align*}
h^2(P_H, \mscr{G}_K) &\leq \frac{3}{\sqrt{2\pi}} \sqrt{2W(1/4eR^2)}R k^{1/2}  \left( 4 e R^2 W(1/4eR^2) \right)^k + \frac{R}{2} \exp\left( - k W(1/4eR^2)\right)\\
&= R \exp\left( - k W(1/4eR^2) \right) \left[ k^{1/2} 3\sqrt{W(1/4eR^2)/\pi} + 1/2\right].
\end{align*}
\begin{itemize}
\item For all $x> 0$, $0<W(x) <x$.
\item For all $x\in(0,1), x(1-x)<W(x)$. Therefore,
\begin{align*}
W(1/4eR^2) &\geq \frac{1}{4eR^2} \left(1-\frac{1}{4eR^2} \right)\\
&\geq \frac{(1-1/4e)}{4eR^2} = \frac{4e-1}{R^2} \geq 9/R^2.
\end{align*}
\end{itemize}
Finally,
\begin{equation*}
h^2(P_H, \mscr{G}_K) \leq R \exp\left( - 9 k /R^2 \right) \left[ k^{1/2} \frac{3}{2 R \sqrt{e\pi}} + 1/2\right].
\end{equation*}
If $k$ is the largest integer such that $k(2k-1)+1\leq K$, i.e. $k = \left\lfloor \frac{1}{4} + \sqrt{(K-7/8)/2} \right\rfloor$, we have
\begin{equation*}
K \in \{ n(2n-1)+1, \dots, (2n+1)(n+1)\} \Rightarrow k=n\geq \sqrt{K} \frac{n}{\sqrt{(2n+1)(n+1)}}
\end{equation*}
and $k\leq 2\sqrt{K/7}$. Since $\frac{x}{\sqrt{(2x+1)(x+1)}}$ is non-decreasing on $[1,+\infty)$, we have $k \geq \sqrt{K} /\sqrt{6}$ for all $K\geq 2$. Finally, we have
\begin{align*}
h^2(P_H, \mscr{G}_K) &\leq R \exp\left( - K^{1/2} \frac{3 \sqrt{3}}{\sqrt{2} R^2} \right) \left[ K^{1/4} \frac{3 \sqrt{2}}{2 R \sqrt{e\pi} 7^{1/4}} + 1/2\right]\\
&= \frac{1}{2} \exp\left( - K^{1/2} \frac{3 \sqrt{3}}{\sqrt{2} R^2} \right) \left[ K^{1/4} \frac{3 \sqrt{2}}{ \sqrt{e\pi} 7^{1/4}} + R \right].
\end{align*}
One can see that $\underline{\sigma}$ does not play a role here and is equivalent to $s$ in the definition of $\mcal{C}(A,R)$. 
\section{Regular parametric models}
\subsection*{Proof of Theorem \ref{th:regular_parametric}}
We apply the results of Ibragimov and Has’minski\u{\i} \cite{Ibragimov} (Chapter 1, Section 7.1 and 7.3) to parametric mixture models. We recall the notation
\begin{equation*}
p(\cdot;\theta) =\suml_{k=1}^{K-1} w_k f_k(\cdot;\alpha_k) + (1-w_1-\dots-w_{K-1})f_K(\cdot;\alpha_K).
\end{equation*}
and $\Theta = \left\{ w \in (0,1)^{K-1}, \suml_{k=1}^{K-1} w_k < 1 \right\} \times A_1 \times \dots \times A_K$. Obviously, $\Theta$ is an open convex subset of $\mbb{R}^{K-1} \times \mbb{R}^d_1 \times \dots\times \mbb{R}^{d_K}$. We first check that Assumption \ref{hyp:regular} implies that the model is regular.
\begin{itemize}
\item{a)} $\Rightarrow \theta \mapsto p(x;\theta)$ is continuous on $\Theta$ for $\mu$-almost all $x\in\mscr{X}$.
\item{b)} $\Rightarrow$ For $\mu$-almost all $x\in\mscr{X}$ the function $u\mapsto p(x;u)$ is differentiable at the point $u=\theta$. For all $k\in\{1,\dots,K\}$ and $j\in\{1,\dots,d_k\}$, we have
\begin{align*}
\int_{\mscr{X}} \left| \frac{\partial p(x;\theta)}{\partial \alpha_{k,j}} \right|^2 \frac{\mu(dx)}{p(x;\theta)} &=  \int_{\mscr{X}} \left| \frac{\partial f_k(x;\alpha_k)}{\partial \alpha_{k,j}} \right|^2 \frac{w_k^2 }{p(x;\theta)} \mu(dx)\\
&\leq \int_{\mscr{X}} \left| \frac{\partial f_k(x;\alpha_k)}{\partial \alpha_{k,j}} \right|^2 \frac{\mu(dx)}{f_k(x;\alpha_k)} < \infty
\end{align*}
(it also works with $k=K$ since we only $\pi$ is fixed here) and for $k\in\{1,\dots,K-1\}$ we get
\begin{align*}
\int_{\mscr{X}} \left| \frac{\partial p(x;\theta)}{\partial w_k} \right|^2 \frac{\mu(dx)}{p(x;\theta)} &=  \int_{\mscr{X}} \left( f_k(x;\alpha_k) - f_K(x,\alpha_K) \right)^2 \frac{\mu(dx)}{p(x;\theta)}\\
&\leq \frac{2}{w_k} \int_{\mscr{X}} f_k^2(x;\alpha_k) \frac{\mu(dx)}{f_k(x;\alpha_k)}\\
&+ \frac{2}{1-w_1-\dots-w_k} \int_{\mscr{X}} f_K^2(x;\alpha_k) \frac{\mu(dx)}{f_K(x;\alpha_K)}\\
&= \frac{2}{w_k} + \frac{2}{1-w_1-\dots-w_k} <\infty.
\end{align*}
\end{itemize}
Therefore, we have a regular statistical experiment (see \cite{Ibragimov}). Since the Fisher's information matrix
\begin{equation*}
I\left(\overline{\theta}\right) = \int_{\mscr{X}} \frac{\partial p\left(x;\overline{\theta}\right)}{\partial \theta} \left( \frac{\partial p\left(x;\overline{\theta}\right)}{\partial \theta} \right)^T \frac{\mu(dx)}{p\left(x;\overline{\theta}\right)}
\end{equation*}
is definite positive. We can apply Theorem 7.6 of Ibragimov and Has’minski\u{\i} \cite{Ibragimov}. There exists a positive constant $c(\overline{\theta})$ such that
$\liminf\limits_{h\rightarrow 0} ||t||^{-2} h^2(P_{\overline{\theta}}, P_{\overline{\theta}+t}) \geq c(\overline{\theta})$. There exists $a>0$ such that
\begin{align*}
 \inf\limits_{\substack{||t||<a\\\overline{\theta}+t}} ||t||^{-2} h^2\left(P_{\overline{\theta}},P_{\overline{\theta}+t}\right) \geq c\left(\overline{\theta}\right)/2.
\end{align*}
Finally, there exists a positive constant $C\left(\overline{\theta}\right)= \frac{c\left(\overline{\theta}\right)}{2} \wedge \inf\limits_{ \substack{||t||\geq a\\ \theta+t \in\Theta}} h^2\left(P_{\overline{\theta}},P_{\overline{\theta}+t}\right) >0$ such that
\begin{equation*}
\forall \theta \in\Theta, \left(1+||t||^{-2}\right) h^2\left(P_{\overline{\theta}},P_{\overline{\theta}+t}\right) \geq C\left( \overline{\theta} \right).
\end{equation*}
So with probability at least $1-e^{-\xi}$ we have
\begin{align*}
\frac{1}{n} \left[ \mbf{h}^2\left(\mbf{P^*},P_{\overline{\theta}}^{\otimes n}\right) + \overline{V}\log(n)+\xi \right] \geq C h^2\left(P_{\overline{\theta}}, P_{\hat{\theta}} \right) &\geq \frac{\left|\left|\overline{\theta}-\hat{\theta}\right|\right|^2}{1+\left|\left|\overline{\theta}-\hat{\theta}\right|\right|^2} C \times C\left( \overline{\theta} \right)\\
&\geq \frac{\left|\left|\overline{\theta}-\hat{\theta}\right|\right|^2\wedge b}{1+b} C'(\overline{\theta}).
\end{align*}
Since $||\overline{w}-\hat{w}||_2^2 \leq K \suml_{k=1}^{K-1} (\overline{w}_k-\hat{w}_k)^2$ and
\begin{align*}
\suml_{k=1}^{K-1} \left( \overline{w}_k-\hat{w}_k\right)^2 +&  \suml_{k=1}^K \left[ \left|\left|\overline{\alpha}_k-\hat{\alpha}_K\right|\right|^2\wedge 1\right]\\
 &\leq   \suml_{k=1}^{K-1} \left( \overline{w}_k-\hat{w}_k\right)^2 + \left[ \suml_{k=1}^K \left|\left|\overline{\alpha}_k-\hat{\alpha}_K\right|\right|^2 \right] \wedge K\\
&\leq \left[ \suml_{k=1}^{K-1} \left( \overline{w}_k-\hat{w}_k\right)^2 +  \suml_{k=1}^K \left|\left| \overline{\alpha}_k-\hat{\alpha}_K\right|\right|^2 \right] \wedge (K+1)\\
&=  \left|\left|\overline{\theta}-\hat{\theta}\right|\right|^2 \wedge (K+1),
\end{align*}
we get
\begin{equation*}
\frac{1}{n} \left[ \mbf{h}^2\left(\mbf{P^*},P_{\overline{\theta}}^{\otimes n}\right) + \overline{V}\log(n)+\xi \right] \geq \left[ \frac{1}{K} || \overline{w}-\hat{w}||^2 + \suml_{k=1}^K \left|\left| \overline{\alpha}_k -\hat{\alpha}_k \right|\right|^2 \wedge 1 \right] \frac{C'(\overline{\theta})}{K+2},
\end{equation*}
with probability at least $1-e^{-\xi}$.
\subsection*{Proof of Theorem \ref{th:gauss_cauchy}}
Assumption \ref{hyp:vc_density_model} is satisfied with Lemma \ref{lem:vc}.
For all $j$ in $\{0,\dots,K\}$, we have $\overline{V}_j=5K$. We apply Theorem \ref{th:selection_emission_models} with $\Delta_j= \log(K+1)$ for all $j\in\{1,\dots,K\}$. This induces a constant penalty function and one can check that this does not modify the definition of $\rho$-estimators compared to a null penalty function. Therefore, the estimator can be computed with a null penalty. There exists a positive constant that does not depend on $P^*$ such that for $n\geq 5K$, any $\rho$-estimator $\hat{P}_{\delta}$ on $\mscr{Q}_{\delta}$ satisfies, with probability at least $1-e^{-\xi}$,
\begin{equation*}
C h^2(P^*,\hat{P}) \leq \frac{K \log\left(n (K+1) \right) + \xi}{n}.
\end{equation*}
The following lemma allow to prove that for $n$ large enough, the estimator $\hat{P}$ belongs to the true model $\mscr{Q}_{j^*}$ with high probability.
\begin{lemme}
\label{lem:aux_identifiability}
Let $j\in\{0,\dots,K\}$ and assume there is a sequence
\begin{equation*}
(P_n)_n = \left( \suml_{k=1}^j w_{k,n} \mcal{N}(z_{k,n},\sigma^2_{k,n}) + \suml_{k=j+1}^K w_{k,n} \cauchy(z_{k,n},\sigma_{k,n}) \right)\in\mscr{Q}_j^{\mbb{N}}
\end{equation*}
such that $\lim\limits_{n\rightarrow\infty} h(P_n,P^*)=0$. Then, $j=j^*$ and there is a subsequence $P_{\psi(n)}$ such that $\lim_n (z_{k,\psi(n)},\sigma_{k,\psi(n)})_{1\leq k\leq K} = (\overline{z}_k,\overline{\sigma}_k)_{1\leq k\leq K}$.  
\end{lemme}
This implies that $\alpha=\inf_{j\neq j^*} h\left(P^*,\mscr{Q}_j\right)>0$. For $n\geq n_0 =\inf\{n\geq 1 : C^{-1} \alpha^{-1} K < n/\log(n(K+1)) \}$ and $0<\xi< \frac{C n \alpha}{K\log(n(K+1))}$, there is an event $\Omega_{\xi,n}$ of probability $1-e^{-\xi}$ such that
\begin{equation*}
C h^2(P^*,\hat{P}) \leq \frac{K \log\left(n (K+1) \right) + \xi}{n} \text{  and  }\hat{P}\in\mscr{Q}_{j^*}.
\end{equation*}
From now, we follow the proof of Theorem \ref{th:gmm} to prove a lower bound on the Hellinger distance $h(P^*,P)$ for $P\in\mscr{Q}_{j^*}$.
\begin{lemme}
\label{lem:aux_regular_cauchy_gauss}
There exists a positive constant $\overline{a}$ such that for all $P=\suml_{k=1}^{j^*} w_k \mcal{N}(z_k,\sigma^2_k) + \suml_{j^*+1}^K w_k \cauchy(z_k,\sigma_k) \in\mscr{Q}_{j^*}$ 
\begin{align*}
h^2(P^*,P_{\theta}) \geq \overline{a} \bigg( ||w-\overline{w}||_2^2 &+ \suml_{k=1}^{j^*} ||(z_k,\sigma^2_k)-(\overline{z}_k,\overline{\sigma}_k^2)||_2^2 \wedge 1\\
&+ \suml_{k=j^*+1}^K ||(z_k,\sigma_k)-(\overline{z}_k,\overline{\sigma}_k)||_2^2 \wedge 1 \bigg).
\end{align*}
\end{lemme}
Finally, there is a constant $\overline{C}$ such that for $\xi$ and $n$, on the event $\Omega_{\xi,n}$, we have
\begin{align*}
\overline{C} \left( ||\hat{w}-\overline{w}||_2^2 + \suml_{k=1}^{j^*} ||(\hat{z}_k,\hat{\sigma}^2_k)-(\overline{z}_k,\overline{\sigma}_k^2)||_2^2 \wedge 1\right. &  + \left. \suml_{k=j^*+1}^K ||(\hat{z}_k,\hat{\sigma}_k)-(\overline{z}_k,\overline{\sigma}_k)||_2^2 \wedge 1 \right)\\
&\leq \frac{K\log(n(K+1))+\xi}{n}.
\end{align*}
\subsection*{Proof of Theorem \ref{th:parameter_recovery_gmm}}
We apply Theorem \ref{th:single_model} and Lemma \ref{lem:aux_regular_cauchy_gauss} with $j^*=K$.
\section{Two-component mixture models}
\subsection*{Proof of Theorem \ref{th:gadat}}
We take $M = ||z^*||_{\infty} + 1$ to have (\ref{eq:aux_gadat}). With Proposition \ref{prop:gadat}, there exists a positive constant $C$ (depending on $\phi$ and $M$) such that for all $z\in[-M,M]^d$, and all $\lambda \in [0,1]$, we have
\begin{equation*}
C(\phi,M) ||z^*||^2 \left( ||z||^2 \left(\lambda^*-\lambda\right)^2 + \left(\lambda^*\right)^2  \left|\left|z^*-z\right|\right|^2 \right) \leq || p_{\lambda^*,z^*} - p_{\lambda,z}||^2.
\end{equation*}
One can prove (using Proposition 2.1 in  \cite{gadat2018parameter} and $\lambda^*\neq 0$) that we have
\begin{equation}
\label{eq:aux_gadat}
\inf_{\substack{z\not\in [-M,M]^d,\\ \lambda\in[0,1]}} || p_{\lambda^*,z^*} - p_{\lambda,z}||^2 > 0.
\end{equation}
Therefore, there is a constant $C(\phi,\lambda^*,z^*)$ such that for all $z\in \mbb{R}^d$ and all $\lambda \in[0,1]$,
\begin{equation*}
C(\phi,\lambda^*,z^*) \left( \left( ||z||^2 \wedge 1\right) \left(\lambda^*-\lambda\right)^2 + \left(\lambda^*\right)^2 \left(  \left|\left|z^*-z\right|\right|^2 \wedge 1 \right) \right) \leq || p_{\lambda^*,z^*} - p_{\lambda,z}||^2.
\end{equation*}
Since $\phi$ is bounded, with inequality (\ref{eq:l2_hellinger}), there is another constant $C(\phi,\lambda^*,z^*)$ such that for all $z\in\mbb{R}^d$ and $\lambda\in[0,1]$ we have
\begin{equation*}
C(\phi,\lambda^*,z^*) \left( \left(||z||^2 \wedge 1 \right) \left(\lambda^*-\lambda\right)^2 + \left(\lambda^*\right)^2  \left( \left|\left|z^*-z\right|\right|^2 \wedge 1 \right) \right) \leq h^2\left(P_{\lambda^*,z^*},P_{\lambda,z}\right).
\end{equation*}
One can check the following
\begin{align*}
h^2(P_{\lambda^*,z^*},P_{\hat{\lambda},\hat{z}})\leq C(\phi,\lambda^*,z^*) \left(\lambda^*\right)^2 \left(||z^*||^2\wedge 1\right)/2 &\Rightarrow ||z^*-\hat{z}||^2 \wedge 1 \leq (||z^*||^2\wedge 1)/4\\
&\Rightarrow ||\hat{z}||\wedge 1 \geq \frac{||z^*||}{2} \wedge 1.
\end{align*}
We use Theorem \ref{th:single_model} for an upper bound on $h^2(P_{\lambda^*,z^*},P_{\hat{\lambda},\hat{z}})$. For $n\geq n_0(\phi,\lambda^*,z^*)$, with 
\begin{equation*}
n_0(\phi,\lambda^*,z^*) := \inf\left\{ n\geq 1+V \bigg| \frac{4(1+V) [1+\log(2n/(1+V))]}{n C (\lambda^*)^2 \left( ||z^*||^2\wedge 1\right)} \leq C(\phi,\lambda^*,z^*) \right\},
\end{equation*}
for $0<\xi\leq \xi_n=(1+V)[1+\log(2n/(1+V))]$, with probability at least $1-e^{-\xi}$ we have
\begin{align*}
C h^2 \left(P_{\lambda^*,z^*},P_{\hat{\lambda},\hat{z}}\right) &\leq  \frac{1}{n} \left\{ (1+V) \left[ 1 + \log\left( \frac{2 n}{(V+1)} \right) \right] + \xi \right\}\\
&\leq C \times C(\phi,\lambda^*,z^*) \left(\lambda^*\right)^2 \left(||z^*||^2\wedge 1\right)/2,
\end{align*}
where $C$ is the constant given in Theorem \ref{th:single_model}. Therefore, there is a new constant $C(\phi,\lambda^*,z^*)$ such that for $n\geq n_0$ and $\xi\in(0,\xi_n)$, with probability at least $1-e^{-\xi}$ we have
\begin{equation*}
C(\phi,\lambda^*,z^*) \left( \left(\lambda^*-\lambda\right)^2 +  \left( \left|\left|z^*-z\right|\right|^2 \wedge 1 \right) \right) \leq \frac{(1+V)\left[ 1+\log(2n/(1+V))\right]+\xi}{n}.
\end{equation*}
\subsection*{Proof of Theorem \ref{th:faster_rate}}
\begin{prop}
\label{prop:faster_rate}
For $\lambda^*\in(0,1]$ and $z^*\neq 0$, there is a positive constant $C(\alpha,\lambda^*,z^*)$ such that for all $z\in\mbb{R}$ and all $\lambda\in[0,1]$, we have
\begin{equation*}
h^2\left( P_{\lambda^*,z^*}, P_{\lambda,z} \right) \geq C(\alpha,z^*,\lambda^*) \left[ \left( \lambda^* \right)^{1/\alpha} \left( 1 \wedge |z-z^*|^{1-\alpha} \right) + \left( \lambda^* - \lambda \right)^2 \left( 1 \wedge |z^*| \right) \right].
\end{equation*}
\end{prop}
Since $s_{\alpha}$ is unimodal, the class of densities $\left\{ x \mapsto s_{\alpha}(x-z), z\in\mbb{R} \right\}$ is VC-subgraph with VC-dimension not larger than $10$ (see Section \ref{sec:performance_estimator}). With Theorem \ref{th:single_model} and Proposition \ref{prop:faster_rate}, there exists a positive constant $C(\alpha,\lambda^*,z^*)$ such that for all $\xi>0$, we have
\begin{equation*}
 C(\alpha,z^*,\lambda^*) \left[ 1 \wedge |\hat{z}-z^*|^{1-\alpha} + \left( \lambda^* - \hat{\lambda} \right)^2 \right] \leq \frac{\log(n) + \xi}{n},
\end{equation*}
with probability at least $1-e^{-\xi}$.
\subsection*{Proof of Proposition \ref{prop:faster_rate}}
We write
\begin{equation*}
f_z(x) = s_{\alpha}(x-z) = \frac{1-\alpha}{2|x-z|^{\alpha}}\mathbbm{1}_{|x-z|\in(0,1]}.
\end{equation*}
We define $g$ by
\begin{equation*}
g(x) = \frac{2}{1-\alpha}\left( \sqrt{ (1-\lambda^*)f_0(x) + \lambda^* f_{z^*}(x)} - \sqrt{ (1-\lambda)f_0(x) + \lambda f_z(x) } \right)^2
\end{equation*}
such that
\begin{equation*}
2 h^2\left( P_{\lambda^*,z^*}, P_{\lambda,z} \right) = \frac{1-\alpha}{2}\displaystyle\int_{-\infty}^{+\infty} g(x) dx.
\end{equation*}
\begin{lemme}
\label{lem:spike_aux_1}
Assuming $z\cdot z^*>0$ and $|z^*-z|\leq \frac{1}{(1-\alpha)^{2/\alpha}}$. There exists $C(\alpha,z^*,\lambda^*)>0$ such that
\begin{equation*}
\int g(x) dx \geq C(\alpha,z^*,\lambda^*) \left[ \left( \lambda^* \right)^{1/\alpha} \left( 1 \wedge |z-z^*|^{1-\alpha} \right) + \left( \lambda^* - \lambda \right)^2 \left( 1 \wedge |z^*| \right) \right].
\end{equation*}
\end{lemme}
\begin{lemme}
\label{lem:spike_aux_2}
For $z \cdot z^* \leq 0$, we have
\begin{equation*}
\int g(x) dx \geq \lambda^* \alpha^2 \frac{1 \bigwedge \left[(\lambda^*)^{(1-\alpha)/\alpha} (1-\alpha)^{2(1-\alpha)/\alpha} |z^*|^{1-\alpha} \right] }{1-\alpha}.
\end{equation*}
\end{lemme}
\begin{lemme}
\label{lem:spike_aux_3}
For $|z-z^*| > \frac{1}{(1-\alpha)^{2/\alpha}}$ and $z^*\cdot z>0$, we have
\begin{equation*}
\int g(x) dx = \lambda^* (1\wedge |z^*|).
\end{equation*}
\end{lemme}
Combining those three lemmas, there exists a positive constant $C(\alpha,z^*,\lambda^*)$ such that
\begin{equation*}
h^2\left( P_{\lambda^*,z^*}, P_{\lambda,z} \right) \geq C'(\alpha,z^*,\lambda^*) \left[ \left( \lambda^* \right)^{1/\alpha} \left( 1 \wedge |z-z^*|^{1-\alpha} \right) + \left( \lambda^* - \lambda \right)^2 \left( 1 \wedge |z^*| \right) \right],
\end{equation*}
for all $\lambda$ in $[0,1]$ and $z$ in $\mbb{R}$. Without loss of generality, we assume $z^* >0$ through the proof of the lemmas.
\subsection*{Proof of Lemma \ref{lem:spike_aux_1}}
Without loss of generality, we consider $z^*>0$ for now.\\
$\bullet$ For $x\in]-1,0[$, we have
\begin{align*}
g(x) &= \frac{1}{|x|^{\alpha}} \left( \sqrt{ 1-\lambda^* + \lambda^* \frac{|x|^{\alpha}}{|x-z^*|^{\alpha}} \mathbbm{1}_{|x-z^*|\in(0,1]}} - \sqrt{ 1-\lambda + \lambda \frac{|x|^{\alpha}}{|x-z|^{\alpha}} \mathbbm{1}_{|x-z|\in(0,1]}} \right)^2.
\end{align*}
If $z^*\wedge z\geq 1$ then,
\begin{align*}
g(x) &= \frac{1}{|x|^{\alpha}} \left( \sqrt{ 1-\lambda^*} - \sqrt{ 1-\lambda} \right)^2
\end{align*}
and
\begin{equation*}
\int_{-1}^0 g(x)dx \geq \left( \sqrt{ 1-\lambda^*} - \sqrt{ 1-\lambda} \right)^2 \frac{1}{1-\alpha}.
\end{equation*}
Otherwise $z^*\wedge z\in(0,1)$ then for $x\in ]-1,z^*\wedge z -1[$,
\begin{equation*}
\int_{-1}^{z^*\wedge z-1} g(x)dx \geq \left( \sqrt{ 1-\lambda^*} - \sqrt{ 1-\lambda} \right)^2 \frac{1-(1-z\wedge z^*)^{1-\alpha}}{1-\alpha}.
\end{equation*}
Finally,
\begin{equation*}
\int_{-1}^0 g(x)dx \geq \left( \sqrt{ 1-\lambda^*} - \sqrt{ 1-\lambda} \right)^2 \frac{1-(1-z\wedge z^*)_ +^{1-\alpha}}{1-\alpha}.
\end{equation*}
$\bullet$ For $x\in]z^*\vee z,z^*\vee z+1[$, we have
\begin{align*}
g(x) &= \frac{1}{|x-z^*\vee z|^{\alpha}} \left( \sqrt{ (1-\lambda^*) \frac{|x-z^*\vee z|^{\alpha}}{|x|^{\alpha}}\mathbbm{1}_{|x|\in(0,1]} + \lambda^* \frac{|x-z^*\vee z|^{\alpha}}{|x-z^*|^{\alpha}}\mathbbm{1}_{|x-z^*|\in(0,1]}}\right.\\
&- \left. \sqrt{ (1-\lambda) \frac{|x-z^*\vee z|^{\alpha} }{ |x|^{\alpha} } \mathbbm{1}_{|x|\in(0,1]} + \lambda \frac{|x-z^*\vee z|^{\alpha}}{|x-z|^{\alpha}} \mathbbm{1}_{|x-z|\in(0,1]} } \right)^2.
\end{align*}
\begin{itemize}
\item If $z < z^*$, with $V < \frac{1}{|z-z^*|}$, for $x\in]z^*,z^*+V|z-z^*|[$, we have
\begin{align*}
\bullet \frac{|x-z^*|}{|x|} &\leq V \frac{|z^*-z|}{z^*}\leq V,\\
\bullet \frac{|x-z^*|}{|x-z|} &\leq \frac{V|z^*-z|}{(1+V)|z^*-z|}\leq V.
\end{align*}
We get
\begin{align*}
\int_{z^*}^{z^*+V|z-z^*|} g(x) dx &\geq \left( \sqrt{ \lambda^*} - \sqrt{V^{\alpha}} \right)^2 \int_{z^*}^{z^*+V|z-z^*|} \frac{dx}{|x-z^*\vee z|^{\alpha}}\\
&= \left( \sqrt{ \lambda^*} - \sqrt{V^{\alpha}} \right)^2 \frac{\left( V |z^*-z| \right)^{1-\alpha}}{1-\alpha}.
\end{align*}
We take $V = (\lambda^*)^{1/\alpha} (1-\alpha)^{2/\alpha}\leq \frac{(\lambda^*)^{1/\alpha}}{|z^*-z|}\leq \frac{1}{|z^*-z|}$, and we have
\begin{align*}
\int_{z^*}^{z^*+V|z-z^*|} g(x) dx &\geq \lambda^* \alpha^2 \frac{ (\lambda^*)^{(1-\alpha)/\alpha} (1-\alpha)^{2(1-\alpha)/\alpha} |z^*-z|
^{1-\alpha}}{1-\alpha}\\
&= \frac{ (\lambda^*)^{1/\alpha} \alpha^2 (1-\alpha)^{2(1-\alpha)/\alpha} |z^*-z|
^{1-\alpha} }{1-\alpha}.
\end{align*}
\item If $z \geq z^*$, we obtain the same way
\begin{equation*}
\int_{z}^{z+1} g(x) dx \geq  \frac{ \lambda^{1/\alpha} \alpha^2 (1-\alpha)^{2 (1-\alpha)/\alpha} |z^*-z|
^{1-\alpha}}{ 1-\alpha}.
\end{equation*}
\end{itemize}
Finally, for any $z^*$ in $\mbb{R}$, using the following inequalities
\begin{equation}
\forall x,y \in[0,1], 1-(1-|x|)_+^{1-\alpha}\geq (1-\alpha) (1\wedge |x|) \text{  and  } \left( \sqrt{x} - \sqrt{y} \right)^2 \geq \left(x-y\right)^2/4,
\end{equation}
we get
\begin{align*}
\int g(x) dx \geq \mathbbm{1}_{|z|\geq|z^*|} &\left[  \frac{ (\lambda)^{1/\alpha} \alpha^2 (1-\alpha)^{2(1-\alpha)/\alpha} |z^*-z|^{1-\alpha} }{ 1-\alpha} + \left( \lambda^* - \lambda \right)^2 (1\wedge|z^*|) \right]\\
\mathbbm{1}_{|z|<|z^*|} &\left[  \frac{ (\lambda^*)^{1/\alpha} \alpha^2 (1-\alpha)^{2(1-\alpha)/\alpha} |z^*-z|^{1-\alpha} }{ 1-\alpha} + \left( \lambda^* - \lambda \right)^2 (1\wedge|z|) \right].
\end{align*}
\begin{itemize}
\item If $|z|\geq |z^*|$:
\begin{itemize}
\item if $\lambda > c\lambda^*$, then
\begin{align*}
\int g(x) dx &\geq \frac{ (\lambda^*)^{1/\alpha} c^{1/\alpha} \alpha^2 (1-\alpha)^{2(1-\alpha)/\alpha} |z^*-z|^{1-\alpha} }{ 1-\alpha} + \left( \lambda^* - \lambda \right)^2 (1\wedge|z^*|)\\
&\geq C_1(\alpha,c) \left[ (\lambda^*)^{1/\alpha} |z^*-z|^{1-\alpha} + (1\wedge|z^*|) \left( \lambda^* - \lambda \right)^2  \right]
\end{align*}
with $C_1(\alpha,c)= 1 \bigwedge \frac{c^{1/\alpha} \alpha^2 (1-\alpha)^{2(1-\alpha)/\alpha}}{1-\alpha}$;
\item otherwise $\int g(x) dx \geq (\lambda^*)^2 (1-c)^2 (1\wedge|z^*|)$,
\begin{equation*}
(\lambda^*)^{1/\alpha} |z^*-z|^{1-\alpha} + (1\wedge|z^*|) \left( \lambda^* - \lambda \right)^2 \leq (\lambda^*)^{1/\alpha} \frac{1}{(1-\alpha)^{2(1-\alpha)/\alpha}} + (1\wedge|z^*|)
\end{equation*}
and finally
\begin{align*}
\int g(x) dx \geq &\frac{(\lambda^*)^2 (1-c)^2 (1\wedge|z^*|)}{(\lambda^*)^{1/\alpha} \frac{1}{(1-\alpha)^{2(1-\alpha)/\alpha}} + (1\wedge|z^*|)}\\
&\times \left[ (\lambda^*)^{1/\alpha} |z^*-z|^{1-\alpha} + (1\wedge|z^*|) \left( \lambda^* - \lambda \right)^2 \right].
\end{align*}
\end{itemize}
\item If $|z|<|z^*|$:
\begin{itemize}
\item if $|z| \geq d |z^*|$, then
\begin{align*}
\int g(x) dx &\geq \frac{ (\lambda^*)^{1/\alpha} \alpha^2 (1-\alpha)^{2(1-\alpha)/\alpha} |z^*-z|^{1-\alpha} }{ 1-\alpha} + \left( \lambda^* - \lambda \right)^2 d(1 \wedge |z^*|)\\
&\geq C_2(\alpha,d) \left[ \left( \lambda^*\right)^{1/\alpha} |z-z^*|^{1-\alpha} + \left( \lambda^*-\lambda \right)^2 \left( 1 \wedge |z^*| \right) \right],
\end{align*}
with $C_2(\alpha,d)= d \wedge \frac{ \alpha^2 (1-\alpha)^{2(1-\alpha)/\alpha} }{ 1-\alpha}$;
\item otherwise $\int g(x) dx \geq \frac{ (\lambda^*)^{1/\alpha} \alpha^2 (1-\alpha)^{2(1-\alpha)/\alpha} |z^*|^{1-\alpha} (1-d)^{1-\alpha} }{ 1-\alpha}$ and
\begin{equation*}
\left( \lambda^*\right)^{1/\alpha} |z-z^*|^{1-\alpha} + \left( \lambda^*-\lambda \right)^2 \left( 1 \wedge |z^*| \right) \leq \left( \lambda^*\right)^{1/\alpha} \frac{1}{(1-\alpha)^{2(1-\alpha)/\alpha}} + \left( 1 \wedge |z^*| \right)
\end{equation*}
and finally
\begin{align*}
\int g(x) dx \geq &\frac{ (\lambda^*)^{1/\alpha} \alpha^2 (1-\alpha)^{2(1-\alpha)/\alpha} |z^*|^{1-\alpha} (1-d)^{1-\alpha} /(1-\alpha)}{(\lambda^*)^{1/\alpha} \frac{1}{(1-\alpha)^{2(1-\alpha)/\alpha}} + (1\wedge|z^*|)}\\
&\times \left[ (\lambda^*)^{1/\alpha} |z^*-z|^{1-\alpha} + (1\wedge|z^*|) \left( \lambda^* - \lambda \right)^2 \right].
\end{align*}
\end{itemize}
\end{itemize}
Finally,
\begin{equation*}
\int g(x) dx \geq C(\alpha,z^*,\lambda^*) \left[ \left( \lambda^* \right)^{1/\alpha} \left( 1 \wedge |z-z^*|^{1-\alpha} \right) + \left( \lambda^* - \lambda \right)^2 \left( 1 \wedge |z^*| \right) \right],
\end{equation*}
with
\begin{align*}
C(\alpha,z^*,\lambda^*) = \min \bigg( &1, \frac{c^{1/\alpha} \alpha^2 (1-\alpha)^{2(1-\alpha)/\alpha}}{1-\alpha}\\
&\frac{(\lambda^*)^2 (1-c)^2 (1\wedge|z^*|)}{(\lambda^*)^{1/\alpha} \frac{1}{(1-\alpha)^{2(1-\alpha)/\alpha}} + (1\wedge|z^*|)},\\
&d, \frac{ \alpha^2 (1-\alpha)^{2(1-\alpha)/\alpha} }{ 1-\alpha},\\
&\frac{ (\lambda^*)^{1/\alpha} \alpha^2 (1-\alpha)^{2(1-\alpha)/\alpha} |z^*|^{1-\alpha} (1-d)^{1-\alpha} /(1-\alpha)}{(\lambda^*)^{1/\alpha} \frac{1}{(1-\alpha)^{2(1-\alpha)/\alpha}} + (1\wedge|z^*|)} \bigg)\\
= \min \bigg( &1, \frac{c^{1/\alpha} \alpha^2 (1-\alpha)^{2(1-\alpha)/\alpha}}{1-\alpha}, d,\\
&\frac{(\lambda^*)^2 (1-c)^2 (1\wedge|z^*|)}{(\lambda^*)^{1/\alpha} \frac{1}{(1-\alpha)^{2(1-\alpha)/\alpha}} + (1\wedge|z^*|)},\\
&\frac{ (\lambda^*)^{1/\alpha} \alpha^2 (1-\alpha)^{2(1-\alpha)/\alpha} |z^*|^{1-\alpha} (1-d)^{1-\alpha} /(1-\alpha)}{(\lambda^*)^{1/\alpha} \frac{1}{(1-\alpha)^{2(1-\alpha)/\alpha}} + (1\wedge|z^*|)} \bigg).
\end{align*}
\subsection*{Proof of Lemma \ref{lem:spike_aux_2}}
Without loss of generality, we take $z^* > 0$.\\
$\bullet$ For $x\in]z^*,z^*(1+a)[$, $a<(z^*)^{-1}$ we have
\begin{align*}
g(x) &= \frac{1}{|x-z^*|^{\alpha}} \left( \sqrt{ (1-\lambda^*) \frac{|x-z^*|^{\alpha}}{|x|^{\alpha}}\mathbbm{1}_{|x|\in(0,1]} + \lambda^* }\right.\\
&- \left. \sqrt{ (1-\lambda) \frac{|x-z^*|^{\alpha} }{ |x|^{\alpha} } \mathbbm{1}_{|x|\in(0,1]} + \lambda \frac{|x-z^*|^{\alpha}}{|x-z|^{\alpha}} \mathbbm{1}_{|x-z|\in(0,1]} } \right)^2.
\end{align*}
and
\begin{equation*}
\frac{|x-z^*|}{|x-z|} \leq \frac{|x-z^*|}{|x|} \leq \frac{a}{1+a} \leq a.
\end{equation*}
We get
\begin{align*}
\int_{z^*}^{z^*+a} g(x) dx &\geq \left( \sqrt{ \lambda^*} - \sqrt{ a^{\alpha} } \right)^2 \int_{z^*}^{z^*+a} \frac{dx}{|x-z^*|^{\alpha}}\\
&= \left( \sqrt{ \lambda^*} - \sqrt{ a^{\alpha} } \right)^2 \frac{\left( a z^* \right)^{1-\alpha}}{1-\alpha}.
\end{align*}
We take $a = (\lambda^*)^{1/\alpha} (1-\alpha)^{2/\alpha}\leq \frac{1}{z^*}$, and we have
\begin{equation*}
\int_{z^*}^{z^*+a} g(x) dx \geq \lambda^* \alpha^2 \frac{ (\lambda^*)^{(1-\alpha)/\alpha} (1-\alpha)^{2(1-\alpha)/\alpha} (z^*)^{1-\alpha}}{1-\alpha}.
\end{equation*}
Otherwise $a=1/z^*\leq (\lambda^*)^{1/\alpha} (1-\alpha)^{2/\alpha}$ and
\begin{equation*}
\int_{z^*}^{z^*+a} g(x) dx \geq \lambda^* \alpha^2 \frac{1}{1-\alpha}.
\end{equation*}
Finally,
\begin{equation*}
\int_{z^*}^{z^*+1} g(x) dx \geq \lambda^* \alpha^2 \frac{1 \bigwedge \left[(\lambda^*)^{(1-\alpha)/\alpha} (1-\alpha)^{2(1-\alpha)/\alpha} (z^*)^{1-\alpha} \right] }{1-\alpha}.
\end{equation*}
\subsection*{Proof of Lemma \ref{lem:spike_aux_3}}
Without loss of generality, we take $z^*\geq 0$.
\begin{itemize}
\item If $z\geq z^*+\frac{1}{(1-\alpha)^{2/\alpha}}$. For $x\in]z^*\vee 1,(z^*+1)\wedge (z-1)[$, we have
\begin{equation*}
g(x) = \frac{\lambda^*}{|x-z^*|^{\alpha}}.
\end{equation*}
One can prove that
\begin{equation*}
|z-z^*|-1 \geq \frac{1}{(1-\alpha)^{2/\alpha}}-1 \geq 1.
\end{equation*}
\begin{itemize}
\item If $z^*\geq 1$, then
We get
\begin{align*}
\int_{z^*}^{z^*+1} g(x) dx &\geq \frac{\lambda^*}{1-\alpha} \left[ 1 \bigwedge |z-z^*|-1 \right]^{1-\alpha}\\
&\geq \frac{\lambda^*}{1-\alpha}.
\end{align*}
\item If $z^*\leq 1$, then
\begin{align*}
\int_1^{ (z^*+1) \wedge (z-1) } g(x) dx &\geq \frac{\lambda^*}{1-\alpha} \left[ 1 \bigwedge (|z-z^*|-1)^{1-\alpha} - \left( 1-z^* \right)^{1-\alpha} \right]\\
&\geq \frac{\lambda^*}{1-\alpha} \left[ 1 - \left(1-z^*\right)^{1-\alpha} \right].
\end{align*}
\end{itemize}
\item If $z^*\geq z+\frac{1}{(1-\alpha)^{2/\alpha}}$, we get
\begin{equation*}
\int_{z^*}^{z^*+1} g(x) dx = \frac{\lambda^*}{1-\alpha}.
\end{equation*}
\end{itemize}
Finally,
\begin{equation*}
\int_{z^*}^{z^*+1} g(x) dx = \frac{\lambda^*}{1-\alpha} \left[ 1 - (1-z^*)_+^{1-\alpha} \right] \geq \lambda^* (1\wedge z^*).
\end{equation*}
\section{Lemmas}
\subsection*{Proof of Lemma \ref{lem:vc}}
The different arguments used in this proof are from Proposition 42 of Baraud \emph{et al} \cite{baraudinventiones} and Lemmas 2.6.15 and 2.6.16 from van der Vaart \& Wellner \cite{VanDerVaart}.
\begin{itemize}
\item We have 
\begin{equation*}
\mcal{C} = \square^{-1} \circ \left\{ x\mapsto \pi\sigma \left[ 1 + \left(\frac{x-z}{\sigma}\right)^2 \right]; \sigma>0, z\in\mbb{R} \right\},
\end{equation*}
where $\square^{-1}$ is the inverse function on $(0,+\infty)$. Since
\begin{equation*}
\left\{ x\mapsto \pi\sigma \left[ 1 + \left(\frac{x-z}{\sigma}\right)^2 \right]; \sigma>0, z\in\mbb{R} \right\} \subset \mbb{R}_2[x] = \left\{ x\mapsto ax^2 + bx + c; (a,b,c)\in\mbb{R}^3 \right\}
\end{equation*}
and $\square^{-1}$ is monotone, we get that $V(\mcal{C})\leq 3 +2$.
\item We have 
\begin{equation*}
\mcal{G} = \exp \circ \left\{ x\mapsto -\frac{1}{2} \log(2\pi \sigma^2) - \frac{(x-z)^2}{2\sigma^2}; \sigma>0, z\in\mbb{R} \right\}.
\end{equation*}
Since $\left\{ x\mapsto -\frac{1}{2} \log(2\pi \sigma^2) - \frac{(x-z)^2}{2\sigma^2}; \sigma>0, z\in\mbb{R} \right\} \subset \mbb{R}_2[x]$ and $\exp$ is monotone, we get that $V(\mcal{C})\leq 3 +2$.
\item We have 
\begin{equation*}
\mcal{L} = \frac{1}{2} \exp \circ -\sqrt{~} \circ \left\{ x \mapsto (x-z)^2; z\in\mbb{R} \right\}.
\end{equation*}
Since $\left\{ x\mapsto (x-z)^2; z\in\mbb{R} \right\} \subset \left\{ x\mapsto ax^2 + bx + c; (a,b,c)\in\mbb{R}^3 \right\}$, $\exp$ and $\sqrt{~}$ are monotone, we get that $V(\mcal{C})\leq 3 +2$. 
\item Azzalini \& Capitanio \cite{azzalini2013skew} proved that the density function $x\mapsto $ is unimodal, therefore the translation family $\mcal{SG}_{\alpha}$ is VC-subgraph with VC-index at most 10 (see Section \ref{sec:performance_estimator}).
\end{itemize}
\subsection*{Proof of Lemma \ref{lem:upper_bound_parameters}}
With Young's inequality, we can easily prove the following inequality
\begin{equation*}
\forall x,y,z\geq0, \left( \sqrt{\suml_{k\in[K]} x_k z_k } - \sqrt{\suml_{k\in[K]} x_k y_k} \right)^2 \leq \suml_{k\in[K]} x_k (\sqrt{z_k}-\sqrt{y_k})^2.
\end{equation*}
Therefore, we can apply to get an upper bound on the Hellinger distance between mixture distributions. For $w,v\in\mcal{W}_K$ and $F_k,G_k\in\mscr{P}$ for all $k\in\{1,\dots,K\}$, we have
\begin{align*}
h\left( \suml_{k \in [K]} w_k F_k, \suml_{k \in [K]} v_k G_k \right) &\leq h\left( \suml_{k \in [K]} w_k F_k, \suml_{k \in [K]} w_k G_k \right)\\
&+ h\left( \suml_{k \in [K]} w_k G_k, \suml_{k \in [K]} v_k G_k \right)\\
\leq &\sqrt{ \suml_{k \in [K]}{ w_k h^2 \left( F_k, G_k \right) }} + h\left( w, v\right)\\
&\leq \max_{k\in[K]} h(F_k,G_k) + h(w,v).
\end{align*}
\subsection*{Proof of Lemma \ref{lem:dist_m_k_delta}}
We prove by induction that
\begin{equation}
\label{eq:prop_dist_m_k_delta_aux1}
\forall \delta\in(0,1/K], \sup\limits_{w\in\mcal{M}_K} h^2(w,\mcal{W}_{K,\delta}) \leq 1-\sqrt{1-(K-1)\delta}.
\end{equation}
$\bullet$ Assume (\ref{eq:prop_dist_m_k_delta_aux1}) holds true for $K\geq 2$. Let $\delta$ be in $(0,1/(K+1))$ and $w$ be in $\mcal{W}_{K+1}$. Without loss of generality we consider $w_1\leq w_2 \leq \dots \leq w_k \leq w_{K+1}$. We define the function $r$ by 
\begin{equation*}
r : \left| \begin{array}{ll}
\mcal{W}_{K+1} &\rightarrow \mcal{W}_K\\
w &\mapsto \left\{ \begin{array}{ll} \left( \frac{w_2}{1-w_1}, \frac{w_3}{1-w_1}, \dots, \frac{w_k}{1-w_1} \right) & \text{  for }w_1\neq 0,\\
\left(\frac{1}{K},\frac{1}{K},\dots,\frac{1}{K}\right)&\text{  for  }w_1=1,
\end{array}\right.
\end{array}\right.
\end{equation*}
and informally $r^{-1}$ by
\begin{equation*}
r^{-1} : \left| \begin{array}{ll}
\mcal{W}_K \times [0,1) &\rightarrow \mcal{W}_{K+1}\\
(w',a) &\mapsto \left(a,(1-a)w'_1,\dots, (1-a) w'_K \right).
\end{array}\right.
\end{equation*}
\begin{itemize}
\item If $w_1\geq\delta$ then $w\in\mcal{M}_{K+1,\delta}$ and $h(\pi,\mcal{W}_{K+1,\delta})=0$.
\item Otherwise $w_1<\delta$ and we build a distribution $v\in\mcal{W}_{K+1,\delta}$ to approximate $w$. Take $\eta=\delta/(1-\delta)\in(0,1/K]$. From (\ref{eq:prop_dist_m_k_delta_aux1}), there exists $v'\in\mcal{M}_{K,\eta}$ such that $h(r(w),v')\leq 1 - \sqrt{1-K\eta} $. Now take $v=r^{-1}(\delta,v')$. We have $v_1=\delta$ and for $j\geq 2$, $v_j=(1-\delta)v'_{j-1}\geq (1-\delta)\eta=\delta$. Therefore $v$ belongs to $\mcal{W}_{K+1,\delta}$. We have 
\begin{align*}
h^2(w,v) &= \frac{1}{2} \left[ \left( \sqrt{w_1} - \sqrt{\delta} \right)^2 + \left( \sqrt{1-w_1} - \sqrt{1-\delta} \right)^2 \right]\\
&+ \sqrt{1-w_1}\sqrt{1-\delta} h^2(r(w),v')\\
&\leq \left[ 1 - \sqrt{1-\delta} \right] + \sqrt{1-\delta} \left[ 1 - \sqrt{1-(K-1)\eta}\right]\\
&= 1-\sqrt{1-(K-1)\delta/(1-\delta)}\\
&\leq 1- \sqrt{1-K\delta}.
\end{align*}
\end{itemize}
$\bullet$ We now prove (\ref{eq:prop_dist_m_k_delta_aux1}) for $K=2$. Let $w$ be in $\mcal{W}_2$ and without loss of generality assume that $w_1\leq w_2$. Once again we only need to consider $w_1<\delta$. Then we take $v=(\delta,1-\delta)$ and we get
\begin{align*}
h^2(w,\mcal{W}_{2,\delta}) &\leq h^2(w,v)\\
&= \frac{1}{2} \left[ \left(\sqrt{w_1}-\sqrt{\delta}\right)^2 + (\sqrt{1-w_1}-\sqrt{1-\delta})^2\right]\\
&\leq  1 - \sqrt{1-\delta}.
\end{align*}
This ends the proof of Lemma \ref{lem:dist_m_k_delta}.
\subsection*{Proof of Lemma \ref{lem:aux_rho_dim}}
The lemma is an intermediate result in the proof of Theorem 2 of Baraud \& Chen \cite{baraudglm}. We follow the proof with $h(x)=a+b\log(1/x)$. We define $\overline{Z}(f) = \sup_{f\in\mcal{F}} \left| \suml_{i=1}^n \epsilon_i f(X_i) \right|$ where $\epsilon_1,\dots,\epsilon_n$ are i.i.d. Rademacher random variables. Everything is the same until equation (42) in the proof. We get
\begin{align*}
\mbb{E} \left[ \overline{Z}(\mcal{F}) \right] &\leq \sqrt{2n}\frac{1+q}{1-q} \int_0^B \sqrt{\log 2 + 2 a + b\log(1/q) + 2b\log(1/u)} du,
\end{align*}
with $B=\sqrt{\sigma^2+\frac{8\mbb{E}\left[\overline{Z}(\mcal{F})\right]}{n}}\wedge 1$. With Lemma 2 \cite{baraudglm}, we get
\begin{align*}
\mbb{E} \left[ \overline{Z}(\mcal{F}) \right] &\leq 16A^2 + A \sqrt{2n\sigma^2},
\end{align*}
with $A=\frac{1+q}{1-q}\left( 1 + \frac{b}{\log 2 + 2 a + b\log(1/q)} \right) \sqrt{\log 2 + 2 a + b\log(1/q) + 2b\log(1/\sigma)}$. Classical symmetrisation arguments imply $\mbb{E} \left[ Z(\mcal{F}) \right] \leq 2 \mbb{E} \left[ \overline{Z}(\mcal{F}) \right]$.
\subsection*{Proof of Lemma \ref{lem:unif_entropy_vc_mix}}
We write $\phi=\psi\left(\sqrt{\cdot/\overline{p}}\right)$.
\begin{lemme}
\label{lem:mix_entropy_phi_delta}
For any probability distribution $R$ on $(\mscr{X},\mcal{X})$, for $w,v\in\mcal{W}_K$ such that $w_k,v_k\geq \delta$ for $k=1,\dots,K$ and for any densities $q_1,\dots,q_K,r_1,\dots,r_K$, we have 
\begin{align*}
&||\phi\circ(w_1 q_1 + \dots + w_k q_K) - \phi\circ(v_1 r_1 + \dots + v_K r_K)||_{2,R}\\
&\leq \frac{1}{\sqrt{\delta}} \suml_{k=1}^K ||\phi\circ q_k - \phi\circ  r_k||_{2,R} + \frac{2}{\delta} ||w-v||_{\infty},
\end{align*}
where $||w-v||_{\infty}=\max\limits_{k\in[K]} |w_k-v_k|$.
\end{lemme}
This lemma implies, for any probability distribution $R$ on $(\pmb{\mscr{X}},\pmb{\mcal{X}})=(\mscr{X}^n,\mcal{X}^{\otimes n})$,
\begin{equation}
\label{eq:mix_entropy_aux_1}
\log N\left(\epsilon,\mcal{F}^{\mcal{Q}_{\delta}},||\circ||_{2,R}\right) \leq \log N\left(\epsilon_{K+1}, \mcal{M}_K,||\cdot||_{\infty}\right) + \suml_{k=1}^K \log N\left(\epsilon_k,\mcal{G}_k,||\cdot||_{2,R} \right),
\end{equation}
where $\mcal{G}_k := \left\{ \phi \circ f \bigg| F \in \mscr{F}_k \right\}$ for $k=1,\dots,K$ and $\epsilon=\frac{\epsilon_1+\dots+\epsilon_K}{\sqrt{\delta}}+\frac{2\epsilon_{K+1}}{\delta}$.
The following lemmas provide bounds on the covering numbers involved in (\ref{eq:mix_entropy_aux_1}). First, we use Proposition 42 \cite{baraudinventiones} (VC + monotone) and Lemma 1 from Baraud \& Chen \cite{baraudglm} to bound the covering numbers on the $\mcal{G}_k,k=1,\dots,K$.
\begin{lemme}
For any probability measure $R$ on $(\mscr{X},\mcal{X})$, for all $\epsilon_k\in(0,2)$,
\begin{equation}
\label{eq:mix_entropy_aux_2}
\log N(\epsilon_k,\mcal{G}_k,R) \leq \log\left( e V_k (8e)^{V_k-1 }\right) + 2(V_k-1)\log(1/\epsilon_k).
\end{equation}
\end{lemme}
We also have a bound on the covering number of $\mcal{M}_K$ given by the following lemma.
\begin{lemme}
\label{lem:mk_entropy_d}
For $\epsilon_{K+1} >0$, we have
\begin{equation}
\label{eq:mix_entropy_aux_3}
\log N\left(\epsilon_{K+1},\mcal{W}_K,||\cdot||_{\infty}\right)\leq K \log\left( \frac{3}{\epsilon_{K+1}} \right).
\end{equation}
\end{lemme}
We combine (\ref{eq:mix_entropy_aux_1}), (\ref{eq:mix_entropy_aux_2}) and (\ref{eq:mix_entropy_aux_3}). For $\epsilon\in(0,2)$ and $\delta\in(0,1/K]$, we take
\begin{equation*}
\epsilon_{K+1} = \epsilon \frac{\delta}{2} \frac{K}{K+\suml_{k=1}^K 2(V_k-1)} \text{  and  } \epsilon_j = \epsilon \sqrt{\delta} \frac{2(V_j-1)}{K+\suml_{k=1}^K 2(V_k-1)},j=1,\dots,K.
\end{equation*}
we get
\begin{align*}
\log N\left(\epsilon,\mcal{F}^{\mcal{Q}_{\delta}},||\cdot||_{2,R}\right) &\leq K \log\left(\frac{3}{\epsilon_{K+1}} \right) + \log\left( e^K \left( \prod_k V_k \right) (8e)^{\suml_k (V_k-1)} \right)\\
&+ \suml_{k=1}^K 2(V_k-1) \log(1/\epsilon_k)\\
&= K\log\left( \frac{6}{\epsilon\delta} \frac{K+\suml_{k=1}^K 2(V_k-1)}{K} \right)\\
&+ \log\left( e^K \left( \prod_k V_k \right) (8e)^{\overline{V}-K} \right)\\
&+ \suml_{k=1}^K 2(V_k-1) \log\left( \frac{1}{\epsilon \sqrt{\delta}} \frac{K+\suml_{j=1}^K 2(V_j-1)}{2(V_k-1)} \right)\\
&= \log\left( \frac{ \left[ K+\suml_{j=1}^K 2(V_j-1) \right]^{K+\suml_{j=1}^K 2(V_j-1)}  }{K^K \times \prod_{k=1}^K [2(V_k-1)]^{2(V_k-1)}} \right)\\
&+  \overline{V} \log\left( \left[\prod_k V_k\right]^{1/\overline{V}} \right)\\
&+ \log\left( \frac{ 6^K e^{\overline{V}} 8^{\overline{V}-K}}{\delta^{\overline{V}}} \right) + (2\overline{V}-K) \log( 1/\epsilon).
\end{align*}
We use the following lemma to obtain a simpler expression.
\begin{lemme}
\label{lem:aux_bound_log}
For all $x_1,\dots,x_n\geq 0$, we have
\begin{equation*}
\log\left( \frac{(x_1+\dots+x_n)^{x_1+\dots+x_n}}{x_1^{x_1}\dots x_n^{x_n}} \right) \leq (x_1+\dots+x_n)\log(n),
\end{equation*}
and
\begin{equation*}
\left( \prod\limits_{i=1}^n x_i \right)^{\frac{1}{x_1+\dots+x_n}} \leq \left( e^{\frac{1}{e}}\right)^{^\frac{1}{n^n}} \leq e^{1/e}.
\end{equation*}
\end{lemme}
Then, we get
\begin{align*}
\log N\left(\epsilon,\mcal{F}^{\mcal{Q}_{\delta}},||\cdot||_{2,R}\right) &\leq \left[ K+\suml_{j=1}^K 2(V_j-1) \right]\log(K+1) +  \overline{V} \log\left(e^{1/e}\right)\\
&+ \log\left( \frac{ e^{\overline{V}} 8^{\overline{V}}}{\delta^{\overline{V}}} \right) + (2\overline{V}-K) \log( 1/\epsilon)\\
&\leq \overline{V} \log\left( \frac{ e^{1+1/e}  8 (K+1)^2 }{ \delta } \right) + 2\overline{V} \log( 1/\epsilon).
\end{align*}
\subsection*{Proof of Lemma \ref{lem:ghosal_dtv}}
The bound is  obtained following the proofs of lemmas in Ghosal \& van der Vaart \cite{ghosal2001}
\begin{itemize}
\item 1st step: \\
For $|x|>a$ we have,
\begin{align}
p_H( x) &= \int \frac{1}{\sqrt{2\pi \sigma^2}} \exp\left( - \frac{(x-z)^2}{2 \sigma^2} \right) dH(z,\sigma)\nonumber\\
&\leq \frac{1}{\sqrt{2\pi\underline{\sigma}^2}} \exp\left( - \frac{(|x|-a)^2}{2 \overline{\sigma}^2} \right). \label{eq:p_h_upper_bound}
\end{align}
\item 2nd step:\\
See Lemma A.1 in  Ghosal \& van der Vaart \cite{ghosal2001}. Take $N=k(2k-1)+1$. There is a discrete distribution $H'$ with at most $K$ support points in $[-a,a]\times[\underline{\sigma},\overline{\sigma}]$ such that
\begin{equation}
\label{eq:eq_param_vdv}
\int z^l \sigma^{-(2j+1)} dH(z,\sigma) =  \int z^l \sigma^{(2j+1)} dH'(z,\sigma)
\end{equation}
for $l=0,\dots,2k-2$ and $j=0,\dots,k-1$.
Because of (\ref{eq:eq_param_vdv}) we get
\begin{align*}
\int\displaystyle{ \suml_{j=0}^{k-1} \frac{(-1)^j \sigma^{-(2j+1)} (x-z)^{2j} }{j!} dH(z,\sigma)} = \int\displaystyle{\suml_{j=0}^{k-1} \frac{(-1)^j \sigma^{-(2j+1)} (x-z)^{2j} }{j!}  dH'(z,\sigma)},
\end{align*}
for $x\in\mbb{R}$. Taylor's expansion of the exponential function (\cite{ghosal2001}),
\begin{equation*}
\left| \exp\left(-\frac{(x-z)^2}{2 \sigma^2}\right) - \suml_{j=0}^{k-1} \frac{\left(- \frac{(x-z)^2}{2 \sigma^2}\right)^j}{j!} \right| \leq \left( \frac{e (x-z)^2 }{k 2 \sigma^2} \right)^k.
\end{equation*}
Therefore,
\begin{align*}
\sqrt{2\pi} \sup_{|x|\leq M} | p_H( x ) & - p_{H'}( x) |\\
&= \sup_{|x|\leq M} \left| \int \frac{1}{\sigma} \exp\left( -\frac{(x-z)^2}{2 \sigma^2} \right) dH(z,\sigma) \right.\\
&\left.- \int \frac{1}{\sigma} \exp\left( -\frac{(x-z)^2}{2 \sigma^2} \right) dH'(z,\sigma) \right|\\
&= \sup_{|x|\leq M} \left| \int \frac{1}{\sigma} \left[  \exp\left( -\frac{(x-z)^2}{2 \sigma^2} \right) - \suml_{j=0}^{k-1} \frac{\left(- \frac{(x-z)^2}{2 \sigma^2}\right)^j}{j!}  \right] dH(z,\sigma) \right.\\
&\left.- \int \frac{1}{\sigma} \left[ \exp\left( -\frac{(x-z)^2}{2 \sigma^2} \right) - \suml_{j=0}^{k-1} \frac{\left(- \frac{(x-z)^2}{2 \sigma^2}\right)^j}{j!} \right] dH'(z,\sigma) \right|\\
&\leq 2 \sup_{\substack{|x|\leq M\\ |z|\leq a\\
\underline{\sigma}\leq\sigma\leq\overline{\sigma}}} \frac{1}{\sigma} \left| \exp\left( -\frac{(x-z)^2}{2 \sigma^2} \right) - \suml_{j=0}^{k-1} \frac{\left(- \frac{(x-z)^2}{2 \sigma^2}\right)^j}{j!} \right|\\
&\leq 2 \sup_{\substack{|x|\leq M\\ |z|\leq a\\
\underline{\sigma}\leq\sigma\leq\overline{\sigma}}} \frac{1}{\sigma} \left( \frac{e (x-z)^2 }{k 2 \sigma^2} \right)^k\\
&\leq \frac{2}{\underline{\sigma}} \left( \frac{e (M+a)^2 }{k 2 \underline{\sigma}^2} \right)^k.\label{eq:p_h_difference}
\end{align*}
\end{itemize}
Obviously, the inequality (\ref{eq:p_h_upper_bound}) holds also for $p_{H'}$. We combine it with the last one we obtained in order to bound the total variation distance. Therefore, for $M=ma$, $m > 1$, we have
\begin{align*}
d_{TV}\left( P_H, P_{H'} \right) &= \frac{1}{2}  \int |p_H(x) - p_{H'}(x) | dx\\
&\leq M \sup_{|x|\leq M} \left| p_H(x) - p_{H'}(x)\right| + \frac{1}{2}\int_{|x|>M} p_H(x)\vee p_{H'}(x) dx\\
&\leq \frac{\sqrt{2/\pi}}{\underline{\sigma}} M \left( \frac{ e (M+a)^2}{2 k \underline{\sigma}^2} \right)^k + \frac{1}{2} \int_{|x|>M} \frac{1}{\sqrt{2\pi\underline{\sigma}^2}} \exp\left( - \frac{(|x|-a)^2}{2 \overline{\sigma}^2} \right) dx\\
&\leq \frac{\sqrt{2/\pi}}{\underline{\sigma}} M   \left( \frac{ e (M+a)^2}{2 k \underline{\sigma}^2} \right)^k + \frac{\overline{\sigma}}{\underline{\sigma}} \int_{x>M} \frac{1}{\sqrt{2\pi\overline{\sigma}^2}} \exp\left( - \frac{(x-a)^2}{2 \overline{\sigma}^2} \right) dx\\
&\leq \frac{\sqrt{2/\pi}}{\underline{\sigma}} am  \left( \frac{ e a^2 (1+m)^2}{2 k \underline{\sigma}^2} \right)^k + \frac{\overline{\sigma}}{2 \underline{\sigma}} \exp\left( - \frac{(m-1)^2 a^2}{2\overline{\sigma}^2}\right). 
\end{align*}
Finally, writing $A=a/\underline{\sigma}$ and $R=\overline{\sigma}/\underline{\sigma}$, we have
\begin{equation*}
d_{TV}\left( P_H, P_{H'} \right) \leq \inf\limits_{m>1} \left\{ \frac{\sqrt{2/\pi}}{\underline{\sigma}} am  \left( \frac{ e a^2 (1+m)^2}{2 k \underline{\sigma}^2} \right)^k + \frac{\overline{\sigma}}{2 \underline{\sigma}} \exp\left( - \frac{(m-1)^2 a^2}{2\overline{\sigma}^2}\right) \right\}. 
\end{equation*}
\subsection*{Proof of Lemma \ref{lem:aux_identifiability}}
Let $j\in \{0,\dots,K\}$ and assume there is a sequence
\begin{equation*}
(P_n)_n = \left( \suml_{k=1}^{j} w_{k,n} \mcal{N}(z_{k,n},\sigma^2_{k,n}) + \suml_{k=j+1}^K w_{k,n} \cauchy(z_{k,n},\sigma_{k,n}) \right)\in\mscr{Q}_j^{\mbb{N}}
\end{equation*}
such that $\lim\limits_{n\rightarrow\infty} h(P_n,P^*)=0$. The mixing weights are bounded so we can assume we are already considering a sequence such that $w_{k,n} \xrightarrow[n\rightarrow\infty]{} w_{k,\infty}$ for all $k\in\{1,\dots,K\}$. For the other parameters, it is always possible to extract a subsequence $P_{\psi(n)}$ such that for all $k$
\begin{equation*}
z_{k,\psi(n)} \xrightarrow[n\rightarrow\infty]{} \begin{cases}
z_{k,\infty}\in\mbb{R},\\
\text{or  } \pm\infty,
\end{cases}
\text{  and   }
\sigma_{k,\psi(n)} \xrightarrow[n\rightarrow\infty]{} \begin{cases}
\sigma_{k,\infty}\in\mbb{R}^+,\\
\text{or  } +\infty.
\end{cases}
\end{equation*}
We now consider the different cases possible (dropping the dependency on $\psi$ in the notation).
\begin{itemize}
\item If $z_{k,n} \xrightarrow[n\rightarrow \infty]{} \pm\infty$ (without loss of generality we consider $+\infty$ in the proof), for $b\in\mbb{R}$, we have
\begin{align*}
P_n([b,+\infty[) &\geq w_{k,n} \big[ \mathbbm{1}_{k \leq j} \mcal{N}(z_{k,n},\sigma^2_{k,n})([b,+\infty[)\\
&+ \mathbbm{1}_{k > j} \cauchy(z_{k,n},\sigma_{k,n})([b,+\infty[) \big]\\
&\geq \frac{w_{k,n}}{2} \text{  for } n \text{  large enough}.
\end{align*}
Assume $w_{k,\infty}>0$. Since $P^*([b,+\infty[) \xrightarrow[b\rightarrow\infty]{} 0$, there exists $b$ such that $P^*([b,+\infty[) \leq w_{j,\infty}/4$. On the other hand we have $P^*([b,+\infty[) = \lim\limits_{n\rightarrow\infty} P_{\theta_n}([b,+\infty[) \geq w_{k,\infty}/2$. Therefore, it means that $w_{k,\infty} = 0$ and it also holds for $z_{k,n}\rightarrow-\infty$. 
\item If $z_{k,n} \xrightarrow[n\rightarrow\infty]{} z_{k,\infty}\in\mbb{R}$ and $\sigma_{k,n}\xrightarrow[n\rightarrow\infty]{} 0$, for $b>0$ we have
\begin{align*}
P_n([z_{k,\infty}-b,z_{k,\infty}+b]) \geq w_{k,n} &\left( \mathbbm{1}_{k \leq j} \mcal{N}(z_{k,n},\sigma^2_{k,n})([b,+\infty[)\right.\\
& \left. + \mathbbm{1}_{k > j} \cauchy(z_{k,n},\sigma_{k,n})([b,+\infty[) \right) \rightarrow w_{k,\infty}.
\end{align*}
Assume $w_{k,\infty}>0$. Since $P^*([z_{k,\infty}-b,z_{k,\infty}+b]) \xrightarrow[b\rightarrow 0]{} 0$, there exists $b>0$ such that $P^*([z_{k,\infty}-b,z_{k,\infty}+b])\leq w_{j,\infty}/2$. On the other hand we have $P^*([z_{k,\infty}-b,z_{k,\infty}+b]) = \lim\limits_{n\rightarrow\infty} P_n([z_{k,\infty}-b,z_{k,\infty}+b]) \geq w_{k,\infty}$. Therefore, it means that $w_{k,\infty} = 0$.
\item If $z_{k,n} \rightarrow z_{k,\infty}\in\mbb{R}$ and $\sigma_{k,n}\rightarrow \infty$, for $a>0$ we have
\begin{align*}
P_n([-a,a]) &\leq (1-w_{k,n})\\
&+ w_{k,n} \left( \mathbbm{1}_{k \leq j} \mcal{N}(z_{k,n},\sigma^2_{k,n})([-a,a]) + \mathbbm{1}_{k > j} \cauchy(z_{k,n},\sigma_{k,n})([-a,a]) \right)\\
&\xrightarrow[n\rightarrow\infty]{} (1-w_{k,\infty}).
\end{align*}
Since $P^*([-a,+a]) \xrightarrow[a\rightarrow +\infty]{} 1$, we get $w_{k,\infty} = 0$
\end{itemize}
This proves that $P_n$ converges to
\begin{equation*}
P_{\infty} = \suml_{\substack{k\leq j(\lambda)\\
w_{k,\infty}>0}} w_{k,\infty} \mcal{N}(z_{k,\infty},\sigma^2_{k,\infty}) + \suml_{\substack{k > j(\lambda)\\w_{k,\infty}>0}} w_{k,\infty} \cauchy(z_{k,\infty},\sigma_{k,\infty}),
\end{equation*}
and necessarily $P^*=P_{\infty}$. Lemma \ref{lem:aux_identifiability} with the assumptions on $P^*$ implies $j=j^*$ and there exist two permutations $\tau_g,\tau_c$ respectively on $\{1,\dots,j^*\}$ and $\{j^*+1,\dots,K\}$ such that $(\overline{\pi}_k,\overline{z}_k,\overline{\sigma}_k) = (w_{\tau_g(k)},z_{\tau_g(k)},\sigma_{\tau_g(k)})$ for $k$ in $\{1,\dots,j^*\}$ and $(\overline{\pi}_k,\overline{z}_k,\overline{\sigma}_k) = (w_{\tau_c(k)},z_{\tau_c(k)},\sigma_{\tau_c(k)})$ for $k$ in $\{j^*+1,\dots,K\}$.
\subsection*{Proof of Lemma \ref{lem:aux_regular_cauchy_gauss}}
\begin{itemize}
\item The map $(z,\sigma^2) \mapsto g(x;z,\sigma^2)=\phi_{\sigma}(x-z)$ is continuous and differentiable on $\mbb{R}\times\mbb{R}^{+*}$ with
\begin{align*}
\partial_z \phi_{\sigma}(x-z) &= \phi_{\sigma}(x-z) \frac{(x-z)}{\sigma^2}\\
\partial_{\sigma^2} \phi_{\sigma}(x-z) &= \phi_{\sigma}(x-z) \left[ \frac{(x-z)^2}{2 \sigma^4} - \frac{1}{2 \sigma^2} \right].
\end{align*}
Similarly $(z,\sigma) \mapsto f(x;z,\sigma)=\frac{1}{\pi\sigma} \frac{1}{c(x;z,\sigma)}$ is continuous and differentiable on $\mbb{R}\times\mbb{R}^{+*}$ with
\begin{align*}
\partial_z f(x;z,\sigma) &= \frac{1}{\pi \sigma^3} \frac{x-z}{c^2(x;z,\sigma)}\\
\partial_{\sigma} f(x;z,\sigma) &= \frac{1}{\pi \sigma^2 c(x;z,\sigma)} \left[ 1 - \frac{2}{c(x;z,\sigma)} \right].
\end{align*}
Moreover, on can check that we have
\begin{align*}
&\int_{\mbb{R}} \left| \partial_z g(x;z,\sigma^2) \right|^2 \frac{dx}{g(x;z,\sigma^2)} = \int_{\mbb{R}} \frac{(x-z)^2}{\sigma^4} \phi_{\sigma}(x-z) dx < \infty\\
&\int_{\mbb{R}} \left| \partial_{\sigma^2} g(x;z,\sigma^2) \right|^2 \frac{dx}{g(x;z,\sigma^2)} = \int_{\mbb{R}}  \left[ \frac{(x-z)^2}{2 \sigma^4} - \frac{1}{2 \sigma^2} \right]^2 \phi_{\sigma}(x-z) dx < \infty\\
&\int_{\mbb{R}} \left| \partial_z f(x;z,\sigma) \right|^2 \frac{dx}{f(x;z,\sigma)} = \int_{\mbb{R}} \frac{(x-z)^2}{\pi \sigma^5 c^3(x;z,\sigma)} dx < \infty\\
&\int_{\mbb{R}} \left| \partial_{\sigma^2} f(x;z,\sigma^2) \right|^2 \frac{dx}{f(x;z,\sigma)} = \int_{\mbb{R}}  \frac{1}{\pi \sigma^3 c(x;z,\sigma)} \left[ 1 - \frac{2}{c(x;z,\sigma)} \right]^2 dx < \infty.
\end{align*}
\item The function $\theta\mapsto \psi(\cdot;\theta)=\frac{\partial}{\partial \theta} p^{1/2}(\cdot;\theta)$, where
\begin{equation*}
p(x;\theta)=\suml_{k=1}^{j^*} w_k \phi_{\sigma_k}(x-z_k) + \suml_{k=j^*+1}^K \frac{1}{\pi\sigma c(x;z,\sigma)}
\end{equation*}
and
\begin{equation*}
\theta=(w_1,\dots,w_{K-1},z_1,\dots,z_K,\sigma^2_1,\dots,\sigma^2_{j^*},\sigma_{j^*+1},\dots,\sigma_K),
\end{equation*}
is continuous in the space $L_2(\mu)$.
\item We apply Theorem 1 of Meijer \& Ypma \cite{meijer}. For $j^*<K$,
\begin{align*}
&\det(I(\theta))=0\\
 &\Rightarrow \exists \lambda\neq 0, \suml_{k=1}^{j^*} \phi_{\sigma_k}(x-z_k) \left( \frac{w_k \lambda_{z_k}(x-z_k)}{\sigma_k^2} + w_k \lambda_{\sigma_k^2} \left[ \frac{(x-z)^2}{2 \sigma_k^4} - \frac{1}{2 \sigma_k^2} \right] + \lambda_{w_k} \right)\\
&+ \suml_{k=j^*+1}^{K-1} \left(  \frac{w_k \lambda_{z_k}(x-z_k)}{\pi\sigma^3 c^2(x;z_k,\sigma_k)} + \frac{w_k \lambda_{\sigma_k}}{\pi \sigma_k^2} \left[ \frac{1}{c(x;z_k,\sigma_k)} - \frac{2}{c^2(x;z_k,\sigma_k)} \right] + \frac{\lambda_{w_k}}{\pi \sigma_k c(x;z_k,\sigma_k)} \right)\\
&+ (1-w_1-\dots-w_{K-1}) \left(  \frac{\lambda_{z_K} (x-z_K)}{\pi \sigma_K^3 c^2(x;z_K,\sigma_K)} + \frac{ \lambda_{\sigma_K}}{\pi \sigma_K^2} \left[ \frac{1}{c(x;z_K,\sigma_K)} - \frac{2}{c^2(x;z_K,\sigma_K)} \right]  \right)\\
&- \frac{1}{\pi\sigma_K c(x;z_K,\sigma_K)} \suml_{k=1}^{K-1} \lambda_{w_k} = 0 \text{  for  }\mu\text{-almost all  }x.
\end{align*}
For $j^*=K$,
\begin{align*}
&\det(I(\theta))=0\\
&\Rightarrow \exists \lambda\neq 0, \suml_{k=1}^{K-1} \phi_{\sigma_k^2}(x-z_k) \left( w_k \lambda_{z_k} \frac{(x-z_k)}{\sigma_k^2} + w_k \lambda_{\sigma_k^2} \left[ \frac{(x-z)^2}{2 \sigma_k^4} - \frac{1}{2 \sigma_k^2} \right] + \lambda_{w_k} \right)\\
&+ \phi_{\sigma_K}(x-z_K) \left\{ (1-w_1-\dots-w_{K-1}) \left( \lambda_{z_K} \frac{(x-z_K)}{\sigma_K^2} + \lambda_{\sigma_K^2} \left[ \frac{(x-z)^2}{2 \sigma_K^4} - \frac{1}{2 \sigma_K^2} \right]  \right)\right.\\
&- \left. \suml_{k=1}^{K-1} \lambda_{w_k} \right\} = 0 \text{  for  }\mu\text{-almost all  }x.
\end{align*}
\begin{lemme}
\label{lem:aux_linear_ind}
Let $(z_1,\sigma_1),\dots,(z_K,\sigma_K)$ be distinct elements of $\mbb{R}\times \mbb{R}^{+*}$. For any integer $n$, the families
\begin{equation*}
A=\left\{ x \mapsto x^j \phi_{\sigma_i}(x-z_i) ; i\in\{1,\dots,K\}, j\in\{0,\dots,n\} \right\}
\end{equation*}
and
\begin{equation*}
B=\left\{ x \mapsto \frac{x^j}{c^l(x;z_i,\sigma_i)}; i\in\{1,\dots,K\}, l\in\{1,2\}, j\in\{0,1\} \right\}
\end{equation*}
are linearly independent. Moreover, the linear spaces $\vect(A)$ and $\vect(B)$ are orthogonal.
\end{lemme}
This proves that $I(\overline{\theta})$ is non singular.
\item We now check $\inf_{ \substack{||\overline{\theta}-\theta||\geq a\\ P_{\theta} \in\mscr{Q}_{j^*}}} h^2(P_{\overline{\theta}},P_{\theta}) > 0,\forall a>0$. It is a direct consequence of Lemma \ref{lem:aux_identifiability}.
\item $\mscr{Q}(\lambda^*)$ is a regular parametric model. We consider the parameter to be $\sigma$ for the Cauchy distribution and $\sigma^2$ for the Gaussian distribution. Obviously, $(z,\sigma) \mapsto g(x;z,\sigma) = \frac{1}{\pi\sigma} \frac{1}{c(x;z,\sigma)}$, with $c(x;z,\sigma)=1+\left(\frac{x-z}{\sigma}\right)^2$ is continuous and differentiable on $\mbb{R}\times\mbb{R}^{+*}$ with
\begin{align*}
\partial_z g(x;z,\sigma) &= \frac{2(x-z)}{\pi\sigma^3 c^2(x;z,\sigma)}\\
\partial_{\sigma} g(x;z,\sigma) &= \frac{1}{\pi\sigma^2 c(x;z,\sigma)} - \frac{2}{\pi \sigma^2 c^2(x;z,\sigma)}.
\end{align*}
Moreover, on can check that we have
\begin{equation*}
\int_{\mbb{R}} \left| \partial_z g(x;z,\sigma) \right|^2 \frac{dx}{g(x;z,\sigma)} = \int_{\mbb{R}} \frac{4 (x-z)^2}{\pi \sigma^3 c^3(x;z,\sigma)} dx < \infty
\end{equation*}
and
\begin{equation*}
\int_{\mbb{R}} \left| \partial_{\sigma} g(x;z,\sigma) \right|^2 \frac{dx}{g(x;z,\sigma)} = \int_{\mbb{R}}  \frac{1}{\pi\sigma^3 c(x;z,\sigma)} \left[ 1 - \frac{2}{c(x;z,\sigma)} \right]^2 dx < \infty.
\end{equation*}
\item With the results of \cite{Ibragimov}, we get that there is a constant $a^*>0$ such that
\begin{equation*}
\forall P_\theta \in\mscr{Q}(\lambda^*), a^* \frac{||\theta-\overline{\theta}||^2}{1+||\theta-\overline{\theta}||^2} \leq h^2(P^*,P_{\theta}).
\end{equation*} 
\end{itemize}
\subsection*{Proof of Lemma \ref{lem:aux_linear_ind}}
\begin{itemize}
\item Let $f$ be any function in $\vect(A)\cap \vect(B)$. Therefore there are constants $(\lambda_{g,i,j})_{\substack{1\leq i\leq K,\\ 0\leq j\leq n}}$ and $(\lambda_{c,i,l,j})_{\substack{1\leq i\leq K,\\0\leq j\leq 1\leq l\leq 2}}$ such that 
\begin{equation*}
f(x) = \suml_{i=1}^K \suml_{j=0}^n \lambda_{g,i,j} x^j \phi_{\sigma_i}(x-z_i)= \suml_{i=1}^K \suml_{l=1}^2 \suml_{j=0}^1 \lambda_{c,i,l,j} \frac{x^j}{c^l(x;z_i,\sigma_i)}.
\end{equation*}
Since $f\in \vect(A)$, we have $f(x) = o_{\pm\infty}(x^{-k}), \forall k\in\mbb{N}$. Therefore $\lambda_{c,i,l,j}=0$ for all $i,j,l$ and $f=0$. This proves $\vect(A)\cap \vect(B) = \{ 0 \}$.
\item One can check that $>$ is a strict total order such that
\begin{equation*}
(z_1,\sigma_1) > (z_2,\sigma_2) \Rightarrow x^j  \phi_{\sigma_2}(x-z_2)/\phi_{\sigma_1}(x-z_1) \xrightarrow[x\rightarrow + \infty]{} 0,
\end{equation*}
for any $j\in\mbb{N}$. Let $\lambda$ be such that $\suml_{i,j} \lambda_{i,j} x^j \phi_{\sigma_i}(x-z_i)=0$ for all $x$. Without loss of generality, we assume $(z_1,\sigma_1)>\dots>(z_K,\sigma_K)$. Therefore,
\begin{align*}
0 &= \suml_{i,j} \lambda_{i,j} x^j \phi_{\sigma_i}(x-z_i)\\
&= \suml_{i,j} \lambda_{i,j} x^j \phi_{\sigma_i}(x-z_i)/\phi_{\sigma_1}(x-z_1) + \suml_{j} \lambda_{1,j}x^j\\
&= \suml_{j} \lambda_{1,j}x^j + o_{+\infty}(1).
\end{align*}
It implies that $\lambda_{1,j}=0$ for all $j$. Then, we have $\suml_{i\geq 2,j} \lambda_{i,j} x^j \phi_{\sigma_i}(x-z_i)=0$. By induction, we get that $\lambda=0$ which proves that the family is indeed linearly independent.
\item The partial fraction decomposition theorem implies that $B$ is linearly independent.
\end{itemize}
\subsection*{Proof of Lemma \ref{lem:mix_entropy_phi_delta}}
The result is just the combination of the two following lemmas and the triangle inequality.
\begin{lemme}
\label{lem:aux_q}
For any nonnegative measurable functions $r,q_1,q_2$ and any $w\in(0,1)$ we have
\begin{align*}
||\phi\circ(w q_1 + (1-w) r) - \phi\circ(w q_2 + (1-w) r)||_{2,Q} &\leq \frac{1}{\sqrt{\pi}} ||\phi\circ q_1 - \phi\circ  q_2||_{2,Q}.
\end{align*}
\end{lemme}
\begin{lemme}
\label{lem:aux_pi_K}
Let $g_1,\dots,g_K$ be $K$ densities. Let $w,v\in\mcal{M}_{K,\delta}$.
\[ ||\phi \circ (w_1 g_1 + \dots + w_K g_K) - \phi\circ (v_1 g_1 + \dots + v_K g_K)||_{2,Q} \leq \frac{2}{\delta} d(w,v),\]
where $d(w,v)=\max\limits_{1\leq i\leq K} |w_i-v_i|$.
\end{lemme}
Indeed,
\begin{align*}
\left|\left|\phi \circ \left( \suml_{j=1}^K w_j f_j \right) - \phi\circ \left( \suml_{j=1}^K v_j g_j \right) \right|\right|_{2,Q} &\leq \left|\left|\phi \circ \left( \suml_{j=1}^K w_j f_j \right) - \phi\circ \left( \suml_{j=1}^K v_j f_j \right) \right|\right|_{2,Q}\\
&+ \suml_{i=1}^K \left|\left|\phi \circ \left( h_{i-1} \right) - \phi\circ \left( h_i  \right) \right|\right|_{2,Q}\\
&\leq \frac{2}{\delta} d(w,v) + \suml_{i=1}^K \frac{1}{\sqrt{\nu_i}} \left|\left|\phi \circ \left( g_i \right) - \phi\circ \left( f_i  \right) \right|\right|_{2,Q}\\
&\leq \frac{2}{\delta} d(w,v) + \frac{1}{\sqrt{\delta}} \suml_{i=1}^K  \left|\left|\phi \circ \left( g_i \right) - \phi\circ \left( f_i  \right) \right|\right|_{2,Q},
\end{align*}
with $h_i = \suml_{j=1}^i v_j g_j + \suml_{j=i+1}^K v_j f_j$.
\subsection*{Proof of Lemma \ref{lem:mk_entropy_d}}
Let $\epsilon\in(0,1)$. Let $N$ be an integer greater than $\frac{1}{\epsilon}$. We define
\[\mcal{W}_{K,N} := \left\{ w\in\mcal{W}_K \bigg| \forall i, \exists d_i\in\mathbb{N}, w_i=\frac{d_i}{N}\right\}.\]
\begin{itemize}
\item One can easily see that there is a bijection between $\mcal{M}_{K,N}$ and
\[ D_{K,N} := \left\{ d_1,\dots,d_K\in\mathbb{N} \bigg| \suml_i d_i = N \right\}.\]
We refer to Barron \& Klusowski (\cite{barron}, section 4) to bound $|\mcal{D}_{k,N}|$. Using the stars and bars argument of Feller (\cite{feller}, page 38), we have
\[ |\mcal{D}_{K,N}| =\binom{N+K-1}{N}.\]
It has the usual bound $\binom{N+K-1}{N} \leq (N+1)^K$.
\item Let $w$ be in $\mcal{W}_K$. We write $a_i=\lfloor N w_i\rfloor\geq 1$. There exists $d\in D_{K,N}$ such that
\[d_i=a_i \text{ or } d_i=a_i+1, \text{ for all } i=1,\dots,K.\]
Then, there is $v$ in $\mcal{W}_{K,N}$, defined by $v_i=\frac{d_i}{N}\geq\delta$, such that
\[\forall i, |w_i-v_i|\leq 1/N,\]
i.e. $d(w,v) \leq 1/N \leq \epsilon$.
\end{itemize}
Therefore,
\[ \log(N(\epsilon,\mcal{W}_K,d))\leq K\log(1+N)\leq K\log\left(\frac{3}{\epsilon}\right).\]
\subsection*{Proof of Lemma \ref{lem:aux_bound_log}}
\begin{itemize}
\item First assume $x_1+\dots+x_n=1$, i.e. $x\in\mcal{W}_n$. Then
\begin{align*}
\log\left( \frac{(x_1+\dots+x_n)^{x_1+\dots+x_n}}{x_1^{x_1}\dots x_n^{x_n}} \right) &= - \suml_{i=1}^n x_i\log(x_i),\\
\left( \prod\limits_{i=1}^n x_i \right)^{\frac{1}{x_1+\dots+x_n}} &= \prod\limits_{i=1}^n x_i.
\end{align*}
One can easily check that the function $x\mapsto \suml_{i=1}^n x_i \log(x_i)$ is bounded on $\mcal{W}_n$ and attains a minimum for $x_1=\dots=x_n=1/n$ such that
\[ \log\left( \frac{(x_1+\dots+x_n)^{x_1+\dots+x_n}}{x_1^{x_1}\dots x_n^{x_n}} \right) \leq \log(n).\]
Similarly, one can verify that the function $x\mapsto \prod\limits_{i=1}^n x_i$ is bounded on $\mcal{W}_n$ and attains a maximum for $x_1=\dots=x_n=1/n$ such that
\[ \prod\limits_{i=1}^n x_i \leq \left(\frac{1}{n}\right)^n.\]
\item Now, for any $x_1,\dots,x_n\geq 0$. If $x_1=\dots+x_n=0$ then the result is obvious. Otherwise $s(x):=x_1+\dots+x_n>0$ and we define $y$ in $\mcal{W}_n$ by
$y_i=x_i/s(x)$. Therefore,
\begin{align*}
\log\left( \frac{(x_1+\dots+x_n)^{x_1+\dots+x_n}}{x_1^{x_1}\dots x_n^{x_n}} \right) &= \suml_{i=1}^n x_i\log\left( \frac{s(x)}{x_i} \right)\\
&= - \suml_{i=1}^n x_i \log(y_i)\\
&= s(x) \times \left[ - \suml_{i=1}^n y_i\log(y_i) \right]\\
&\leq (x_1+\dots+x_n)\log(n).
\end{align*}
\begin{align*}
\left( \prod\limits_{i=1}^n x_i \right)^{\frac{1}{x_1+\dots+x_n}} &= s(x)^{1/s(x)} \times \left( \prod\limits_{i=1}^n y_i \right)^{1/s(x)}\\
&\leq  s(x)^{1/s(x)} \times \left( \frac{1}{n^n} \right)^{1/s(x)}\\
&= \left[ \left(\frac{s(x)}{n^n}\right)^{\frac{n^n}{s(x)}} \right]^{\frac{1}{n^n}}.
\end{align*}
Now, we also use that
$\forall x>0, x^{1/x}\leq e^{1/e}$, so finally
\[ \left( \prod\limits_{i=1}^n x_i \right)^{\frac{1}{x_1+\dots+x_n}} \leq \left( e^{\frac{1}{e}}\right)^{^\frac{1}{n^n}}.\]
\end{itemize}
\subsection*{Proof of Lemma \ref{lem:aux_q}}
First, computation gives
\begin{align}
|\phi\circ q_1 -\phi\circ q_2| &= \left| \psi\left(\sqrt{\frac{q_1}{\overline{p}}}(x)\right)-\psi\left(\sqrt{\frac{q_2}{\overline{p}}}(x)\right)\right|\nonumber\\
&= \left| \frac{\sqrt{\frac{q_1}{\overline{p}}}-1}{\sqrt{\frac{q_1}{\overline{p}}}+1} - \frac{\sqrt{\frac{q_2}{\overline{p}}}-1}{\sqrt{\frac{q_2}{\overline{p}}}+1} \right|\nonumber\\
&= \frac{ 2 \left| \sqrt{\frac{q_1}{\overline{p}}} - \sqrt{\frac{q_2}{\overline{p}}} \right| }{\left(\sqrt{\frac{q_1}{\overline{p}}}+1\right) \left(\sqrt{\frac{q_2}{\overline{p}}}+1\right)}\nonumber\\
&= \frac{2 \left| \frac{q_1-q_2}{\overline{p}} \right|}{\left(\sqrt{\frac{q_1}{\overline{p}}}+1\right) \left(\sqrt{\frac{q_2}{\overline{p}}}+1\right) \left( \sqrt{\frac{q_1}{\overline{p}}} + \sqrt{\frac{q_2}{\overline{p}}} \right)}.\label{eq:psi_aux}
\end{align}
For any $x,y_1,y_2\geq 0$,
\begin{align}
\frac{ \sqrt{w} \left(\sqrt{y_1}+1\right)}{ \left(\sqrt{w y_1 + (1-w) x} + 1 \right) } &\times  \frac{ \sqrt{w} \left(\sqrt{y_2}+1\right) }{ \left( \sqrt{w y_2 + (1-w) x} + 1 \right) }\nonumber\\
&\times \frac{ \sqrt{w} \left( \sqrt{y_1} + \sqrt{y_2} \right)}{   \left( \sqrt{w y_1 + (1-w) x} + \sqrt{w y_2 + (1-w) x} \right) } \leq 1.\label{eq:sqrt_delta}
\end{align}
Using (\ref{eq:psi_aux}) and (\ref{eq:sqrt_delta}), we have
\begin{align*}
&|\phi\circ (w q_1 + (1-w) r) - \phi \circ (w q_2 + (1-w) r)|\\
&= \frac{2 w \left| \frac{q_1-q_2}{\overline{p}} \right|}{ \left(\sqrt{\frac{w q_1 + (1-w) r}{\overline{p}}} + 1 \right) \left( \sqrt{\frac{w q_1 + (1-w) r}{\overline{p}}} + 1 \right) \left( \sqrt{\frac{w q_1 + (1-w) r}{\overline{p}}} + \sqrt{\frac{w q_2 + (1-w) r}{\overline{p}}} \right)}\\
&= \frac{2 \left| \frac{q_1-q_2}{\overline{p}} \right| }{ \left(\sqrt{\frac{q_1}{\overline{p}}}+1\right) \left(\sqrt{\frac{q_2}{\overline{p}}}+1\right) \left( \sqrt{\frac{q_1}{\overline{p}}} + \sqrt{\frac{q_2}{\overline{p}}} \right)}\\
&\times \frac{w \left(\sqrt{\frac{q_1}{\overline{p}}}+1\right) \left(\sqrt{\frac{q_2}{\overline{p}}}+1\right) \left( \sqrt{\frac{q_1}{\overline{p}}} + \sqrt{\frac{q_2}{\overline{p}}} \right)}{ \left(\sqrt{\frac{w q_1 + (1-w) r}{\overline{p}}} + 1 \right) \left( \sqrt{\frac{w q_2 + (1-w) r}{\overline{p}}} + 1 \right) \left( \sqrt{\frac{w q_1 + (1-w) r}{\overline{p}}} + \sqrt{\frac{w q_2 + (1-w) r}{\overline{p}}} \right) }\\
&\leq \left| \phi \circ q_1 - \phi \circ q_2 \right| \frac{1}{ \sqrt{w}} \leq \frac{1}{\sqrt{\delta}} |\phi \circ q_1 - \phi \circ q_2 |.
\end{align*}
Then, you just have to take the $L_2(Q)$ norm and it gives the result.
\subsection*{Proof of Lemma \ref{lem:aux_pi_K}}
Let $g_1,\dots,g_K$ be $K$ densities. Let $w,v\in\mcal{W}_{K,\delta}$. If $w=v$ the proof is obvious. From now on, we consider $w\neq v$. We define
\begin{align*}
t_1 &:= \max_i \frac{w_i-v_i}{\mathbbm{1}_{w_i>v_i}-v_i}\in[0,1],\\
t_2 &:= \max_i \frac{v_i-w_i}{\mathbbm{1}_{v_i>w_i}-w_i}\in[0,1].
\end{align*}
Since $w\neq v$, we have $t_1,t_2>0$. We define
\begin{align*}
F_1 &:= \suml_i \left[ v_i + \frac{w_i-v_i}{t_1} \right] g_i\\
F_2 &:= \suml_i \left[ w_i + \frac{v_i-w_i}{t_2} \right] g_i.
\end{align*}
Then, one can check that
\begin{align*}
\suml_i w_i g_i &= \frac{t_1}{t_1+t_2(1-t_1)} F_1 + \frac{t_2(1-t_1)}{t_1+t_2(1-t_1)} F_2,\\
\suml_i \nu_i g_i &= \frac{t_1(1-t_2)}{t_2+t_1(1-t_2)} F_1 + \frac{t_2}{t_2+t_1(1-t_2)} F_2.
\end{align*}
We can now use the following lemma.
\begin{lemme}
\label{lem:aux_pi_2}
For any $r,q$ densities and $w,v\in(0,1)$, we have
\begin{align*}
&|\phi\circ (w r + (1-w) q) - \phi \circ (v r + (1-v) q )|\nonumber\\
&\leq  2 \left( \frac{|(1-w)^{1/4}-(1-v)^{1/4}|}{ (1-w)^{1/4}+(1-v)^{1/4}} \bigvee \frac{|w^{1/4}-v^{1/4}|}{ w^{1/4} + v^{1/4} }\right).
\end{align*}
\end{lemme}
It gives
\begin{align*}
&||\phi \circ (w_1 g_1 + \dots + w_k g_K) - \phi\circ (v_1 g_1 + \dots + v_K g_K)||_{2,Q}\\
&= \left|\left|\phi\circ \left(\frac{t_1}{t_1+t_2(1-t_1)} F_1 + \frac{t_2(1-t_1)}{t_1+t_2(1-t_1)} F_2 \right) \right.\right.\\
&- \left.\left. \phi\circ \left( \frac{t_1(1-t_2)}{t_2+t_1(1-t_2)} F_1 + \frac{t_2}{t_2+t_1(1-t_2)} F_2 \right) \right|\right|_{2,Q}\\
&\leq 2 \left( \frac{ \left| \left(\frac{t_2(1-t_1)}{t_1+t_2(1-t_1)}\right)^{1/4} - \left(\frac{t_2}{t_2+t_1(1-t_2)}\right)^{1/4}\right|}{ \left(\frac{t_2(1-t_1)}{t_1+t_2(1-t_1)}\right)^{1/4} + \left(\frac{t_2}{t_2+t_1(1-t_2)}\right)^{1/4} } \bigvee \frac{ \left| \left(\frac{t_1}{t_1+t_2(1-t_1)}\right)^{1/4} - \left(\frac{t_1(1-t_2)}{t_2+t_1(1-t_2)}\right)^{1/4}\right|}{ \left(\frac{t_1}{t_1+t_2(1-t_1)}\right)^{1/4} + \left(\frac{t_1(1-t_2)}{t_2+t_1(1-t_2)}\right)^{1/4} }\right)\\
&= 2 \left( \frac{ \left| \left(t_2(1-t_1)\right)^{1/4} - \left(t_2\right)^{1/4}\right|}{ \left(t_2(1-t_1)\right)^{1/4} + \left(t_2\right)^{1/4} } \bigvee \frac{ \left| \left(t_1\right)^{1/4} - \left(t_1(1-t_2)\right)^{1/4}\right|}{ \left(t_1\right)^{1/4} + \left(t_1(1-t_2)\right)^{1/4} }\right)\\
&= 2 \left( \frac{ \left| (1-t_1)^{1/4} - 1 \right|}{ (1-t_1)^{1/4} + 1 } \bigvee \frac{ \left| 1 - (1-t_2)^{1/4}\right|}{ 1 + (1-t_2)^{1/4} }\right)\\
&= 2 \left( \frac{ t_1 }{ \left((1-t_1)^{1/4} + 1\right)^2 \left((1-t_1)^{1/2} + 1\right) } \bigvee \frac{ t_2 }{ \left((1-t_2)^{1/4} + 1\right)^2 \left((1-t_2)^{1/2} + 1\right) } \right)\\
&\leq 2(t_1\vee t_2)\\
&\leq \frac{2}{\delta} d(w,v).
\end{align*}
\subsection*{Proof of Lemma \ref{lem:aux_pi_2}}
Using (\ref{eq:psi_aux}), we get
\begin{align*}
&|\phi\circ (w r + (1-w) q) - \phi \circ (v r + (1-v) q )|\\
&= \frac{2 |w-v| \left| \frac{r-q}{\overline{p}} \right| }{ \left(\sqrt{\frac{w r + (1-w) q}{\overline{p}}} + 1 \right) \left( \sqrt{\frac{v r + (1-v) q}{\overline{p}}} + 1 \right) \left( \sqrt{\frac{w r + (1-w) q}{\overline{p}}} + \sqrt{\frac{v r + (1-v) q}{\overline{p}}} \right) }\\
&\leq  \begin{cases}
\frac{2 |w-v| \left| \frac{r-q}{\overline{p}} \right| }{ \left(\sqrt{\frac{w |r-q| +(1-w) q}{\overline{p}}} + 1 \right) \left( \sqrt{\frac{v |r-q| + (1-v) q}{\overline{p}}} + 1 \right) \left( \sqrt{\frac{w |r-q| + (1-w) q}{\overline{p}}} + \sqrt{\frac{v |r-q| + (1-v) q}{\overline{p}}} \right) },\quad\text{if } r\geq q\\
\frac{2 |w-v| \left| \frac{r-q}{\overline{p}} \right| }{ \left(\sqrt{\frac{(1-w) |q-r| + w r}{\overline{p}}} + 1 \right) \left( \sqrt{\frac{(1-v) |q-r| + v r}{\overline{p}}} + 1 \right) \left( \sqrt{\frac{(1-w) |q-r| + w r}{\overline{p}}} + \sqrt{\frac{(1-v) |q-r| + v r}{\overline{p}}} \right) },\quad\text{if r<q.} 
     \end{cases} 
\end{align*}
One can easily check that the following function
\[x \mapsto \frac{x}{ \left(\sqrt{\alpha x } + 1 \right) \left( \sqrt{\beta x} + 1 \right) \left( \sqrt{\alpha x} + \sqrt{\beta x} \right) }\]
is bounded above by
\[\frac{1}{ \left( \alpha^{1/4} + \beta^{1/4} \right)^2 \left( \alpha^{1/2} + \beta^{1/2} \right) }. \]
Then we get
\begin{align*}
&|\phi\circ (w r + (1-w) q) - \phi \circ (v r + (1-v) q )|\\
&\leq  2 |w-v| \times\\
&\left( \frac{1}{\left( (1-w)^{1/4}+(1-v)^{1/4} \right)^2 \left( (1-w)^{1/2} + (1-v)^{1/2} \right)} \bigvee \frac{1}{\left( w^{1/4} + v^{1/4} \right)^2 \left( w^{1/2} + v^{1/2} \right)}\right)\\
&= 2 \left( \frac{|\sqrt{1-w}-\sqrt{1-v}|}{\left( (1-w)^{1/4}+(1-v)^{1/4} \right)^2} \bigvee \frac{|\sqrt{w}-\sqrt{v}|}{\left( w^{1/4} + v^{1/4} \right)^2 }\right)\\
&= 2 \left( \frac{|(1-w)^{1/4}-(1-v)^{1/4}|}{ (1-w)^{1/4}+(1-v)^{1/4}} \bigvee \frac{|w^{1/4}-v^{1/4}|}{ w^{1/4} + v^{1/4} }\right).
\end{align*}
\end{document}